%% file: Cellulaire.tex
\documentclass[10pt]{smfart}
\usepackage{a4wide}
\usepackage{amsmath,amssymb,smfthm}
\usepackage[francais]{babel}
\usepackage[latin1]{inputenc}
\usepackage[all]{xy}
\usepackage{graphicx}
\usepackage{epic,eepic}
\usepackage{epsfig}
\usepackage{color}

\def\R{\mathbb{R}}

\def\C{\mathbb{C}}

\def\Cp{\mathbb{C}_p}
\def\Fp{\mathbb{F}_p}
\def\Fq{\mathbb{F}_q}
\def\Fqb{\overline{\mathbb{F}}_q}
\def\Fpb{\overline{\mathbb{F}}_p}
\def\L{\mathcal{L}}
\def\La{\Lambda}

\def\O{\mathcal{O}}
\def\F{\mathcal{F}}

\def\*{^\times }
\def\dpt{\displaystyle}

\def\Gm{\mathbb{G}_m}         
\def\l{\lambda}
\def\g{\gamma}

\def\a{\alpha}
\def\b{\beta}

\def\s{\sigma}
\def\ph{\varphi}

\def\e{\epsilon}
\def\ssi{\Leftrightarrow}

\def\impl{\Rightarrow}
\def\limpl{\Longrightarrow}
\def\drt{\rightarrow}
\def\ldrt{\longrightarrow}
\def\Q{\mathbb{Q}}

\def\Qp{\mathbb{Q}_p}
\def\Qpb{\overline{\mathbb{Q}}_p}
\def\Zp{\mathbb{Z}_p}
\def\Z{\mathbb{Z}}
\def\N{\mathbb{N}}
\def\Zl{\mathbb{Z}_\ell}
\def\Hom{\text{Hom}}

\def\={\! = \!}

\def\spec{\text{Spec}}
\def\spf{\text{Spf}}

\def\E{\mathcal{E}}

\def\limp{\underset{\longleftarrow}{\text{ lim }}\;}
\def\limi{\underset{\longrightarrow}{\text{ lim }}\;}

\def\iso{\xrightarrow{\;\sim\;}}

\def\End{\text{End}}

\def\GL{\hbox{GL}}

\def\xrig{\xrightarrow}

\def\M{\mathcal{M}}
\def\Mf{\widehat{\mathcal{M}}}

\def\X{\mathfrak{X}}
\def\GG{\Gamma}
\def\Ext{\text{Ext}}
\def\bc{\backslash}
\def\AA{\mathcal{A}}

\def\<{<\hspace{-1mm}}
\def\>{\hspace{-1mm}>}

\def\Lie{\text{Lie}}

\def\dem{{\it D{\'e}monstration. }} 
\def\xs{\underline{x}}

\def\et{\text{ét}}

\author{  Laurent Fargues} 
\address{CNRS-université Paris-Sud-IHES} 
\email{laurent.fargues@math.u-psud.fr} 
\date{}

\begin{document}

\newtheorem{Fait}{Fait}

\title{L'isomorphisme entre les tours de Lubin-Tate et de Drinfeld  : 
Décomposition cellulaire de la tour de Lubin-Tate}
\maketitle

\begin{abstract}
Cet article est le premier d'une série visant à construire un
isomorphisme entre les tours $p$-adiques de Lubin-Tate et de Drinfeld,
décrire cet isomorphisme et en donner des applications. Nous-y 
construisons un modèle entier $p$-adique équivariant en niveau infini de la
tour de Lubin-Tate. Ce schéma formel $p$-adique sera comparé plus tard
à un autre associé à la tour de Drinfeld. 
\end{abstract}

\begin{altabstract}
This article is the first one of a series aiming to construct an
isomorphism between the $p$-adic Lubin-Tate and Drinfeld towers,
describe this isomorphism and give applications.  We construct  a
$p$-adic equivariant integral model of the Lubin-Tate tower with
infinite level. This formal scheme will be later  compared to another
one associated to the Drinfeld tower.
\end{altabstract}

\section*{Introduction générale à la série d'articles sur les tours de
Lubin-Tate et de Drinfeld}

Cet article est le premier d'une série visant à construire un
isomorphisme entre les tours $p$-adiques de Lubin-Tate et de Drinfeld,
décrire cet isomorphisme et en donner des applications. 
\\

Voici une bref description des différents articles :
\begin{itemize}
\item Le présent article pour lequel on renvoit à l'introduction qui suit.
\item L'article \cite{Points} intitulé ``L'isomorphisme entre les tours de
  Lubin-Tate et de Drinfeld au niveau des points'' dans lequel on
  construit l'isomorphisme entre les deux tours au niveau des points à
  valeurs dans des corps valués du type de ceux intervenant dans la
  théorie des espaces de Berkovich.
\item L'article \cite{Rami} intitulé ``Application de Hodge-Tate duale d'un groupe
  de Lubin-Tate, immeuble de Bruhat-Tits du groupe linéaire et
  filtrations de ramification''. On y décrit entre autres  précisément
  l'isomorphisme entre les deux tours au niveau des squelettes des
  espaces de Lubin-Tate et de Drinfeld.
\item L'article \cite{iso4} intitulé ``L'isomorphisme entre les tours
  de Lubin-Tate et de Drinfeld : démonstration du résultat
  principal''. On y démontre l'existence de l'isomorphisme entre les
  éclatés du schéma formel construit dans le présent article et un
  autre schéma formel construit à partir de la tour de Drinfeld.
\item L'article \cite{iso5} intitutlé ``L'isomorphisme entre les tours
  de Lubin-Tate et de Drinfeld : comparaison de la cohomologie des
  deux tours''. On y démontre que l'isomorphisme construit dans
  \cite{iso4} induit un isomorphisme entre les complexes de
  cohomologie à support compact équivariants des tours de Lubin-Tate
  et de Drinfeld. On y démontre également l'existence d'une
  correspondances de Jacquet-Langlands géométrique entre faisceaux
  équivariants sur l'espace des périodes de Gross-Hopkins et faisceaux
  équivariants sur l'espace de Drinfeld. 
\end{itemize}
\vspace{4mm}

Pour le présent article, les articles \cite{Points} et \cite{iso4} nous
suivons le plan fourni par \cite{Faltings8} et il apparaîtra comme clair
au lecteur que nous nous sommes largement inspirés des travaux de Gerd Faltings.

 Les articles \cite{Rami} et
\cite{iso5} eux sont complétement nouveaux et indépendants de \cite{Faltings8}.

\section*{Introduction à cet article} 

Cette première partie vise à $p$-adifier la tour de Lubin-Tate : les
modèles entiers usuels de la tour de Lubin-Tate (\cite{Drinfeld1}) sont donnés par le
spectre formels d'anneaux du type $\Zp [[x_1,\dots,x_{n-1}]]$, un
idéal de définition étant $(p,x_1,\dots,x_{n-1})$. Les modèles entiers
naturels 
de la tour de Drinfeld eux sont $p$-adiques : l'idéal de définition
des schémas formels associés est l'idéal engendré par $p$, les anneaux
associés étant du type $\Zp <x_1,\dots,x_d>/\text{Idéal}$.  
Pour pouvoir comparer ces deux tours nous devons donc d'abord modifier
les modèles  usuels de la tour de Lubin-Tate. 

Utilisant
certains domaines fondamentaux quasicompacts pour l'action des correspondances de
Hecke sphériques on reconstruit un modèle entier $p$-adique en niveau
infini de la tour de Lubin-Tate, modèle qui pourra être comparé à la
tour de Drinfeld. Ce schéma formel $p$-adique est obtenu par recollement d'itérés sous
les correspondances de Hecke d'un modèle entier $p$-adique du domaine
fondamental.  

\begin{figure}[h]
   \begin{center}
      \input{DrLT1_immeuble2.pstex_t}
   \end{center}
\caption{\footnotesize La décomposition cellulaire de l'arbre de
  $PGL_2$ indexée par les sommets}
\label{arbreDrLT1}
\end{figure}
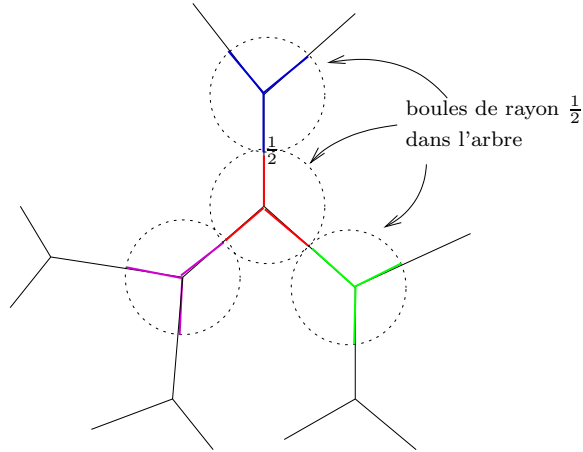

Ces correspondances de Hecke sont paramétrées par un immeuble de
Bruhat-Tits et le schéma formel est obtenu par recollement équivariant de
``cellules'', des schémas formels affines, indexées par les sommets de
 l'immeuble. Il est à noter que, contrairement à celle de l'espace de Drinfeld, 
cette décomposition cellulaire est indexée par les sommets de
l'immeuble  et non les simplexes. Ainsi on verra par exemple  dans les
articles \cite{Rami} et \cite{iso4} que dans le
cas de $GL_2$ 
l'image de ces cellules dans l'immeuble de Bruhat-Tits paramétrant l'espace de
Drinfeld consiste en les boules de rayon $1/2$ centrées
en les sommets de l'immeuble (figure \ref{arbreDrLT1}).

Le schéma formel $p$-adique construit n'est pas un schéma formel
topologiquement de type fini, il vit en niveau infini. Il n'est pas
possible d'effectuer sa construction en niveau fini (par niveau
on entend un sous-groupe compact ouvert de $\GL_n (\Zp)$). Plus
précisément si $r\in [0 1[\cap \Q$ et $\mathbb{B}^{n-1} (0,r)$ désigne
la boule fermée de rayon $r$ dans $\mathring{\mathbb{B}}^{n-1} (0,1)$,
la base de l'espace de Lubin-Tate (sans niveau), on peut alors 
construire un tel espace en niveau fini pour l'image réciproque de
cette boule dans la tour de Lubin-Tate. Mais lorsque $r\drt 1$, ou si
l'on veut lorsque l'on tend vers le bord de l'immeuble de Bruhat-Tits
i.e. on sort de tout compact, le niveau doit tendre vers l'infini. Cela est relié au fait suivant. Soit $\mathcal{I}$ l'immeuble de Bruhat-Tits de $\text{PGL}_n (\Qp)$ et $A\subset \mathcal{I}$ un sous-ensemble simpliciale fini. Il existe alors un sous-groupe compact ouvert $K_A\subset \GL_n (\Zp)$ tel que le composé $A\hookrightarrow \mathcal{I} \twoheadrightarrow \mathcal{I}/K_A$ soit un plongement. Mais lorsque $A$ grandit $K_A\ldrt \{Id \}$.
\\

Décrivons les différentes parties de l'article : 

\begin{itemize}
\item Dans les sections 2,3 et 4 nous démontrons des résultats
  généralisant ceux de \cite{HopkinsGross} concernant l'application
  des périodes de l'espace de Lubin-Tate vers l'espace projectif et
  les domaines fondamentaux pour les correspondances de Hecke
  non-ramifiées dans ces espaces. Nous y utilisons pour cela la théorie des
  Displays qui donne des formules matricielles  pour
  l'application des périodes.
  Certains de ces résultats sont déjà dans l'article \cite{Yu} qui utilise le point 
de vue de \cite{HopkinsGross} des quasi-logarithmes, nous préférons utiliser la
théorie des Displays. 
 Puis nous étudions en détails les
  domaines fondamentaux  pour l'action des correspondances de Hecke
  non-ramifiées sur la base de la tour de Lubin-Tate définis par Gross-Hopkins dans \cite{HopkinsGross}.
  Ces résultats sont placés dans un cadre plus conceptuel et
  généralisés dans \cite{Rami}. Néanmoins les résultats plus complexes
  de \cite{Rami}
  ne sont pas nécessaires à une première compréhension de
  l'isomorphisme entre les deux tours. C'est pourquoi ils ne sont pas
  inclus dans cet article. Par exemple dans \cite{Rami} on montre comment
 construire de façon systématique des domaines fondamentaux
 généralisant ceux de
 Gross-Hopkins et comment comprendre la façon dont ils se recollent 
avec leurs itérés sous des correspondances de Hecke.
\item Les sections 5 à 8 sont le coeur de l'article. Nous y utilisons
  les domaines fondamentaux de Gross-Hopkins pour décomposer
  cellulairement la tour de Lubin-Tate. Nous y construisons à la fin un
  schéma formel $p$-adique cellulairement décomposé au dessus de
   l'immeuble de Bruhat-Tits du groupe linéaire. La cellule au dessus
   d'un sommet de l'immeuble  est le spectre formel d'une algèbre
   $p$-adique du type la boule unité dans une algèbre de Banach
   $p$-adique obtenue par complétion d'une union croissante
   $\bigcup_{k\geq 0} \mathcal{A}_k$ où  $\mathcal{A}_k$ est une
   algèbre affinoïde munie de sa norme infinie et $\mathcal{A}_k \drt
   \mathcal{A}_{k+1}$ est fini. L'algèbre $\mathcal{A}_0$ est la fibre
   générique d'un
   espace de déformations de groupe $p$-divisible avec contraintes sur
   le polygone de Newton de ses points de $p$-torsions. Les algèbres
   $\mathcal{A}_k$ sont obtenues en ajoutant des points de torsions
   (structures de niveau).
\item Dans l'appendice A on discute de la normalisation des schémas
  formels $p$-adiques admissibles dans un revêtement fini de leur fibre générique
  vue comme espace rigide. On utilise ces résultats de façon cruciale
  dans la construction de nos schémas formels. 
Cet appendice est né de diverses questions
  de géométrie rigide que s'est posées l'auteur et pour lesquelles il n'a
  pas trouvé de références dans la littérature.   
\item Soit $\O$ l'anneau des entiers d'une extension de degré fini de
  $\Qp$ d'uniformisante $\pi$. 
Dans l'appendice B nous établissons une théorie de la
  déformation pour les $\O$-modules  $\pi$-divisibles
  relativement aux immersions nilpotentes définies par des idéaux
  munis de $\pi$-analogues des puissances divisées (le cas $\O=\Zp$ et
  $\pi=p$
  étant le cas ``classique'' traité par la théorie de Messing). Cette
  section n'est pas strictement nécessaire pour démontrer les
  résultats auxquels nous nous intéressons mais permet néanmoins
  d'interpréter agréablement certains objets intervenant dans la
  définition de l'application des périodes, et accessoirement
  d'améliorer les résultats d'intégralité de cette application.  
  Les résultats de \cite{Faltings7} sont plus généraux que ceux de cet
  appendice puisqu'ils s'appliquent aux groupes plats finis. Néanmoins
  dans le cas des groupes $p$-divisibles la méthode utilisée dans
  l'appendice B est plus simple que celle de \cite{Faltings7}. 
\end{itemize}

\vspace{8mm}

{\it Prérequis : Concernant les espaces de Lubin-Tate on supposera le
  lecteur familier avec les chapitres 1 et 4 de \cite{Drinfeld1}, le
  chapitre II de \cite{Har4} et \cite{HopkinsGross}. 
Concernant \cite{HopkinsGross} nous n'aurons à utiliser que le
corollaire 23.26 de cet article, les autres résultats étant
redémontrés et généralisés par d'autres méthodes (la théorie des
Displays). Néanmoins la lecture de \cite{HopkinsGross} est fortement
conseillée car elle fournit une
introduction ``concrète'' à l'étude détaillée des espaces de Lubin-Tate.
Nous utiliserons le langage des espaces
  de Rapoport-Zink (\cite{RZ}) qui est le langage naturel de ce type de problème. On pourra
  consulter \cite{Har5} pour une reformulation de \cite{Har4} dans ce
  langage là. On suppose également le lecteur familier avec la théorie
  de la déformation de Messing (\cite{Messing1}).  
\\

Avertissements : Le lecteur uniquement intéressé par la construction du
schéma formels $p$-adique qui sera relié plus tard à l'espace de
Drinfeld peut sauter les chapitres 2 et 4. 
Néanmoins le résultat concernant l'application des périodes sera
utilisé plus tard dans la construction de l'isomorphisme entre les
deux tours.
L'appendice B  est réservé
aux ``experts''.
}
\\

{\it Enfin il apparaîtra comme clair au lecteur que l'auteur de cet
  article s'est largement inspiré des travaux de Gerd Faltings
  \cite{Faltings8}.} 

\section*{Introduction bis : Construction du schéma formel dans le cas
  de $\GL_2 (\Qp)$ par
  éclatements des courbes modulaires}

Soit pour $N\geq 3$ $Y(N)$ la courbe modulaire ouverte sur $\spec
(\Z)$ classifiant les courbes elliptiques $E$ munies d'une structure
de niveau $N$, 
$N^{-1}\Z/\Z\ldrt E[N]$, au sens de Katz-Mazur. On a 
$$
Y(N) (\C) = \coprod_{(\Z/N\Z)^\times} \GG (N) \bc \mathbb{H}^{\pm}
$$

Fixons un nombre premier $p$ ainsi qu'un
entier $N_0\geq 3$ tel que $(N_0,p)=1$. Soit $E$ la courbe elliptique
universelle sur $Y_0 (N)\otimes \Fp$. Considérons le Verschiebung $V:
E^{(p)}\ldrt E$ c'est à dire l'isogénie duale du Frobenius $F:E\ldrt
E^{(p)}$. Le morphisme induit au niveau des formes différentielles
invariantes est
$$
V^* :\omega_E\ldrt \omega_{E^{(p)}} \simeq \omega_E^{\otimes p}
$$
où l'on a fixé un isomorphisme (non canonique) entre fibrés en droites
 $ \omega_{E^{(p)}} \simeq \omega_E^{\otimes p}$. Il fournit une forme
 modulaire mod $p$ de poids $k-1$ 
$$
H\in \GG ( Y(N_0)_{\Fp},\omega_E^{\otimes (p-1)})
$$
qui est l'invariant de Hasse. Cela définit un diviseur de Cartier
réduit $D$ dans $Y(N_0)_{\Fp}$ de support le lieu supersingulier, un
nombre fini de points de la fibre spéciale. 
Si $y\in Y (N_0) (\Fpb)_{s.s.}$ 
$$
Y (N_0)^{\widehat{\;}}_{/\{ y\}} \simeq \spf ( W (\Fpb) [[x ]])
$$
où l'on peut choisir $x \equiv H \text{ mod p}$ (si $p\neq 2,3$ $\; x$
peut être choisi égal à une série d'Eisenstein $E_{p-1}$). 

Désormais on notera $Y(N_0)$ pour $Y (N_0)\otimes_\Z
\widehat{\Zp^{nr}}$ et $D$ pour le diviseur de Cartier précédent
étendu à 
$Y(N_0)_{\Fpb}$. Soit 
$$
Y(p^\infty N_0) = \underset{k\geq 0}{\limp} Y(p^k N_0)
$$
un $\widehat{\Zp^{nr}}$-schéma et
$$
\pi : Y(p^\infty N_0)\ldrt Y(N_0)
$$
la projection. C'est un morphisme plat
totalement ramifié au dessus des points supersinguliers. 
 Le schéma $Y(p^\infty N_0)$
est muni d'une action de Hecke de $\GL_2 (\Qp)$. Le morphisme $\pi$
est $\GL_2 (\Zp)$-invariant. 

Considérons le diviseur de Cartier 
$$
D_\infty = \pi^* D
$$
Il est ``infiniment ramifié'' au sens où le sous-schéma fermé de
$Y(p^\infty N_0)_{\Fpb}$ définit par $D_\infty$ possède des nilpotents
d'ordre quelconque. Néanmoins $\text{supp} ( \pi^* D)$ est un nombre
fini de points fermés 
$$
\pi : \text{supp} (\pi^* D) \iso \text{supp} (D) = Y
(N_0)_{s.s.} (\Fpb)
$$
Le morphisme $\pi$ étant $\GL_2 (\Zp)$-invariant $D_\infty$ est $\GL_2
(\Zp)$-invariant 
$$
\forall g\in \GL_2 (\Zp)\;\; g^* D_\infty =D_\infty
$$
Considérons les itérés de $D_\infty$ par les correspondances de Hecke 
$$
\{ g^* D_\infty \;|\; g\in \GL_2 (\Qp)/\GL_2 (\Zp) \}
$$
ensemble paramétré par l'arbre de $\GL_2$. Bien sûr $\forall g\in
\GL_2 (\Qp)/\GL_2 (\Zp) \; \text{supp} ( g^* D_\infty ) =\text{supp}
(D_\infty)$ puisque les correspondances de Hecke laissent invariantes
le lieu supersingulier. Néanmoins lorsque $g$ varie dans $\GL_2
(\Qp)/\GL_2 (\Zp)$ ces diviseurs de Cartier sont distincts. 

Fixons un entier $\a\geq 1$. Pour $A\subset GL_2 (\Qp)/\GL_2 (\Zp)$ un
sous-ensemble fini de sommets de l'arbre soit
$$
X_A = \text{L'éclatement de } Y(p^\infty N_0)\text{ le long du
  diviseur } \prod_{g\in A} g^* D_\infty^\a \text{ de la fibre
  spéciale } Y(p^\infty N_0)_{\Fpb}
$$
(un diviseur de la fibre spéciale définit un sous-schéma fermé de
codimension $2$ du modèle entier) et $U_A\subset X_A$ l'ouvert où $p$
engendre le diviseur exceptionnel. En fibre générique, i.e. après
inversion de $p$, $(U_A)_\eta = (X_A)_\eta = Y(p^\infty N_0)_\eta$. 

Si $A\subset B$ notons $\Pi_{A,B} : X_B \ldrt X_A$. On a
$\Pi_{A,B}^{-1} (U_A)\subset U_B$ et $\Pi_{A,B}$ induit un
isomorphisme
$$
\Pi_{A,B} : \Pi_{A,B}^{-1} (U_A) \iso U_A
$$
et donc via $\Pi_{A,B}$, $U_A\subset U_B$. Soit
$$
Z= \underset{A\subset GL_2 (\Qp)/GL_2 (\Zp)\atop fini}{\limi} U_A
$$
où les sous-ensembles finis de l'arbre sont ordonnés par
l'inclusion. C'est un schéma sur $\spec (W(\Fpb))$ muni d'une action
de $\GL_2 (\Qp)$ puisque pour $g\in \GL_2 (\Qp)$ on a $g: U_A \ldrt
U_{g.A}$. De plus $Z_\eta = Y (p^\infty N_0)_\eta$. 

Soit maintenant $Z'$ le normalisé de $Z$ dans sa fibre générique
$Z_\eta$. Soit $\mathfrak{Z}$ le schéma formel sur $\spf ( W(\Fpb))$ 
égal au complété $p$-adique de $Z'$. Il est muni d'une action de
$\GL_2 (\Qp)$. 

Pour $\a=2$
le schéma formel $\mathfrak{Z}$ est ``à quelques détails près'' le schéma formel
$p$-adique  associé à la tour de
Lubin-Tate que nous allons construire en général. Nous ne le construirons pas de cette manière. 
 Pour $\GL_2 (\Qp)$ ce  schéma formel $\X_\infty$ sera un schéma formel
$p$-adique muni d'une action de $\GL_2 (\Qp)\times D^\times$ où
$D|\Qp$ est une algèbre de quaternions. On aura alors 
$$
\mathfrak{Z} \simeq \coprod_{i\in I} \X_{\infty}/ \pi^{a_i\Z}
$$
où $I$ est un ensemble fini et 
 $a_i\in \N\setminus \{0 \}$. Bien sûr sur la construction précédente
 on ne voit pas l'action de $D^\times$, il faut utiliser le théorème
 de Serre-Tate pour la voir.

\section{Hypothèses et notations }

Soit $F|\Qp$ une extension de degré fini, d'uniformisante $\pi$ et
de corps résiduel $k=\mathbb{F}_q = \O_F/\pi \O_F$. 
 On note
$\breve{F}= \widehat{F^{nr}}$ 
le complété de l'extension maximale non-ramifiée de $F$
dans une clôture algébrique de celui-ci et 
$F^0$ l'extension maximale non-ramifiée de $\Qp$ dans $F$.
 On note parfois $\O$ pour $\O_F$ et $\breve{\O}$ pour $\O_{\breve{F}}$.
On fixe un isomorphisme entre le corps résiduel de $\O_{\breve{F}}$ et
$\Fqb$ une clôture algébrique de $\Fq$. 
\\

Si $\mathfrak{Z}$ est un schéma formel d'idéal de définition
$\mathfrak{I}$ par définition la catégorie des groupes $p$-divisibles
sur $\mathfrak{Z}$ est la catégorie limite projective 
$$
\underset{k\geq 1}{\limp} \left ( \text{Groupes }p\text{-divisibles sur }\spec (
\O_\mathfrak{Z}/\mathfrak{I}^k)\right )
$$
de la catégorie fibrée des groupes $p$-divisibles sur le système de schémas ($\mathfrak{Z}
\text{ mod } \mathfrak{I}^k)_{k\geq 1}$. Cela signifie concrètement
que se donner un groupe $p$-divisible sur $\mathfrak{Z}$ est
équivalent à se donner une famille de groupes $p$-divisibles
$(G_k)_{k\geq 1}$ sur les $\spec (
\O_\mathfrak{Z}/\mathfrak{I}^k)_{k\geq 1}$ munis d'isomorphismes 
$\forall k\;G_{k+1} \text{ mod }\mathfrak{I}^k \iso G_k$
satisfaisant une condition de cocyle évidente. Rappelons que si
$\mathfrak{Z} =\spf (A)$ est affine il y a  une équivalence de
catégories entre groupes $p$-divisibles sur $\mathfrak{Z}$ et groupes
$p$-divisibles sur $\spec
(A)$ (il s'agit de l'analogue du théorème d'algébrisation de
Grothendieck). 

Une quasi-isogénie $\ph:G_1\ldrt
G_2$  entre deux groupes $p$-divisibles sur $\mathfrak{Z}$ est un
système compatible de quasi-isogénies sur les  $\spec (
\O_\mathfrak{Z}/\mathfrak{I}^k)_{k\geq 1}$. On prendra garde que la notion
de quasi-isogénies sur  $\spec (A)$ 
est beaucoup plus forte que celle sur $\spf (A)$. 
\\

Si $S$ est un $\O_F$-schéma ou bien un $\spf (\O_F)$-schéma formel un
$\O$-module  $\pi$-divisible sur $S$ est un groupe $p$-divisible
$H$ sur $S$ muni d'une action de $\O$ induisant l'action canonique sur
son algèbre de Lie (cf. appendice B pour plus de détails). Il sera dit
formel si ses fibres géométriques en tous les points de la base ne
possèdent pas de partie étale. 

\subsection{Espaces}

\begin{defi}
Nous notons $\mathbb{H}_0$ un $\O$-module $\pi$-divisible formel de dimension $1$ et 
hauteur $n$ sur $\mathbb{F}_q$. 
 Nous notons $\mathbb{H} = \mathbb{H}_{0/\Fqb}$.
\end{defi}

\begin{defi}
Nous notons $\X_0$ l'espace de Lubin-Tate 
des déformations par isomorphismes de $\mathbb{H}_0$, ou encore
l'espace des {$*$}-déformations d'une loi de groupe formel associée
comme dans \cite{HopkinsGross}.
 Nous notons $H_0$ la déformation universelle. On note de même $\X$, $H$  les objets étendus à $\breve{\O}$. 
\end{defi}

Le  $\O_F$-schéma formel $\X_0$ est non-canoniquement isomorphe à
$\spf (\O [[x_1,\dots,x_{n-1}]])$. Il représente le foncteur qui à une
$\O$-algèbre locale complété $R$ d'idéal maximal $\mathfrak{m}$ et de
corps résiduel $\Fq$ associe les classes d'isomorphisme de couples
$(H_0,\rho_0)$ où $H_0$ est un $\O$-module $\pi$-divisible sur $R$ et 
$$
\rho_0 : \mathbb{H}_0 \iso H_0\otimes_R R/\mathfrak{m}
$$

\begin{defi}
On note $\Mf$ l'espace de Rapoport-Zink associé sur $\spf (\breve{\O})$ des déformations par quasi-isogénies de $\mathbb{H}$.
On note $\M$ l'espace rigide fibre générique associé. 
\end{defi}

Le schéma formel $\widehat{\mathcal{M}}$ représente le foncteur qui à un $\spf
(\breve{\O})$-schéma formel $\mathcal{S}$ associe les classes
d'isomorphismes de couples $(H,\rho)$ où $H$ est un $\O$-module
$\pi$-divisible sur $\mathcal{S}$ et 
$$
\rho : \mathbb{H}\times_{\spec ( \Fqb)} (\mathcal{S}\text{ mod } \pi)
\ldrt H\times_{\mathcal{S}} ( \mathcal{S}\text{ mod } \pi)
$$
est une quasi-isogénie.

\subsection{Action}

Le groupe $\mathbb{H}_0$ étant défini sur $\mathbb{F}_q$, $\mathbb{H}= \mathbb{H}^{(q)}$ et donc le morphisme de Frobenius définit
un élément 
$$
\Pi = \text{Frob}_q\in \End ( \mathbb{H})
$$
et alors 
$$
 \End ( \mathbb{H}) =  \End ( \mathbb{H}_{0/\mathbb{F}_{q^n}})
=\O_{F_n} [\Pi] = \O_D
$$
l'ordre maximal dans l'algèbre à division $D$ d'invariant $\dpt{\frac{1}{n}}$ où $F_n | F$ désigne l'extension non-ramifiée de degré $n$. 
L'action de $F_n$ n'est pas définie sur $\Fq$ mais $\mathbb{F}_{q^n}$ et on a 
$\End ( \mathbb{H}_0)= \O_F[\Pi]$ l'anneau des entiers de 
l'extension totalement ramifiée de $F$  définie par le polynôme d'Eisenstein $X^n - \pi$. 

\begin{defi}
On muni $\X$ de l'action {\it à gauche} 
de $\O_D^\times$  et $\Mf$ de celle de $D^\times$ en posant 
$$
d.(H,\rho)=(H,\rho\circ d^{-1})
$$
pour $H$ un groupe $p$-divisible et $\rho$ l'isomorphisme ou quasi-isogénie définissant la déformation. 
\end{defi}

\subsection{Scindage de l'espace de Rapoport-Zink}

Rappelons que toute quasi-isogénie de degré zéro entre $\O$-modules $\pi$-divisibles
formels de dimension $1$ définis sur un schéma réduit annulé par $p$ est un
isomorphisme.  
Il y a donc une décomposition
$$
\Mf=\coprod_{i\in \Z} \Mf^{[i]}
$$
où $\Mf^{[i]}$ désigne l'ouvert-fermé de $\Mf$ où la quasi-isogénie universelle est de hauteur $i$. De plus,
\begin{eqnarray*}
\Pi^{-i} : \Mf^{[0]} & \iso & \Mf^{[i]} \\
(H,\rho) & \mapsto & (H,\rho\circ \text{Frob}_q^i)
\end{eqnarray*}
et 
$$
\Mf^{[0]} \simeq \X
$$

\subsection{Donnée de descente de Rapoport-Zink}

Elle est définie dans la section 3.48 de \cite{RZ}.
Soit $\s$ le Frobenius arithmétique de $F^{nr}|F$. Cette donnée de
descente $\a$ est donnée dans notre cas par le diagramme suivant 
$$
\xymatrix@C=14mm{
\Mf \ar[r]^\a \ar[d]^\sim  & \Mf^{(\s)} \ar[d]^\sim  \\ 
\dpt{\coprod_{i\in \Z}} \X \ar[r]^{\Pi^{-1} \otimes \s} & \dpt{\coprod_{i\in \Z}} \X
}
$$
où l'opérateur $\Pi^{-1}$ est l'identité de $\X$ décalée de $+1$ dans
$\dpt{\coprod_{i\in \Z}}\X $ et $1\otimes \s$ est donné par
l'isomorphisme $\X=\X_0\hat{\otimes}_{\O} \breve{\O}$. 
On voit en particulier qu'elle n'est pas effective puisqu'elle décale 
de $1$ la hauteur de l'isogénie. Cependant elle devient effective 
sur les quotients $\Mf/\pi^{a \Z}$ pour tout entier $a\in \N^*$
(ce qui revient du point de vue cohomologique à prendre les représentations ayant un caractère central  d'ordre fini).
L'action de $D^\times$ commute à cette donnée de descente $\a$.

\begin{rema}
La donnée de descente $1\otimes \s$  définissant $\X_0$ sur
$\X=\X_0\hat{\otimes} \breve{\O}$  n'est pas la bonne puisque par
exemple elle ne commute pas à l'action de $\O_D^\times$. Dans
\cite{HopkinsGross} Gross et Hopkins  utilisent cependant cette donnée
de descente quitte à tordre l'action de $\O_D^\times$ en le remplaçant
par le groupe étale localement constant $\text{Aut}
(\mathbb{H}_0)$ ($\text{Aut} (\mathbb{H}_0)$ est un groupe sur
$\spec (\Fq)_{\et}$ qui devient constant sur $\spec (\mathbb{F}_{q^n})$, il
définit dont un groupe sur $(\X_0)_{\et}$ qui devient constant sur
$\X_0\otimes_{\O_F} \O_{F_n}$, $F_n|F$ désignant l'extension
non-ramifiée de degré $n$). 
 Nous préférons cependant adopter le point de vue de
Rapoport-Zink.
\end{rema}

\subsection{Polygone de Newton des points de torsion}

Soit $L|F$ une extension valuée complète pour une valuation à valeurs
dans $\R\cup \{+\infty \}$ (il s'agit des corps intervenant dans la
théorie des espaces analytiques de Berkovich). Si $H$ est un
$\O$-module  $\pi$-divisible formel sur $\O_L$ de groupe formel
associé $\widehat{H}$ il y a une application ``valuation''
$$
v: \widehat{H} ( \O_L) \ldrt \R\cup \{+\infty \}
$$
qui est définie en fixant un isomorphisme de $\O_L$-schémas formels
pointés $\spf ( \O_L [[T]]) \iso \widehat{H}$, l'un étant pointé par
la section $T=0$ et l'autre par sa section unité (c'est à dire en fixant
une loi de groupe formel associée),  par la ``valuation de la coordonnée
$T$''. Cette application ne dépend pas du choix d'un tel isomorphisme puisque tout
automorphisme  du schéma formel $(\mathbb{A}^1)^{\widehat{\;}}_{/\{ 0\}}$ 
envoyant l'origine sur l'origine induit en fibre générique un automorphisme de
$\mathring{\mathbb{B}}^{1}$ conservant la distance à l'origine. 

On s'intéressera en fait à la valuation des points de torsion
$H (\O_L)= \widehat{H}[\pi^\infty] (\O_L)$. La ``valuation'' définit
une filtration appelée filtration de ramification inférieure 
par des
sous-$\O$-modules sur ces points de torsion, formée des sous-modules où la
valuation
est supérieure à un nombre donné. 
 Elle est
étudiée en détails dans \cite{Rami}. Bornons nous à quelques rappels. 

Il existe un système de coordonnées formelles $(x_1,\dots,x_{n-1})$ sur $\X_0$ i.e. un isomorphisme 
$$
\X_0\simeq \spf (\O [[x_1,\dots,x_{n-1}]])
$$ 
et une loi de groupe formel universelle $F^{univ}$ telle que 
$$
[\pi]_{F^{univ}} = \pi u_0 T+ x_1 u_1 T^q+\dots + x_{n-1}
u_{n-1} T^{q^{n-1}} 
+ u_n T^{q^n}
$$
où $u_0,\dots ,u_n \in \breve{O}[[x_1,\dots,x_{n-1}]][[T]]^\times$
sont des unités (cf. \cite{Lubin1} chapitre 1 page 106, cela se déduit
également de la théorie de Cartier cf. \cite{HopkinsGross}). 
Nous fixons un tel isomorphisme. Pour $\xs= (x_1,\dots,x_{n-1})\in
\X^{rig}_0 (\overline{F})$, ou plus généralement $\X_0^{an} (L) :=
\X_0 (\O_L)$ avec $L|F$ une extension valuée complète comme
précédemment, 
on note  $H_{\xs}$ la spécialisation du
groupe $p$-divisible universel.

\begin{figure}[htbp]
   \begin{center}
      \includegraphics{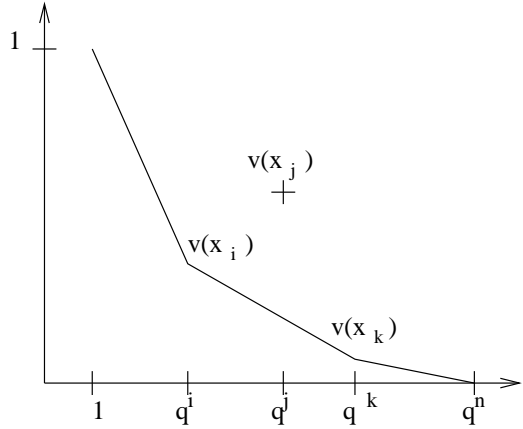}
   \end{center}
   \caption{\footnotesize Le polygone de Newton de
    la multiplication par $\pi$}
   \label{OSFshshftzr36zrp}
\end{figure}

On vérifie alors aisément que le polygone de Newton de la
multiplication par $\pi$ sur la loi de groupe formel universelle est
l'enveloppe convexe des $(q^i,
v(x_i))_{0\leq i \leq n}$ où $x_0=0$ et $x_n=1$ (cf. figure
\ref{OSFshshftzr36zrp}). Ses pentes non infinies $\l_1\geq \dots \geq
\l_n$, où $\l_i$ apparaît $q^i-q^{i-1}$ fois entre les abscisses
$q^{i-1}$ et $q^i$, sont les valuations des points de $\pi$-torsion non-nuls
dans $H_{\xs} [\pi] (\O_{\overline{L}})\setminus \{0 \}$. On peut donc lire
sur ce polygone de Newton la filtration de ramification de $H_{\xs} [\pi]
(\O_{\overline{L}})$. 

\begin{rema}
Notons encore $H_0$ le groupe $p$-divisible sur $\spec (\O [[\xs]])$
associé à celui sur $\spf (\O [[\xs]])$. Le choix de coordonnées
précédentes implique que la stratification de Newton de $H_0 \text{
  mod } \pi$ sur $\spec (k [[\xs]])$ est donnée par $$ V(x_1,\dots,
x_{n-1})\subset \dots \subset V(x_1, x_2)\subset V(x_1) $$
où $V(x_1,\dots,x_{n-1})$ est le lieu supersingulier (le point dont on
est parti et qu'on a déformé) et $V(x_1)^c$ le lieu ordinaire, 
et c'est essentiellement la seule propriété dont nous auront besoin
(on trouvera dans \cite{Oort1} la démonstration de l'existence de systèmes de coordonnées formelles sur des espaces de déformation plus généraux tels que la stratification de Newton soit définie par des sous-espaces
linéaires).  On vérifie en effet que cette seule propriété sur la
stratification de Newton implique que le polygone de Newton de $H_0
[\pi]$ est celui donné précédemment. 
\end{rema}

L'ensemble des points de $\X^{rig}$ où le polygone de Newton de
$H[\pi]$ est au dessus d'un polygone donné est un ouvert admissible
quasicompact. Cela permet de stratifier l'espace rigide $\X^{rig}$ par
de tels polygones. En quelques sortes les ``bonnes coordonnées'' sur
l'espace $\X^{rig}$ ne sont pas les $(v(x_i))_{1\leq i\leq n-1}$ mais
celles données par le polygone de Newton précédent.  

\section{Application des périodes}

\subsection{Définition}

\begin{defi}
On note $\breve{\pi}_1: \M^{[0]} \ldrt \mathbb{P}^{n-1}_{\breve{F}}$
l'application des périodes telle que définie dans la section 23 de
\cite{HopkinsGross} ou plus généralement le chapitre 5 de \cite{RZ}.
\end{defi}

Ce morphisme étale $\O_{D}^\times$-équivariant d'espaces rigides
est défini de la façon suivante. Soit $\E_1$ le cristal de Messing de 
$H$ comme objet de $\left ( (\X\otimes (\breve{\O}/p \breve{\O})\, /\,\spec (\breve{\O})\right )_{NCRIS}$
(le cristal algèbre de Lie de l'extension vectorielle universelle de
\cite{Messing1} sur le site cristallin nilpotent). Soit
$\E^{rig}_1$ l'isocristal convergent associé sur $\X^{rig}$.  Soit 
 $\E_2$ le cristal de Dieudonné de $\mathbb{H}$ sur $(\spec (\Fqb )/\spec (W (\Fqb) )
)_{NCRIS}$.

 Il y a un plongement $$\iota : W(\Fqb)\hookrightarrow \O_{\breve{F}}$$ fixé par le choix
de l'isomorphisme  entre le corps résiduel de $\breve{F}$ et $\Fqb$. 
Soit $$f: \X\ldrt \spec (\O_{\breve{F}}) \ldrt \spec  (W (\Fqb))$$ le morphisme composé.
Il y a  un diagramme 
$$
\xymatrix{ 
\X\otimes (\O_{\breve{F}}/p \O_{\breve{F}}) \ar[r]\ar[d]_{f \text{ mod } p} & \spec (\O_{\breve{F}}) \ar[d] \\
\spec (\Fqb) \ar[r] & \spec (W (\Fqb))
}
$$
qui induit donc un morphisme de topos
 $$f^{CRIS} : \left ( (\X\otimes (\O_{\breve{F}}/p \O_{\breve{F}})\, /\,\spec (\O_{\breve{F}})\right )_{NCRIS} \ldrt (\spec (\Fqb )/\spec (W (\Fqb) )
)_{NCRIS}$$

 La quasi-isogénie universelle $\rho : \mathbb{H}\times_{\Fqb}
 (\X\text{ mod } p) \ldrt H\times_\X (\X \text{
   mod p})$ sur l'espace des
déformations induit une quasi-isogénie de cristaux 
$$
\mathbb{D} (\rho) : f^{CRIS*} \E_2 \ldrt \E_1
$$
et donc un isomorphisme d'isocristaux convergents 
\begin{equation}\label{iso_rig1}
f^{rig *} (\E_2^{rig}) =(f^* \E_2)^{rig} \iso \E_1^{rig}
\end{equation}
Si $\mathbb{D} (\mathbb{H})$ désigne le module de Dieudonné covariant 
 ``classique'' de $\mathbb{H}$,
c'est à dire l'évaluation de $\E_2$ sur l'épaississement $\spec
(\Fqb)\hookrightarrow \spec (W (\Fqb))$, 
et 
$E(H)$ l'extension vectorielle universelle du groupe p-divisible universel $H$ (un $\O_{\X}$-module
libre de rang $n [F:\Qp]$), il y a donc un isomorphisme $\O_D^\times$-équivariant 
\begin{equation}\label{iso_riga2}
\mathbb{D} (\mathbb{H})_{\Q} \otimes_{W(\Fqb)_\Q,\iota } \O_{\X^{rig}} \simeq 
\text{Lie } (E(H))^{rig}
\end{equation}
 et de plus 
$$
\mathbb{D} (\mathbb{H})_{\Q}\otimes_{W(\Fqb)_\Q,\iota } \breve{F} = \left (\text{Lie} (E(H))^{rig} \right )^{\nabla =0}
$$
où $\nabla$ désigne la connexion de Gauss-Manin induite par la structure cristalline de
l'extension vectorielle universelle; si 
$$
\xymatrix{
 & \Delta^{(2)} \ar[ld]^{\text{pr}_1} \ar[rd]_{\text{pr}_2} \ar@{^(->}[r] & \X\times \X  \\
\X & & \X
}
$$
désigne le voisinage infinitésimal d'ordre $2$ de la diagonale de $\X$, la nature cristalline
de l'extension vectorielle universelle induit un isomorphisme 
$$
\Theta : \text{pr}_1^* E(H) \iso \text{pr}_2^* E(H)
$$
et alors 
$$
\nabla (m) = \text{Lie}(\Theta) (1\otimes m) - m\otimes 1 
$$

On renvoie à la section 5.3 de \cite{DeJong1} ou la proposition 2.3.26 de \cite{Laurent1} 
pour la construction des isomorphismes (\ref{iso_rig1}) et  (\ref{iso_riga2})  
qui repose sur la construction de Berthelot-Ogus
du foncteur des F-cristaux à isogénie près vers les isocristaux convergents. Cette construction
repose elle-même sur  l'astuce de Dwork; 
la structure de Frobenius permet d'agrandir les domaines de définition des solutions
de l'équation différentielle $\nabla=0$ où $\nabla$ est la connexion de Gauss-Manin. 

Il est également construit sans recours à la structure de F-cristal mais en utilisant la rigidité
des quasi-isogénies dans \cite{RZ} (cf. proposition 5.15 de \cite{RZ}).
\\

L'isomorphisme (\ref{iso_riga2}) est un isomorphisme de $F\otimes_{\Qp} \breve{F}$-modules via l'action de $\O_F$ sur $H$ et $\mathbb{H}$. Considérons la décomposition isotypique 
$$
\mathbb{D} ( \mathbb{H} )_\Q = \bigoplus_{\tau :F^0 \hookrightarrow W(\Fqb)_\Q}
\mathbb{D} (\mathbb{H} )_{\Q,\tau}
$$ 
Le morphisme $\iota : W(\Fqb)_{\Q} \hookrightarrow \breve{F}$ induit un plongement 
$\tau_0 : F^0\hookrightarrow  W(\Fqb)_{\Q}$. Celui-ci induit un
isomorphisme $\O\otimes_{\O_{F^0,\tau_0}} W(\Fqb)\iso W_\O (\Fqb)$ 
($W_\O$ désigne les vecteurs de Witt ramifiés, cf. \cite{DrinfeldOmega}).
Notons alors
$$
 \mathbb{D}_\O ( \mathbb{H} )_\Q =
\mathbb{D} ( \mathbb{H} )_{\Q,\tau_0}
$$
où $\O$ est là pour $\O_F$. C'est un $W_\O (\Fqb)$-module
 muni d'un Frobenius $\ph$ qui est
semi-linéaire relativement au Frobenius de $W_\O$.
 C'est un cristal relativement à $\O$, les cristaux usuels correspondant au cas $\O=\Zp$.

 Soit $\a:  F\otimes_{\Qp} \breve{F} \twoheadrightarrow \breve{F}$. Le plus grand quotient de $\mathbb{D} ( \mathbb{H} )_\Q\otimes_{W (\Fqb)_\Q,\iota} \breve{F}$ à travers lequel $F$ agit via $F\subset \breve{F}$ est
\begin{eqnarray*}\left ( 
\mathbb{D} ( \mathbb{H} )_\Q\otimes_{W (\Fqb)_\Q,\iota} \breve{F}\right ) \otimes_{F\otimes_{\Qp} \breve{F},\a} \breve{F}
 &=& \mathbb{D} ( \mathbb{H} )_\Q \otimes_{F\otimes W(\Fqb)_\Q} \breve{F} \\
&=& \mathbb{D}_\O ( \mathbb{H} )_\Q
\end{eqnarray*} 

Notons $\left ( \Lie\, (E(H))^{rig} \right )'$ le plus grand quotient de $\Lie\, (E(H))^{rig} $
sur lequel $\O$ agit à travers $F\subset \breve{F}$. L'isomorphisme $\O_D^\times$-équivariant
(\ref{iso_riga2}) induit alors un isomorphisme entre fibrés $\O_D^\times$-équivariants sur $\X^{rig}$
$$
\mathbb{D}_\O ( \mathbb{H})_\Q\otimes_{W_\O (\Fqb)} \O_{\X^{rig}} \simeq \left ( \Lie\, (E(H))^{rig} \right )'
$$
La filtration localement facteur directe 
$$
V(H) \subset \Lie\, E(H)
$$
définie par la partie vectorielle de l'extension universelle est telle que sur son quotient
$\Lie\, H$, $\O$ agisse à travers $\O\hookrightarrow \O_{\breve{F}}$. Le $F\otimes_{\Qp} \O_{\X^{rig}}$-module $\Lie\, E(H)^{rig}$ est libre d'après l'isomorphisme (\ref{iso_riga2}). On en déduit
(cf. par exemple la démonstration la proposition \ref{exOver} de l'appendice \ref{MdcO}) que si $I=\ker ( \O\otimes_{\Zp} \O \twoheadrightarrow \O)$ alors 
$$
\left ( V (H)^{rig} \right )' := V(H)^{rig}/I.\Lie\, E(H)^{rig} \hookrightarrow (\Lie\, E(H)^{rig})' 
$$
est une filtration localement facteur directe de codimension $1$. D'où une filtration 
$$
\left ( V(H)^{rig}  \right )' \subset \mathbb{D}_\O (\mathbb{H}) \otimes \O_{\X^{rig}}
$$
localement facteur directe de codimension un. Cette filtration définit
l'application des périodes
$$
\X^{rig} \ldrt \mathbb{P} (\mathbb{D}_\O (\mathbb{H}))
$$
Bien sûr cette application s'étend sur tout l'espace de Rapoport-Zink
en un morphisme étale $D^\times$-équivariant
$$
\M\ldrt \mathbb{P} (\mathbb{D}_\O (\mathbb{H}))
$$
défini de la même façon que précédemment en remplaçant $\X$ par
$\widehat{\mathcal{M}}$. L'espace des périodes $ \mathbb{P}
(\mathbb{D}_\O (\mathbb{H}))$ est en quelques sortes l'espace $\M$
quotienté par la relation d'isogénies, c'est à dire le quotient de la
tour de Lubin-Tate par le groupe $\GL_n (F)$. 

\subsection{Interprétation en termes du cristal $\O$-extension
  vectorielle universelle}

Les résultats de l'appendice \ref{MdcO} permettent de construire 
directement l'application des périodes précédente en interprétant 
$\mathbb{D}_\O (H)$ et $\Lie (H)/I.\Lie (H)$ comme évaluations d'un cristal
algèbre de Lie de la $\O$-extension vectorielle universelle sur un
site cristallin défini en termes de $\O$-puissances divisées.

Ils permettent par exemple de voir que l'application des périodes
est ``entière''sur le polydisque $\forall i\; v(x_i)\geq p-1$ alors
que la théorie ``classique''montre que ce n'est le cas que sur le
polydisque $\forall i\; v(x_i)\geq e(p-1)$ où $e$ est l'indice de
ramification de $F|\Qp$. 

\subsection{La donnée de descente sur l'espace des périodes}

Le Verschiebung sur le cristal $\mathbb{D}_\O (\mathbb{H})$ induit un
isomorphisme 
$$
\mathbb{D}_\O (\mathbb{H}) \iso \mathbb{D}_\O (\mathbb{H})^{(\s)}
$$
qui induit une donnée de descente compatible à l'action de $D^\times$
sur l'espace des périodes 
$$
 \mathbb{P} (\mathbb{D}_\O (\mathbb{H})) \iso  \mathbb{P} (\mathbb{D}_\O (\mathbb{H}))^{(\s)}
$$
L'application des périodes est compatible à la donnée de descente de
Rapoport-Zink sur $\M$ et à cette donnée de descente.
On remarquera que cette donnée de descente est effective et que la
variété descendue sur $F$ est la variété de Severi-Brauer associée à
l'algèbre à division $D$.

\subsection{Formules explicites pour l'application des périodes et applications}

Nous donnons dans cette section des formules explicites pour l'application des périodes
en utilisant la théorie des displays de \cite{Zink2}. Ces formules
sont plus simples à manipuler que celles utilisées dans
\cite{HopkinsGross} basées elles sur la théorie des
quasi-logarithmes. Bien que la théorie de \cite{Zink2} concerne des
groupes p-divisibles, les résultats de l'appendice
\ref{MdcO} permettent d'étendre celle-ci au cas des $\O$-modules
formels en remplaçant les vecteurs de Witt par les vecteurs de Witt
ramifiés 
(l'auteur n'affirme pas avoir revérifié chaque démonstration de
\cite{Zink2}, mais il en a revérifié suffisamment pour se convaincre
que cela marchait et il invite le lecteur à en faire de même). 

\subsubsection{Display universel sur $\X$}

Considérons le module de Cartier (pour la théorie de Cartier des $\O$-modules formels) 
 sur $\breve{\O}[[x_1,\dots,x_{n-1}]]$ ayant pour
$V$-base $e$ et comme équation structurelle 
$$
F.e= [x_1] e + V[x_2] e +\dots + V^{n-2}[x_{n-1}] e + V^{n-1} e
$$
La loi de $\O$-module formel associée est la loi universelle considérée dans \cite{HopkinsGross}.
Ce module de Cartier provient d'un $\O$-Display $(P,Q,F,V^{-1})$ où
$$
P= L\oplus T\;\;\; T=<\e_1>\;\;\; L=<\e_2,\dots,\e_{n}>\;\;\; Q= L\oplus I_{W_\O} T
$$
$ I_{W_\O}$ désignant l'idéal d'augmentation des vecteurs de Witt, $T$
est un relèvement de l'espace tangent 
et la matrice de $F\oplus V^{-1}$ dans la base $(\e_i)_i$ est 
$$
\left ( 
\begin{matrix}
[x_1] &  [x_2] & \dots & [x_{n-1}] & 1 \\
1 & 0 & \dots & 0 & 0 \\
0 & 1 &  & 0 & 0 \\
 & 0 & \ddots & \vdots & \vdots \\
 & & & 1 & 0
\end{matrix}
\right ) \in \GL_n\left  ( W_\O ( \breve{\O} [[ x_1,\dots, x_{n-1}]])\right )
$$
Cette matrice s'écrit 
$$
\left ( 
\begin{matrix}
 1 & \vline & [x_1] & \dots & [x_{n-1}] \\
\hline 0  &  \vline\\
 \vdots & \vline &&  \text{\Large I}_{n-1}  \\
0 & \vline
\end{matrix}
\right )
\left ( 
\begin{matrix}
0 & \dots & 0 & 1 \\
1 &       &   &  0 \\
  & \ddots &  & \vdots  \\
  &        &  1 & 0
\end{matrix}
\right )
$$
comme dans la formule (86) page 174 de \cite{Zink2}, où la matrice de
droite est ici la matrice du cristal généralisé (cf. la section
\ref{okklu} de l'appendice \ref{MdcO}) du $\O$-module
formel $\mathbb{H}$ que l'on déforme. La matrice de l'opérateur $F$ s'écrit alors
$$
\left ( 
\begin{matrix}
[x_1] &  \pi [x_2] & \dots & \pi [x_{n-1}] & \pi \\
1 & 0 & \dots & 0 & 0 \\
0 & \pi &  & 0 & 0 \\
 & 0 & \ddots & \vdots & \vdots \\
 & & & \pi & 0
\end{matrix}
\right )
$$
 
\subsubsection{Formule pour l'application des périodes}

Notons 
$$
A= 
\left ( 
\begin{matrix}
x_1 &  \pi x_2 & \dots & \pi x_{n-1} & \pi \\
1 & 0 & \dots & 0 & 0 \\
0 & \pi &  & 0 & 0 \\
 & 0 & \ddots & \vdots & \vdots \\
 & & & \pi & 0
\end{matrix}
\right ) 
$$
la réduite modulo l'idéal d'augmentation de $W_\O$  de la matrice précédente 
et $\forall i\in \N \;\; A^{(\s^i)}$ la
matrice obtenue en remplaçant $\forall k\; x_k$ par $x_k^{q^i}$. Soit 
$$
B=
\left ( 
\begin{matrix}
0 & 0 & \dots & 0 & \pi \\
1 &       &   &  0 & 0 \\
  & \pi &  &  & \vdots  \\
  &      &\ddots && \vdots  \\
 && &  \pi & 0
\end{matrix}
\right )
$$
la matrice de l'opérateur $F$ du cristal de $\mathbb{H}$. Alors,
d'après la proposition 71 de \cite{Zink2} 
\begin{eqnarray}\label{foper}
\underset{k\drt + \infty}{\lim} (1, 0, \dots , 0). A A^{(\s)} \dots A^{(\s^{k-1})} B^{-k}
\end{eqnarray}
existe dans $\GG ( \X^{rig}, \O_{\X^{rig}})^n$ où l'anneau $\GG ( \X^{rig}, \O_{\X^{rig}})$
est muni de sa topologie de Frechet usuelle de la convergence uniforme
sur
toutes les boules fermées de rayon plus petit que un. Notons $(f_0,\dots,  f_{n-1})$ cette limite.
Il s'agit également d'une limite au sens de la topologie $(x_1,\dots,x_{n-1})$-adique.
 Plus précisément 
$$
(f_0,\dots,f_{n-1}) \equiv (1,0\dots 0).
 A \dots A^{(\s^{k-1})} B^{-k} \;\;\text{mod } (x_1^{q^{k}}, \dots,
 x_{n-1}^{q^{k}})
$$
On peut également récrire la limite précédente sous la forme 
$$
\underset{l\drt +\infty}{\lim} \;\pi^{-l (n-1)}\left (1, 0 \dots 0\right ) A A^{(\s)} \dots A^{(\s^{ln -1})}
$$
Alors, toujours d'après la proposition 71 de \cite{Zink2}, 
$$
\breve{\pi}_1 = [f_0 : \dots :  f_{n-1}]
$$
et on peut récrire 
$$
\breve{\pi}_1 = \underset{l\drt +\infty}{\lim} [1 : 0 : \dots : 0]. A A^{(\s)}\dots A^{(\s^{ln-1})}
\in \mathbb{P}^{n-1} \left ( \GG (\X^{rig}, \O_{\X^{rig}}) \right )
$$
qui exprime $\breve{\pi}_1$ comme une limite d'orbite ``$\s$-linéaire''
dans $\mathbb{P}^{n-1}$  sous l'action de $A\in \text{PGL}_n$.

\begin{exem}
Lorsque $n=2$, si 
$$
<a_1,\dots,a_i>=\frac{1}{a_1+ \frac{1}{a_2+ \dots \frac{1}{a_i}}}
$$
alors $\breve{\pi}_1 = [1: f(x)]$ où 
$$
f = \underset{k \drt +\infty}{\lim} <\frac{x^{q^{2k}}}{\pi},{x^{q^{2k-1}}}, 
 \frac{x^{q^{2k-2}}}{\pi},\dots, x^q, \frac{x}{\pi}>
$$
\end{exem}

\subsubsection{Application : généralisation d'un théorème de Gross-Hopkins}

Nous généralisons ici le corollaire 23.15 de l'article
\cite{HopkinsGross} en utilisant les formules précédentes. 
Ce théorème a déjà été obtenu dans \cite{Yu} en
utilisant les quasi-logarithmes et n'est donc pas nouveau. 

\begin{theo}[Yu]\label{periso}
Le morphisme $\breve{\pi}_1$ induit un isomorphisme entre les ouverts
admissibles 
$$
\{\; (x_1,\dots,x_{n-1})\in \X^{rig}\;|\; \forall 1\leq i\leq n\;\;
\forall 0\leq j \leq n-1 \;\; \frac{1-v(x_i)}{q^n (q^i-1)}< \frac{v(x_j)}{q^n-q^j}
 \; \}
$$
et
$$
\{\; [w_0:\dots : w_{n-1}] \in \mathbb{P}^n \;|\;  w_0\neq 0 \text{ et
  } \forall 1\leq i\leq n\;\;
\forall 0\leq j \leq n-1 \;\; \frac{1-v(\widetilde{w}_i)}{q^n (q^i-1)}< \frac{v(\widetilde{w}_j)}{q^n-q^j}
 \;\}
$$
où on a posé $x_0=\pi, x_n=1, \widetilde{w}_n=1, \widetilde{w}_0=\pi$
et $\forall 1\leq i\leq n-1\;\; \widetilde{w}_i = \frac{w_i}{w_0}$.  
De plus sur ces domaines 
$$
\forall 1\leq i\leq n-1 \;\; v(x_i) = v\left (\frac{f_i(x_1,\dots,x_{n-1})}{f_0(x_1,\dots,x_{n-1})} \right )
$$
\end{theo}
\dem 
On utilise la formule limite (\ref{foper}) pour l'application des périodes. Remarquons 
qu'avec les notations de la sous-section précédente
 $A=CB$ où 
$$
C=
\left ( 
\begin{matrix}
 1 & \vline & x_1 & \dots & x_{n-1} \\
\hline 0  &  \vline\\
 \vdots & \vline &&  \text{\Large I}_{n-1}  \\
0 & \vline
\end{matrix}
\right )
$$
et que donc, si l'on pose 
$$
(f_0^{(k)},\dots, f_{n-1}^{(k)}) = (1,0 \dots 0) . A A^{(\s)}\dots A^{(\s^{k-1})} B^{-k}
$$
alors on a la formule de récurrence 
\begin{eqnarray*}
(f_0^{(k+1)},\dots, f_{n-1}^{(k+1)})
 & = & 
 (f_0^{(k)},\dots, f_{n-1}^{(k)}). B^k A^{(\s^k)} B^{-k-1} \\
&=& (f_0^{(k)},\dots, f_{n-1}^{(k)}). B^b C^{(\s^k)} B^{-b}
\end{eqnarray*}
où $k=an+b$, $b\in\{ 0,\dots, n-1 \}$. On en déduit les formules de récurrence suivantes :
$$
\text{Si } b=0 \;\;\left \{ 
\begin{array}{ll}
\hspace{1cm} f_0^{(k+1)} & =  f_0^{(k)} \\
\forall i>0\; f_i^{(k+1)} &= f_i^{(k)} + x_i^{q^k} f_0^{(k)}
\end{array}
 \right.
$$
Si $b\neq0$, en posant $x_0=1$ et $\forall i\in \Z \; x_i=x_j$ où $j\equiv i \text{ mod } n,
j\in \{ 0,\dots, n-1\}$
$$
\left \{
\begin{array}{ll}
 \hspace{1cm} f_b^{(k+1)} &= f_b^{(k)} \\
\forall i\neq b \; f_i^{(k+1)} &= f_i^{(k)} + \pi^{\a (b,i)} x_{n-b+i}^{q^k} f_b^{(k)}
\end{array}
\right.
$$
où 
$$
\a (b,i) = \left \{ \begin{array}{lll}
-1 &\text{ si }& i=0  \\
\;\;\; 0 &\text{ si }& 1 \leq i \leq b-1  \\
-1 &\text{ si }& b+1 \leq i \leq n-1 
\end{array}
\right.
$$
Pour un  $\e$ tel que  $0<\e<1$ et pour une fonction
$g\in \breve{\O}[[x_1,\dots,x_{n-1}]][\frac{1}{\pi}]$ posons 
$$
V_\e (g) = \sup \;\{\; v ( g (\xs))\; |\;
 \xs \in \X^{rig}\;\;\forall 1\leq i\leq n\;\;
\forall 0\leq j \leq n-1 \;\; \frac{1-v(x_i)}{q^n (q^i-1)}\leq \e \frac{v(x_j)}{q^n-q^j}
 \;\}
$$
Il suffit de montrer que 
$$
\forall \e\;
\forall k>0 \; \left \{ 
\begin{array}{ll}
 V_\e (f_0^{(k)}-1) > 0 \\ \text{et } \forall 1\leq i \leq n \;
V_\e (f_i^{(k)}-x_i) > V_\e (x_i)
\end{array}
\right.
$$
($\breve{\O}$ étant de valuation discrète cela implique qu'il existe $\a>0$ tel que
$\forall k>0\;  V_\e (f_i^{(k)}-x_i)> \a + V_\e (x_i)$). Cela ne pose pas 
de problème en utilisant les formules de récurrence données précédemment.
\qed

\begin{rema}
L'ouvert admissible précédent est l'ensemble des points  $\xs\in\X^{rig}$ où le polygone de Newton
de $H_{\xs}[\pi]$ de pentes $\l_1\geq \dots \geq \l_n$
vérifie 
$$
\frac{\l_1}{q^n} < \l_n
$$
ou encore les valuations des éléments de $H_{\xs} [\pi^2]\setminus
H_{\xs} [\pi]$ sont strictement plus petites que celles de $H_{\xs}
[\pi]$. 
Cela se déduit des égalités 
$$
\l_1 = \text{sup} \{ \frac{1-v(x_i)}{q^i-1}\; |\; 1\leq i\leq n \}
\text{ et } \l_n= \inf \{ \frac{v(x_j)}{q^n-q^j}\; |\; 0\leq j\leq n-1 \}
$$
Il s'agit en fait d'un domaine fondamental ``ouvert''  pour les isogénies déformant une puissance de $\pi$
i.e. pour les opérateurs de Hecke non-ramifiés de degré un multiple de
$n$,  le domaine fondamental de Gross-Hopkins défini dans la section
suivante étant quant à lui un domaine fondamental pour les isogénies quelconques i.e. tous les opérateurs de Hecke non-ramifiés. 
Partant de ce point de vue on  donne une autre démonstration plus
conceptuelle d'une partie de ce théorème dans la section \ref{erggi}. 

Plus généralement dans l'article \cite{Rami} nous étudions plus en
détails l'application des périodes en dehors de cet ouvert admissible
comme dans \cite{Yu}.
\end{rema}

\section{Domaine fondamentale de Lafaille/Gross-Hopkins (\cite{HopkinsGross})}

\begin{defi}
Posons
$$
\mathcal{D}=\{ \xs \in \X^{rig}
\; |\; \forall 1\leq i \leq n-1\;\;
v(x_i)\geq 1- \frac{i}{n}\;\}
$$
un ouvert admissible quasi-compact dans la boule unité ouverte de dimension $n-1$. 
\end{defi}

Après extension des scalaires à $\breve{F}( \pi^{1/n})$ l'espace 
$\mathcal{D}$ devient isomorphe à la boule unité fermée de rayon $1$ et de dimension $n-1$.

\begin{figure}[htbp]
   \begin{center}
      \includegraphics{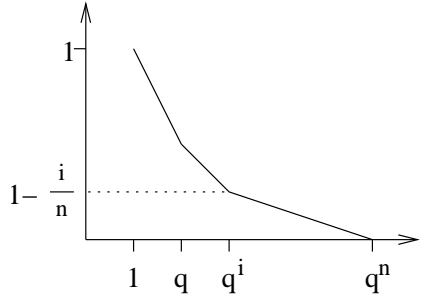}
   \end{center}
   \caption{\footnotesize Le polygone de Newton 
      bordant le domaine fondamental de Gross-Hopkins}
   \label{NewGH}
\end{figure}

\begin{rema}
Le domaine $\mathcal{D}$ admet la description intrinsèque (i.e.  indépendante du choix de coordonnées) modulaire
suivante : 
$$
\mathcal{D} =\{ \xs\in\X^{rig}
\;|\; \text{polygone de Newton de }
 H_{\xs} [\pi ] \text{ est } \geq \text{ à celui de la figure }\ref{NewGH} \;\}
$$
En particulier il est stable sous $\O_D^\times$ puisqu'il ne dépend
pas de la déformation $\rho$ mais seulement de la classe
d'isomorphisme de $H_{\xs}$. 
\end{rema}

Rappelons la proposition suivante qui est un cas particulier du
théorème \ref{periso} :

\begin{prop}[Gross-Hopkins  \cite{HopkinsGross}]
Le morphisme étale $\breve{\pi}_1$ restreint à $\mathcal{D}$ 
induit un isomorphisme entre 
$\mathcal{D}$ et l'ouvert suivant 
$$
\breve{\pi}_1 (\mathcal{D})=\{ \; [y_0:\dots: y_n]\;|\; 
\forall i\;\;
v(\frac{y_i}{y_0})\geq 1 - \frac{i}{n}\; \}
$$
\end{prop}

Pour toute extension valuée $L | F$  
les fibres du morphisme $\breve{\pi}_1$ sur les $L$-points 
correspondent aux 
$\O$-modules $\pi$-divisibles isogénes via une isogénie déformant une puissance 
de $\pi$ en fibre spéciale. Cependant il existe des points distincts de $\mathcal{D}$ correspondant à des $\O$-modules $\pi$-divisibles isogénes via des isogénies déformant une puissance de $\Pi$, 
 ce que précisent les définitions et propositions suivantes. En
 d'autres termes on cherche à comprendre comment se recolle
 $\mathcal{D}$ avec ses itérés sous les correspondances de Hecke sphériques.

\begin{defi}\label{borde}
Pour $1\leq i \leq n-1$ posons 
$$
\partial_i \mathcal{D} = \{\; \xs\in \mathcal{D}\; |\;
v(x_i)=1-\frac{i}{n}\;\}
$$
un domaine de Laurent dans $\mathcal{D}$. 
\end{defi}
Remarquons que sur $\partial_i\mathcal{D}$ le polygone de Newton de $H_{\xs}[\pi]$ possède un point 
de rupture en $q^i$ ce qui correspond à un l'existence d'un cran de
rang $i$ dans la filtration de ramification du schéma en groupes
$H_{\xs}[\pi]$ ou encore à l'existence d'un sous-groupe canonique
``généralisé''.

\begin{prop}[Faltings]\label{pofgh}
Soient $\xs\in\mathcal{D}$, $\xs'\in \X^{an}$  
et $f:H_{\xs}\drt H_{\xs'}$ une isogénie
qui n'est pas un isomorphisme et qui ne se factorise pas par la
multiplication par $\pi$. 
 Le
points $\xs'$ appartient à $\mathcal{D}$ ssi $\exists i\;
\xs\in \partial_i \mathcal{D}$ et $f$ est  une isogénie de noyau 
 le sous-groupe de rang $i$ de $F_{\xs}[\pi]$ formé des $q^i$-points de plus grande valuation (un sous-groupe canonique généralisé).
\end{prop}
\dem

Rappelons que si $L|F$ est une extension valuée complète et
$\ph:H_1\ldrt H_2$ est une isogénie entre deux $\O$-modules
$\pi$-divisibles formels de dimension $1$ alors 
$$
\forall x\in \widehat{H}_1 (\O_L)\;\; v(\ph (x))=\sum_{\a\in \ker \ph } v( x-\a)
$$
Supposons que $\xs\in \partial_i\mathcal{D}$ et soit $f:H_{\xs} \drt
H_{\xs'}$ l'isogénie définie par le quotient des $q^i$-points de plus
grande valuation. Notons $\l_1\geq \dots\geq \l_n$ les pentes du
polygone de Newton de $H_{\xs}[\pi]$ ($\l_k$ est la pente entre les
abscisses $q^{k-1}$ et $q^k$).

 Les éléments de $H_{\xs'}[\pi]$ sont de deux types :
\begin{itemize}
\item de la forme $f(\a)$ où $\a\in F_{\xs}[\pi]\setminus \ker f$. Alors, leurs valuations sont (puisque $\forall \beta\in\ker f\; v(\beta)>v(\a)$) $v(f(\a))=q^i v(\a)$ et varient donc dans $(q^i\l_{i+1},\dots,q^i \l_n)$.
\item de la forme $f(\a)$ où $\a$ est un zéro de la série formelle 
$[\pi]_{F_{\xs}} -\beta$ avec $\beta\in \ker f\setminus \{ 0 \}$
et $F_{\xs}$ désigne une loi de groupe formel associée à $H_{\xs}$. 
 Le polygone de Newton de cette série formelle est obtenu en prenant l'enveloppe convexe de $(0,v (\beta))$ et celui de $[\pi]_{F_{\xs}}$
\end{itemize}

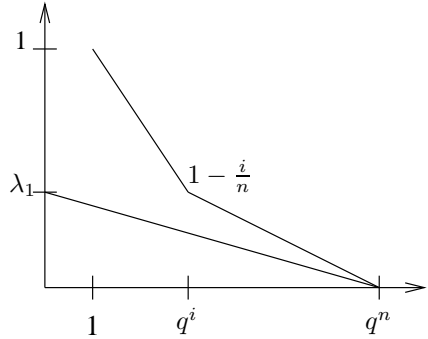
\begin{figure}[h]
   \begin{center}
      \input{Newton2.pstex_t}
   \end{center}
   \caption{\footnotesize Le polygone de Newton de
     $[\pi]_{F_{\xs}}-\beta$, $v(\beta)=\l_1$ }
   \label{figop}
\end{figure}

Mais, 
$\dpt{\underset{\beta\in\ker f\setminus \{ 0\}}{\text{sup}} v(\beta)=\l_1}$. Or, $\dpt{\frac{\l_1}{q^n}<\l_n}$. En effet, 
$$
v(x_1)\geq 1-\frac{1}{n} \impl \l_1 \leq \frac{1}{n(q-1)} \text{ et }
v(x_n)\geq \frac{1}{n} \impl \l_n \geq \frac{1}{n (q^n -q^{n-1})}
$$
Donc, $\forall \beta\in \ker f \setminus \{0\} \; \text{Newt} ([\pi]_{F_{\xs}}-\beta)$ est le segment joignant $(0,v (\beta))$ à $(q^n ,0)$ avec $v(\b)\in \{\l_1,\dots,\l_i \}$ (figure \ref{figop}).
 Les valuations des éléments de $[\pi]_{F_{\xs}}^{-1} (\ker f \setminus \{ 0\})$ sont donc $\dpt{(\frac{\l_1}{q^n},\dots,\frac{\l_i}{q^n})}$. Or, $\dpt{\frac{\l_1}{q^n}<\l_i}$ puisque $\l_1 \leq \dpt{\frac{1}{n(q-1)}}$ et l'égalité $v(x_i)=\dpt{1-\frac{i}{n}}$ implique 
que $\l_i\dpt{\geq \frac{1}{n (q^i-q^{i-1})}}$. Donc, $$\forall \beta \in \ker f \setminus \{0 \}\;\forall \a\in [\pi]_{F_{\xs}}^{-1} (\ker f \setminus \{ 0\}) \;\; v(\a)< v(\beta) 
$$
On en déduit que $v(f(\a))=q^i v(\a)$ qui sont donc les $(\dpt{\frac{\l_1}{q^{n-i}},\dots, \frac{\l_i}{q^{n-i}}})$. 
Or, 
$$
q^i \l_{i+1} \geq \dots \geq q^i \l_n > \frac{\l_1}{q^{n-i}} \geq \dots \geq \frac{\l_i}{q^{n-i}}
$$
Donc, $\xs'\in \partial_{n-i} \mathcal{D}$.  
\\
    
Réciproquement, soit $f:H_{\xs}\ldrt H_{\xs'}$ une isogénie comme dans l'énoncé telle que
$\xs,\xs'\in\mathcal{D}$. 
\\
Commençons par remarquer que les analyses précédentes montrent que les
valuations des éléments non nuls  de $H_{\xs} [\pi^k]$ sont 
$$
\l_1\geq \dots \geq \l_n > \frac{\l_1}{q^n} \geq\dots\geq
\frac{\l_n}{q^n}> \dots >\frac{\l_1}{q^{n(k-1)}} \geq \dots \geq \frac{\l_n}{q^{n(k-1)}}
$$
où il y a $(q^i-q^{i-1})q^{n(k-1)}$ éléments de valuations
$\dpt{\frac{\l_i}{q^{n(k-1)}}}$. 
\\
Soit $M=\ker f$. Notons 
$$
k=\text{sup}\{ j\;|\; M[\pi^j]\neq M[\pi^{j-1}]\; \}
$$
Pour $j$ entre $1$ et $k$ notons 
$$
r_j=\dim_{\Fq} M[\pi^j]/M[\pi^{j-1}]
$$
On a donc 
$$
n-1 \geq r_1\geq r_2\geq \dots\geq r_k >0
$$
Posons pour simplifier les notations $r=r_k$. Le
$\Fq$-e.v. $\pi^{k-1}.M\subset H_{\xs}[\pi]$ est de
dimension $r$. Choisissant un drapeau complet de $\pi^{k-1}.M$
raffinant la filtration
de ramification inférieure (celle donnée par la valuation) on en
déduit que les valuations des éléments de
$[\pi^{k-1}]_{F_{\xs}}(M)\setminus \{ 0 \}$ sont de la forme 
$$
\underbrace{\l_{a_1}}_{q-1\text{ élts.}} \geq \dots \geq
\underbrace{\l_{a_r}}_{(q^r - q^{r-1})\text{ élts.}}
$$
où $\forall j\; a_j \geq j$. Soit $M'\subset H_{\xs} [\pi^{k+1}]$ un
sous-$\O_K$-module tel que $M=M'[\pi^k]$ et $\pi M' =M$. Les
valuations des éléments de $M'\setminus M$ sont 
$$
\frac{\l_{a_1}}{q^{nk}} \geq \dots \geq \frac{\l_{a_r}}{q^{nk}}
$$
et sont en particulier strictement inférieures à celles des éléments de
$M$. On en déduit que $f(M')\subset H_{\xs'}[\pi]$ est un $\Fq$-e.v. de
dimension $r$ dont les valuations des éléments non-nuls sont 
\begin{eqnarray}\label{fvap}
\underbrace{\frac{\l_{a_1}}{q^{nk- (r_1+\dots+ r_k)}}}_{(q-1)\text{ élts.}} \geq \dots \geq
\underbrace{\frac{\l_{a_r}}{q^{nk - (r_1+\dots + r_k)}}}_{(q^r  - q^{r-1})\text{élts.}}
\end{eqnarray}
Et que donc 
\begin{eqnarray*}
\sum_{\a\in f (M')\setminus \{0 \}} v(\a)  &=&
\frac{1}{q^{nk-(r_1+\dots + r_k)}} \left ( (q-1)\l_{a_1} + \dots +
  (q^r - q^{r-1}) \l_{a_r} \right ) \\
& \leq & \frac{1}{q^{nk - (n-1)(k-1) -r } } \underbrace{\left ( (q-1)\l_1+ \dots +
  (q^r - q^{r-1}) \l_r \right )}_{\leq \frac{r}{n} \text{car } \xs \in
\mathcal{D}} \\
&\leq & \frac{r}{n q^{n+ (k-1) -r}}
\end{eqnarray*}
Mais pour un sous-$\Fq$-ev. $N$ de $H_{\xs'} [\pi]$ de dimension $r$
on a 
$$
\sum_{\a\in N\setminus \{ 0 \}} v(\a) \geq  \frac{r}{n q^{n-r}}
$$
car si $\underbrace{\l_{b_1}}_{(q-1)\text{ élts.}} \geq \dots\geq
\underbrace{\l_{b_r}}_{(q^r-q^{r-1})\text{ élts.}}$ sont les
valuations des éléments de $N$ 
\begin{eqnarray*}
\sum_{\a\in N\setminus \{0 \}} v(\a)  &=& (q-1) \l_{b_1} + \dots +
(q^r - q^{r-1}) \l_{b_r} \\
&\geq & (q-1) \l_{n-r+1} + \dots + (q^r -q^{r-1}) \l_n \\
 &=& \frac{1}{q^{n-r}} \left [(q^{n-r+1}-q^{n-r}) \l_{n-r+1} + \dots +
   (q^n - q^{n-1}) \l_n \right ] \\
& \geq & \frac{r}{n q^{n-r}} \text{ car } \xs'\in\mathcal{D}
\end{eqnarray*}
De tous cela on déduit nécessairement que $k=1$. Si de plus il existe
un indice $j$ tel que $a_j>j$ alors les deux inégalités ci-dessus sont
strictes. On en déduit donc que $\forall j \; a_j=j$.
\qed

\begin{rema}
Dans \cite{Rami} nous expliquons plus généralement comment construire
des domaines fondamentaux comme celui de Gross-Hopkins à partir de
domaines fondamentaux pour l'action du groupe des rotations engendré
par le cycle $(1 \dots n)$ dans le simplexe de sommets $1,\dots, n$.
\end{rema}

\subsection{Lien entre le domaine fondamental et les points C.M.}

\begin{figure}[h]
   \begin{center}
      \includegraphics{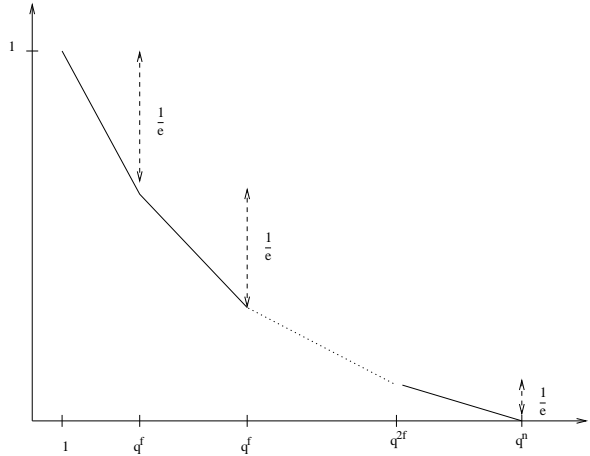}
   \end{center}
   \caption{ Le polygone de Newton d'un groupe C.M. }
   \label{CMf}
\end{figure}

\begin{prop}
Soit $E | F$ une extension de degré $n$. Il existe une unique classe
d'isogénie de groupes $p$-divisible  $H_{\xs}$, $\xs \in \X^{rig}$,
ayant multiplication complexe par un ordre dans $\O_E$. Les
représentants de cette classe d'isogénie dans $\mathcal{D}$ sont
exactement ceux ayant multiplication complexe par 
l'ordre maximal $\O_E$. De plus les $\xs \in \mathcal{D}$ tels que
$H_{\xs}$ ait multiplication complexe par $\O_E$ forment une
$\O_D^\times$-orbite. 
\begin{itemize}
\item
Si $E|F$ est non-ramifiée cette orbite est $\pi \O_E^{n-1}\subset
\mathcal{D}$. Le polygone de Newton associé a une seule pente.
\item Si $E|F$ est ramifié de degré $e>1$ alors cette orbite est
contenue dans 
$$
\bigcap_{k=1}^{e-1} \partial_{ k \frac{n}{e}} \mathcal{D}
$$
Le polygone de Newton a alors comme pentes 
$$
\forall k\in \{0,\dots, e-1\}\;\; \l_{kf+1}=\l_{kf+2}=\dots 
=\l_{(k+1)f-1} = \frac{1}{e (q^{(k+1)f}- q^{kf})}
$$
cf. la figure \ref{CMf}
\end{itemize}
\end{prop}
\dem
Facile.
\qed

\section{Une autre démonstration d'une partie du théorème \ref{periso}}\label{erggi}

Soit 
\begin{eqnarray*}
\mathcal{H} &=& \{ \;\xs\in \X^{rig}\; |\; \frac{\l_{1,\xs}}{q^n} <\l_{n,\xs} \;\}
\end{eqnarray*}
où $\l_{1,\xs}\geq \dots \geq \l_{n,\xs}$ sont les pentes de
$\text{Newt} (H_{\xs}[\pi])$. 
Il s'agit du domaine fondamental introduit dans l'énoncé du théorème
\ref{periso}. Nous donnons une autre démonstration de la partie
``isomorphisme sur son image'' du théorème \ref{periso}. 
 Cette démonstration  utilise
moins de calculs explicites et est plus conceptuelle. Cependant elle
ne permet pas de démontrer les propriétés ``métriques''de l'application
des périodes, par exemple calculer l'image par celle-ci des
polydisques ``fermés''inclus dans $\mathcal{H}$.
\\

Commençons par un lemme :

\begin{lemm}
Soit $K$ un corps valué complet non-archimédien et $f: X\ldrt Y$ un
morphisme étale entre $K$-espaces analytiques de Berkovich. C'est un
isomorphisme sur son image ssi $\forall L|K$ une extension de corps
valués le morphisme $f$ induit une injection de $X(L)$ dans $Y(L)$.
\end{lemm}
\dem
C'est une conséquence du fait qu'un morphisme étale entre espaces
analytiques est un isomorphisme local ssi il induit un isomorphisme
au niveau des extensions de corps résiduels.
\qed

D'après le lemme précédent le théorème \ref{periso} résulte de la
proposition suivante :

\begin{prop}
Soit $L|F$ une extension valuée complète.
Soient $\xs,\xs'\in \mathcal{H}(L)$ tels qu'il existe une isogénie 
$f : F_{\xs}\ldrt F_{\xs'}$ déformant une puissance de $\pi$. Alors 
$\xs = \xs'$. 
\end{prop}
\dem 
Commençons par constater que $\forall \xs \in \mathcal{H} (L)$, si 
$\l_1\geq \dots\l_n$ sont les pentes de
$\text{Newt}(H_{\xs}[\pi])$ alors les valuations des points de 
$H_{\xs} [\pi]\setminus \{ 0 \}, H_{\xs}[\pi^2]\setminus H_{\xs}
[\pi],\dots, H_{\xs} [\pi^k]\setminus H_{\xs}[\pi^{k-1}]$ sont ``strictement 
ordonnées''et valent 
$$
\l_1\geq \dots\geq \l_n >\frac{\l_1}{q^n}\geq\dots \geq\frac{\l_n}{q^n}
>\dots > \frac{\l_1}{q^{(k-1)n}} \geq \dots \geq\frac{\l_n}{q^{(k-1)n}} 
$$
Reprenons les notations de la seconde partie de la démonstration de la
proposition\ref{pofgh} : on considère $f,M,M'$... on obtient alors
d'après la formule \ref{fvap} qu'il existe des éléments dans $H_{\xs'}
[\pi]$ de valuation 
$$
\frac{\l_i}{q^{nk - (r_1+\dots + r_k)}} \text{   pour un }
i\in\{1,\dots, n\}
$$
où rappelons que $1\leq r_k\leq \dots \leq r_1 \leq n-1$. Par 
hypothèse 
$$
n | \text{ht} f = r_1+\dots + r_k
$$
Donc,
$$
nk - ( r_1+\dots + r_k ) = n  (\underbrace{ k - \frac{ r_1+\dots +
  r_k}{n}}_{\geq k (1-\frac{n-1}{n}) =\frac{k}{n}>0})
 \geq n
$$
On en déduit que
$$
\frac{\l_i}{q^{nk - (r_1+\dots + r_k)}} \leq \frac{\l_i}{q^n}
$$
Mais on vérifie que 
\begin{eqnarray*}
\forall \xs\in \mathcal{H}\;\; \forall j\; \l_j &\leq
&\frac{q}{q^n-1} \\
\text{ et } \l_j &>& \frac{1}{q^{n-1} (q^n-1)}
\end{eqnarray*}
\qed

\begin{rema}
Dans \cite{Rami} nous donnerons une démonstration conceptuellement
plus satisfaisante de ce théorème basée sur de la combinatoire dans
un appartement de l'immeuble de Bruhat-Tits de $\text{PGL}_n (F)$. 
\end{rema}

\section{L'espace des paramètres de la décomposition cellulaire}

Il s'agit de l'espace $D^\times \times GL_n (F)$-équivariant
 qui va indexer les cellules et qui est l'immeuble de Bruhat-Tits 
noté $\mathcal{I}$ du groupe $p$-adique suivant sur $F$ 
$$
(\GL_{n/F} \times D^\times)/\Gm 
$$
où $\Gm$ agit diagonalement via $z\mapsto z \text{Id}$ dans $\text{GL}_n (F)$ et 
$z\mapsto z^{-1}$ dans $D^\times$. Il admet la description concrète suivante :
$$
\mathcal{I} = \{\; (\La, M) \;\}/\sim \;\text{ où } \; (\La, M)\sim (\pi \La, \pi^{-1} M )
$$
avec $\La$ un réseau dans $F^n$ et $M$ un réseau $\O_D$-stable dans $D$. On notera $[\La,M]$ la classe de $(\La,M)$. 
En d'autres termes 
$$
\mathcal{I}= \mathcal I (\GL_n \times D^\times)/\pi^\Z
$$
avec 
\begin{itemize}
\item $\mathcal{I} (\GL_n) = \{ \text{ réseaux } \La \text{ dans le module de Tate } \}$
\item $\mathcal{I} (D^\times) = \{ \text{ cristaux } M \text{ dans le module de Dieudonné rationnel } \mathbb{D}_\O (\mathbb{H})[\frac{1}{p}]\; \}$
\item L'action de $\pi$ sur le module de Tate est la m\^eme que celle de
$\pi^{-1}$ sur le module de Dieudonné.
\end{itemize}
On ne considérera pas toute la structure simpliciale de cette immeuble 
(i.e. pas toutes les relations d'incidence) mais seulement les
arêtes  orientées suivantes : 

\begin{defi}
Soient $a,a'\in \mathcal{I}$. On note 
$$
a' \ldrt a \text{ s'il existe des représentants }a=[\La,M] \text{ et } a'=[\La',M] \text{  
tels que }$$
$$
\La'\varsubsetneq \La  \varsubsetneq \pi^{-1} \La'\;\text{ et }\; M'=\Pi^{-[\La : \La']} M
$$
\end{defi}

\begin{defi}
L'action de $\GL_n (F)\times D^\times$ sur $\mathcal{I}$
sera
$$
\forall (g,d)\in \GL_n (F)\times D^\times\;\;\; (g,d).[\La,M]=
	[g^{-1}\La, d.M]
$$
\end{defi}

\section{Les cellules rigides en niveau fini}

\subsection{Digression philosophique} 

Supposons que l'on veuille reconstruire  $\GL_n (F)$ à partir 
de l'immeuble $\GL_n (F)/\GL_n (\O_F)$ et de la cellule à l'origine
$\GL_n (\O_F)$ (cette situation correspond à  l'espace de Rapoport-Zink des déformations du groupe étale $(F/\O_F)^n$). 
 Cela n'est pas possible ! En effet, pour $g\in \GL_n (F)$ les groupes $g \GL_n (\O_F) g^{-1}$ et $\GL_n (\O_F)$ sont isomorphes
mais {\it non-canoniquement}, cela dépendant du choix de $g$ qui n'est pas canoniquement déterminé. Cependant on peut reconstruire $\GL_n (F)$ à partir de l'immeuble de $\GL_n$ 
 et des cellules $\text{Iso}_{\O_F} (\La,\O_F^n)$ pour $\La$ un réseau de $F^n$,  cellules qui sont des $\GL_n (\O_F)$-torseurs non canoniquement triviaux. Un élément $g\in \GL_n (F)$
induit un isomorphisme entre la cellule indexée par $\La$ et celle 
indexée par $g.\La$, d'où une action de $\GL_n (F)$ sur l'union disjointe de ces cellules.   
 On a alors la décomposition cellulaire $\GL_n (F)$-équivariante 
indexée par l'immeuble 
$$
\coprod_{\La} \text{Iso}_{\O_F} (\La,\O_F^{n}) \iso \GL_n (F)
$$
Pour un sous-groupe compact ouvert $K$ dans $\GL_n (\O_F)$ on peut reconstituer cellulairement une partie de l'espace $\GL_n (F)/K$. Plus précisément soit $A\subset \GL_n (F)/\GL_n (\O_F)$  un sous-ensemble fini
tel que $\forall \La \in A\; K\subset \GL (\La)$
i.e. $K$ stabilise $\La$. On peut pour $\La\in A$ définir une cellule en niveau $K$ : $\text{Iso}_{\O_F} (\La,\O_F^n)/K$ où $K$ agit sur $\La$. Alors, 
$$
\coprod_{\La \in A} \text{Iso}_{\O_F} (\La,\O_F^n)/K \iso \{ gK\in \GL_n (F)/K\;|\; g \GL_n (\O_F) \in A\; \}
$$

\subsection{Structures de niveau}

\begin{defi}
Soit $H$ le groupe p-divisible universel au dessus de $\Mf$. 
Pour $\La$ un réseau de $F^n$ et $K\subset \GL_n (F)$ un sous-groupe compact ouvert tel que $K$ stabilise $\La$ on pose 
$$
\M_{\La,K} = \text{Isom}_{\O} (\pi^{-n}\underline{\La} / \underline{\La},  H[\pi^n]^{rig})/K\;\;
\text{ pour } n>>0 
$$
comme faisceau étale quotient au dessus de $\M$.  Il est représenté par un
espace rigide étale fini au dessus de $\M$.
\end{defi}

Soit $U$ un espace rigide quasicompact. 
On notera $(I,\rho,\eta)$ pour une section de $\M_{\La,K}$ sur $U$. Cela signifie que
l'on se donne un modèle entier $\mathcal{U}$ de $U$
 $$\mathcal{U}^{rig}\iso U$$
puis une section 
$$
(I,\rho)\in \Mf (\mathcal{U})
$$
qui définit un élément de $s\in\M (U)$ et enfin une section $\eta$ du produit fibré du diagramme suivant 
$$
\xymatrix{
\M_{\La,K} \ar[r] & \M \\
 & U \ar[u]^s
}
$$ 
que l'on notera parfois en
$$
\eta : \La \iso T_p (I) \;[K]
$$
une structure de niveau $K$ sur la fibre générique de $I$. 

\begin{exem}
L'espace noté $\M_K$ dans \cite{RZ} avec $K\subset \GL_n (\O_F)$ n'est rien d'autre que $\M_{\O^n,K}$. Bien s\^ur 
$\M_{\La,K}\simeq \M_{K'}$ pour un $K'\subset \GL_n (\O_F)$, mais {\bf non-canoniquement}. 
\end{exem}

Lorsque $\La$ est fixé et
 $K$ varie on obtient ainsi une tour d'espaces rigides dont les morphismes de transition sont étales finis et qui est munie d'une action de $\GL (\La)\times D^\times$. 

\subsection{Fonctorialité de Hecke des $\M_{\La,K}$ }
\subsubsection{Première fonctorialité} \label{pretyfo}

Il y a des isomorphismes canoniques 
$$
\forall g\in \GL_n (F) \; \forall \La \; \forall K\subset \GL (\La)\;\;\;
g:\M_{\La,K}\iso \M_{g^{-1} \La, g^{-1} K g}
$$
définis de façon modulaire par 
$$
(I, \rho,\eta) \longmapsto (I,\rho,\eta\circ g)
$$

\subsubsection{Seconde fonctorialité}\label{secft}
\paragraph{Cas général : }

Il s'agit d'isomorphismes 
$$
\forall \La,\La'\;\; \forall K\subset \GL(\La)\cap \GL (\La')\;\;\;
\M_{\La, K } \iso \M_{\La',K}
$$
définis de la façon suivante. Soit $U$ un espace rigide quasicompact et $(I,\rho,\eta)\in \M_{\La,K} (U)$ où $U=\mathcal{U}^{rig}$ et $(I,\rho)\in \Mf (\mathcal{U})$. Soient $n\in\N$ et $N\in \Z$  tels que 
$$
\La \subset  \pi^N\La' \subset \pi^{-n} \La
$$
La structure de niveau $\eta$ induit
$$
\pi^{-n}\La/\La \iso I[\pi^n] \text{ mod } K
$$
(on entend par là une section du faisceau étale quotient $\text{Isom}_\O ( \pi^{-n}\underline{\La}/\underline{\La}, I[\pi^n]^{rig} )/K$). Le sous-$\O$-module $\pi^N \La'/\La\subset \pi^{-n} \La/\La $ est stable sous l'action de $K$. Il existe donc un unique sous-espace rigide en groupes fini localement libre sur $U$, 
$$
J\subset I[\pi^n]^{rig}
$$
tel que le diagramme suivant commute
$$
\xymatrix{
\pi^N \La'/\La \ar[r]^\sim \ar@{^(->}[d] & J  \ar@{^(->}[d]   & \hspace{-8mm} \text{mod } K    \\
\pi^{-n} \La /\La \ar[r]^\sim & I[\pi^n]^{rig} & \hspace{-8mm} \text{mod } K
}
$$
D'après la section 5.4 de \cite{Ray2} il existe un éclatement formel admissible $\widetilde{\mathcal{U}}\ldrt \mathcal{U}$
donné par certains idéaux de Fittings, induisant donc un isomorphisme $\widetilde{\mathcal{U}}^{rig}\iso \mathcal{U}^{rig}$, tel que si $\widetilde{J}$ désigne l'adhérence schématique de $J$ dans 
$I[\pi^n]\times_{\mathcal{U}} \widetilde{\mathcal{U}}$ alors $\widetilde{J}$ soit un sous-groupe fini localement libre de $I[\pi^n]\times_{\mathcal{U}} \widetilde{\mathcal{U}}$. Par adhérence schématique on entend la chose suivante : si $\mathcal{I}\subset \O_{I[\pi^n]^{rig}}$ est l'idéal cohérent définissant $J$ comme fermé dans $I[\pi^n]^{rig}$, si 
$$
sp: I[\pi^n]^{rig} \ldrt I[\pi^n]\times_{\mathcal{U}} \widetilde{\mathcal{U}}
$$
est le morphisme de spécialisation associé au modèle formel $I[\pi^n]\times_{\mathcal{U}} \widetilde{\mathcal{U}}$ alors 
$$
\widetilde{J}= \spf \left ( \O_{I[\pi^n]\times_{\mathcal{U}} \widetilde{\mathcal{U}}} / (
sp_* \mathcal{I}\cap  \O_{I[\pi^n]\times_{\mathcal{U}}\widetilde{\mathcal{U}}} )
\right )
$$
Alors, le morphisme $\M_{\La,K}\ldrt \M_{\La',K}$ est défini par
$$
(I,\rho,\eta) \longmapsto (I\times_{\mathcal{U}} \widetilde{\mathcal{U}}/\widetilde{J}, q \circ \widetilde{\rho} \circ \pi^{-N} , \eta')
$$
où $\widetilde{\rho}$ est le changement de base de $\rho$ sur $\mathcal{U}$ à $\widetilde{\mathcal{U}}$, 
$$
q : I\times_{\mathcal{U}} \widetilde{\mathcal{U}} \twoheadrightarrow I\times_{\mathcal{U}} \widetilde{\mathcal{U}}/\widetilde{J}
$$
et $\eta'$ fait commuter le diagramme suivant 
$$
\xymatrix{
 \La ' \ar@{^(->}[d]_{\pi^{-N}}  \ar[r]_\sim^{\eta'} & T_p (I/J)   \ar@{^(->}[d]& \hspace{-5mm}  \text{mod } K  \\
   \La \ar[r]_\sim^\eta & T_p (I)  & \hspace{-5mm} \text{mod } K
}
$$

\begin{exem}
Si $K$ et $g^{-1} K g$ sont contenus dans $\GL_n (\O)$ le morphisme 
classique (section 5.43 de \cite{RZ}, section 2.3.9.3 page 39 de \cite{Laurent1}) définissant les correspondances de Hecke
$$
\M_K \xrig{\; g \;} \M_{g^{-1} K g} 
$$
n'est rien d'autre que le composé des deux fonctorialités précédentes
$$
\M_K=\M_{\La_0,K}\xrig{\; g\;} \M_{g^{-1}\La_0, g^{-1} K g} \xrig{2^{\text{éme}}\text{ fonctorialité}} \M_{\La_0, g^{-1} K g} = \M_{g^{-1} K g}
$$
où on a posé $\La_0= \O^n$.
\end{exem}

\paragraph{Cas des sous-groupes de congruence principaux : }

Dans le paragraphe précédent faisons de plus l'hypothèse que :
$$
\exists n\in \N \;\exists N\in\Z\;\;\;\; K\subset \text{Id} + \pi^n \End (\La) \text{ et } \La \subset \pi^N\La' \subset \pi^{-n} \La
$$
Alors, dans la définition du morphisme 
$$
\M_{\La,K} \ldrt \M_{\La',K}
$$
il n'est pas nécessaire d'effectuer l'éclatement $\widetilde{\mathcal{U}}\ldrt \mathcal{U}$. En effet, $\M_{\text{Id} + \pi^n \End (\La)}$ possède un modèle entier défini en utilisant des structures de niveau de Drinfeld (cf. la section II.2 de \cite{Har4}). On conclu alors grâce au lemme clef suivant :
\begin{lemm}[lemme II.2.4 de \cite{Har4}]\label{letHT}
Soit $H$ un $\O$-module formel de dimension $1$  et hauteur $h$ sur un $\O$-schéma $S$ muni d'une structure de niveau de Drinfeld 
$$
\eta: (\pi^{-n} \O/\O)^h \ldrt H[\pi^n](S)
$$
Soit $M\subset  (\pi^{-n} \O/\O)^h$ un sous-$\O$-module. Il existe
alors un unique sous-groupe  fini localement libre $G\subset H[\pi^n]$
tel que $\forall m\in M \;\eta (m) \in G(S)$ et les $(\eta (m))_{m\in
  M}$ forment un ensemble plein de sections de $G$ au sens de Katz-Mazur.
\end{lemm}

\begin{rema}
Le cas des sous-groupes de congruence principaux est suffisant pour définir le morphisme $\M_{\La,K}\ldrt \M_{\La',K}$ en général. En effet, $\forall K\subset \GL (\La)\cap \GL (\La')\;\exists K_1\lhd K$ vérifiant les hypothèses de ce paragraphe. On peut donc définir le morphisme $\M_{\La,K_1}\ldrt \M_{\La',K_1}$ comme expliqué ci-dessus en utilisant les structures de niveau de Drinfeld, vérifier que ce morphisme est $K$-équivariant puis définir $\M_{\La,K}\ldrt \M_{\La',K}$ comme 
$$
\M_{\La,K_1}/K \ldrt \M_{\La',K_1}/K
$$
\end{rema}

\subsection{Les cellules}

\begin{defi}\label{dfgty}
Soit $[\La, M]\in \mathcal{I}$ et $K\subset \GL_n (F)$ 
un sous-groupe compact ouvert
tel que $K$ stabilise $\La$. On pose 
\begin{eqnarray*}
\mathcal{D}_{[\La,M],K} &=& \text{ domaine fondamental de
 Gross-Hopkins dans la fibre générique de 
 l'espace} \\
&&\text{des déformations } (H,\rho) \text{ 
 de } \mathbb{H} \text{ par des quasi-isogénies de hauteur } [M : \O_D]
 \\
& & \; + \text{ une structure de niveau } K, \eta: \La\iso T_p H\; [K]
\end{eqnarray*}
Plus précisément, soit $k=[M : \O_D ]$ et 
$$
\mathcal{D}^{[k]} 
$$
l'ouvert admissible quasicompact dans $\M^{[k]}$ où le polygone de
Newton
de la multiplication par $\pi$ sur la loi de groupe formelle
universelle
est au dessus de celui de la figure  \ref{NewGH}. Alors
$\mathcal{D}_{[\La,M],K}$ est défini par le diagramme cartésien
suivant 
$$
\xymatrix{
\mathcal{D}_{[\La,M],K} \ar@{^(->}[r] \ar[d] & \M^{[k]}_{\La, K}
\ar[d] \\
\mathcal{D}^{[k]}  \ar@{^(->}[r] & \M^{[k]} 
}
$$
\end{defi}

\begin{rema}
Si $K\subset \GL (\La)$ et $g\in \GL_n (F)$ sont tels que $g^{-1} \La= \O_F^n$
  il y a alors un isomorphisme 
$$
g\times\Pi^{[M : \O_D]} : \mathcal{D}_{[\La,M],K} \iso \mathcal{D}_{g^{-1} K g} \subset \M_{g^{-1} K g}^{[0]} 
$$
où $\mathcal{D}_{g^{-1} K g}$ est le revêtement étale fini au dessus 
de $\mathcal{D}$ défini en mettant des
structures de niveau $g^{-1} K g$. Cependant cet isomorphisme 
{\bf n'est pas canonique } puisqu'il dépend du choix de $g$. 
\end{rema}

Le premier type de fonctorialité définie dans la section \ref{pretyfo}
induit une action naturelle de 
 $\GL_n (F)\times D^\times$  sur les 
cellules de façon compatible à son action sur $\mathcal{I}$ : 
$$
\forall (g,d)\in \GL_n (F) \times D^\times\;\;\;
g\times d : \mathcal{D}_{[\La,M],K} \iso  \mathcal{D}_{[g^{-1}\La,d.M],g^{-1}Kg}
$$

\subsection{Bord des cellules}

\begin{defi}
Soient $[\La,M]\in \mathcal{I}$ et $K$ tel que $K \subset \text{Id}+ \pi\End_{\O_F} (\La)$. On définit pour $1\leq i\leq n-1$ 
 les ouverts admissibles quasicompacts
(des domaines de Laurent) $\partial_i \mathcal{D}_{[\La,M],K}$ 
 de $ \mathcal{D}_{[\La,M],K}$ par le diagramme cartésien suivant 
$$
\xymatrix{
\partial_i \mathcal{D}_{[\La,M],K} \ar@{^(->}[r] \ar[d] & 
 \mathcal{D}_{[\La,M],K} \ar[d] \\
\partial_i \mathcal{D}_{[\La,M],GL(\La)} \ar@{^(->}[r] &
 \mathcal{D}_{[\La,M],GL (\La)}
}
$$
où $\mathcal{D}_{[\La,M],GL (\La)}$ est la cellule sans structures de
niveau notée $\mathcal{D}^{[[ M : \O_D ]]}$ dans la définition
\ref{dfgty} 
 et l'ouvert $\partial_i$  est défini de façon modulaire
en termes du polygone du Newton de la série formelle $[\pi]_{F^{univ}}$ comme dans
la définition \ref{borde}.
\end{defi}

\begin{lemm}
Soient $\mathfrak{Z}$ un schéma formel admissible sur $\spf
(\breve{\O})$ et $H$ un $\O$-module $\pi$-divisible formel de
dimension $1$ et hauteur $n$. Supposons $H$ muni d'une structure de
niveau de Drinfeld
$$
\eta : \pi^{-1}\La/\La \ldrt H[\pi]
$$
Soit $i$ un entier tel que $1\leq i\leq n-1$. 
Supposons que $\forall x\in\mathfrak{Z}^{an}$ le polygone de Newton
de la multiplication par $\pi$ sur une loi de groupe formel associée à
$H_x$ possède un point de rupture en $q^i$. 
La donnée $\eta$ induit
un scindage de l'espace rigide
$$
\mathfrak{Z}^{rig} = \coprod_{E\subset \pi^{-1}\La/\La \atop
  \dim_{\Fq} E =i} (\mathfrak{Z}^{rig})_E
$$
où $(\mathfrak{Z}^{rig})_E= \{ x\in \mathfrak{Z}^{rig} \;|\; \eta_x
(E) \subset H_x [\pi]\text{ est un sous-groupe canonique de rang } i\}$.
\end{lemm}
\dem Soit $\mathfrak{Z}^{an}$ l'espace de Berkovich fibre générique
associé à $\mathfrak{Z}$. On a 
$$
(\mathfrak{Z}^{an})_E= \{x\in |\mathfrak{Z}^{an}|\;|\; \forall w\in
E\; \forall w'\in \pi^{-1}\La/\La\setminus E\;\; v(\eta_x
(x))>v(\eta_x (x')) \}
$$
qui définit bien un ouvert de $\mathfrak{Z}^{an}$ car localement sur
$\mathfrak{Z}$, si l'on fixe une loi de groupe formel $\spf
(\O_{\mathfrak{Z}}[[T]])\iso \widehat{H}$ alors $\forall w\in
\pi^{-1}\La/\La\; \;\eta (w)\in \O_{\mathfrak{Z}}$ et la valuation de
$\eta (w)$ est donnée par la valuation de cet élément de $\O_Z$.
\qed

\begin{rema}
Le lemme précédent est une version rigide de ``l'astuce de Boyer''
(cf. \cite{Har4}). 
\end{rema}

\begin{prop}
Les bords des cellules se scindent de la façon suivante 
$$
\partial_i \mathcal{D}_{[\La,M],K} = \coprod_{E\subset \pi^{-1} \La/\La \atop rg\, E=i} \partial_{i,E} \mathcal{D}_{[\La,M],K} 
$$
où $E$ est un sous-$\Fq$-ev. de $\pi^{-1} \La/\La$ de rang $i$ et 
$ \partial_{i,E} \mathcal{D}_{[\La,M],K}$ est l'ouvert admissible de 
$\partial_i \mathcal{D}_{[\La,M],K}$ où $\eta (E) \subset H [\pi]_\eta$ est l'ensemble des $q^i$-points de plus grande valuation (un sous-groupe canonique généralisé). 
\end{prop}
\dem
Cela résulte de ce que sur $\partial_i\mathcal{D}$  le polygone de Newton de $[\pi]_{F^{univ}}$ 
possède un point de rupture en $q^i$ et du lemme précédent.
\qed

\begin{rema}
Bien s\^ur cette décomposition en niveau $K$ est obtenue par image réciproque 
de celle en niveau $\text{Id}+\pi \End (\La)$. 
\end{rema}

\begin{rema}
On vérifie que tout cela ne dépend que de $[\La,M]$ i.e. il y a des isomorphismes canoniques $\partial_{i,E}\mathcal{D}_{[\La,M],K} \iso
\partial_{i,\pi E}\mathcal{D}_{[\pi \La,\pi^{-1} M],K}$.
\end{rema}

L'action de $\GL_n (F)\times D^\times$ conserve le bord des 
cellules en permutant les composantes indexées par les sous-espaces $E$
: $\forall g\in \GL_n (F)\;\; g
: \partial_{i,E} \drt \partial_{i, g^{-1}.E}$.  

\subsection{Donnée de recollement}\label{donrec}

Soit $p: \pi^{-1} \La \twoheadrightarrow \pi^{-1} \La /\La$.
Supposons de plus  que $K\subset \text{Id} + \pi \End (p^{-1} (E))$. 
Le quotient par le sous-groupe canonique généralisé
$\eta (E)$  induit alors une immersion ouverte 
$$
\partial_{i,E} \mathcal{D}_{[\La,M],K} \hookrightarrow 
\mathcal{D}_{[p^{-1} (E), \Pi^{-i} M],K}
$$
où le $\Pi^{-i}$ provient de ce que le quotient par le sous-groupe canonique est une
déformation de $Frob_q^i$. Ce morphisme est induit par le second type
de fonctorialité de la section \ref{secft} restreint à $\partial_{i,E}$.

\subsection{Réécriture en termes des arr\^etes orientées de l'immeuble}
\label{secprec}

Soit $a\ldrt a'$ une arr\^ete orientée de $\mathcal{I}$ avec 
$a=[\La,M]$ et $a'=[\La',M']$ tels que 
$$
\La \varsubsetneq \La' \varsubsetneq \pi^{-1} \La \text{ et } 
M'=\Pi^{[\La: \La']} M 
$$
\begin{defi}
Supposons que $K\subset \text{Id} + \pi \End (\La) \cap \text{Id} + \pi \End (\La')$. 
On pose 
$$
\mathcal{D}_{a\drt a',K}  = \partial_{i,\La'/\La} \mathcal{D}_{[\La,M],K}
$$
où $i=\text{dim}_{\Fq} \La'/\La$
\end{defi}

Il y a alors deux immersions ouvertes
$$
\xymatrix@C=6mm@R=6mm{
& \mathcal{D}_{a\drt a',K}
 \ar@{^(->}[ld] \ar@{_(->}[rd] \\ 
 \mathcal{D}_{a,K} & & \mathcal{D}_{a',K} 
}
$$
où l'application de gauche est l'inclusion canonique et celle
de droite est celle définie en \ref{donrec}. 

\begin{rema}
On peut montrer que la seconde application induit un isomorphisme 
$$
\mathcal{D}_{a\drt a',K}
 \iso \mathcal{D}_{a'\drt a,K}
$$
(cf. la section sur les modèles entiers où tout cela est démontré
plus généralement sur les modèles entiers). 
\end{rema}

Bien s\^ur toutes ces applications sont équivariantes sous l'action de
$\GL_n(F)\times D^\times$ au sens où il y a des isomorphismes
se composant naturellement
$$ 
\forall (g,d)\in \GL_n (F)\times D^\times \;\; \mathcal{D}_{a\drt a',K}
 \iso \mathcal{D}_{(g,d).a \drt
  (g,d).a',g^{-1} K g} 
$$
et compatibles aux deux immersions ouvertes ci-dessus.

\section{Décomposition cellulaire des espaces rigides en niveau fini}\label{decellrig}

Soit $A\subset \mathcal{I}$ un sous ensemble d'image finie 
dans $\mathcal{I} (\text{PGL}_n)= \mathcal{I} (\GL_n)/\pi^\Z$ et d'image $\mathcal{I} (D^\times)/\pi^\Z$ ($\simeq \Z/n\Z$)
sur la seconde composante.
Soit $K\subset \GL_n (F)$ un sous-groupe compact ouvert tel que
$\forall [\La,M]\in A$ le groupe  $K$ stabilise $\La$ et $K\subset \text{Id}+ \pi
\End (\La)$ (l'image de $A$ 
 dans $\mathcal{I} (\text{PGL}_n)$ 
étant finie il existe toujours un tel $K$). 

Considérons le diagramme 
d'espaces rigides fibré au dessus de l'immeuble $\mathcal{I}$ :
$$
\xymatrix{ \displaystyle{
X_{1,A,K}= \coprod_{a,a'\in A \atop a\drt a'}   
\mathcal{D}_{a\drt a',K}}
 \ar@<2ex>[r]\ar@<0ex>[r]  & \displaystyle{\coprod_{a\in A\atop \,}   \mathcal{D}_{a,K} = X_{0,A,K} }}
$$
défini par les deux applications de faces de 
\ref{secprec}. Celui-ci est
$\GL_n (F) \times D^\times$ équivariant pour des $A$ et $K$ variant
(le diagramme associé à $A$ et $K$ est envoyé naturellement sur celui
associé à $(g,d).A$ et $g^{-1} K g$).  

Rappelons que l'on note $\M_{K} : = \M_{\O^n,K}$ l'espace de Rapoport-Zink usuel. 
L'application naturelle définie par le second type de fonctorialité 
$\M_{\La,K}\iso \M_K$ induit un morphisme 
$X_{0,A,K}\ldrt \M_K$. 
  
\begin{prop}\label{rigpi}
L'image de $X_{0,A,K}$ dans $\M_K$ est un ouvert quasicompact
$U_{A,K}$ et la suivante est 
exacte
dans la catégorie des espaces rigides
$$ 
\xymatrix{
X_{1,A,K} \ar@<1ex>[r]\ar@<-1ex>[r]& X_{0,A,K} \ar[r] & U_{A,K} 
}
$$
Ces suites sont $\GL_n (F) \times D^\times$-équivariantes lorsque $A$
et $K$ varient. 
De plus lorsque $A$ grandi et $K$ est de plus en plus petit les $U_{A,K}$ recouvrent toute la tour de Lubin-Tate au sens où les deux ind-pro systèmes d'espaces rigides $(U_{A,K})_{A,K}$ et celui des ouverts de $\M_K$ d'image quasicompacte dans $\M_K/\pi^\Z$ pour $K$ variant 
sont équivalents. 
\end{prop}
\dem
Les applications  étant toutes des unions disjointes d'immersions ouvertes d'espaces quasi-compacts 
il suffit de vérifier ces assertions au niveau des points ce qui résulte du corollaire 23.26 de \cite{HopkinsGross} et de  la proposition \ref{pofgh}.
\qed

\begin{exem}
Si $A$ est  l'image réciproque d'un sous-ensemble fini $B$ de 
$\mathcal{I} (\text{PGL}_n)$  alors l'ensemble $B$ définit un
sous-ensemble de correspondances de Hecke sphériques de $\GL_n (F)$ et
$U_{A,K}$ est l'image réciproque dans $\M_K$ de l'itération par ces
correspondances sphériques de l'ouvert $\mathcal{D} \subset \M^{[0]}$, un ouvert de $\M$. 
\end{exem}

\begin{rema}\label{re6}
Bien s\^ur on peut écrire les relations d'incidence supérieures en
termes de facettes orientées $a_1\drt \dots \drt a_d$ de l'immeuble 
et de cellules $ \mathcal{D}_{a_1\drt \dots \drt a_d,K}$. 
\end{rema}

\subsection*{Application Cohomologique}

Dans cette section la cohomologie étale à support compact des espaces rigides
désigne celle définie par Huber. Pour $\M_K$ elle coïncide avec celle définie par Berkovich. Cependant ce n'est pas le cas pour les $U_{A,K}, \mathcal{D}_{a,K}$... 

\begin{coro} Pour $\La$ un anneau de $\ell$-torsion ou bien $\La=\Zl$
il y a des isomorphismes 
dans la catégorie dérivée $\GL_n (F)\times D^\times \times W_F$-équivariante-lisse
$$
\underset{K}{\limi} R\GG_c (\M_K,\La) \simeq \underset{A,K}{\limi}
R\GG_c (U_{A,K},\La)
$$
$$
\underset{K}{\limi} R\GG_c (\M_K,\La) \simeq  \underset{A,K}{\limi}
R \GG_c 
( X_{\bullet,A,K}, \La)
$$
où $X_{\bullet,A,K}$ désigne le diagramme d'espaces rigides de type  le diagramme des facettes orientées $a_1\drt \dots \drt a_d$ défini précédemment (remarque \ref{re6})(cf. thèse Illusie tome II pour la cohomologie des diagrammes de topos, ici la cohomologie à support compact
désignant la cohomologie de $j_!\La$ où $j$ désigne l'inclusion de
$X_{\bullet,A,K}$ dans le diagramme compactifié universel d'espaces adiques
construit par Huber (\cite{Hu1} théorème 5.1.5). 
\end{coro}

\begin{rema}
La deuxième égalité dans le corollaire précédent fournit 
des résolutions équivariantes de la cohomologie 
à support compact de la tour de
Lubin-Tate
par des induites compactes de la cohomologie des cellules et de leur
bord. 
\end{rema}

\section{Modèles entiers des cellules}
\subsection{Niveau fini}
\subsubsection{Préliminaires}

Rappelons qu'un $\breve{\O}$-schéma formel admissible est 
un  $\breve{\O}$-schéma formel  quasi-séparé de type fini sans $\pi$-torsion.

\begin{prop}
Soient $\X$ un  $\breve{\O}$-schéma formel admissible normal (cf. 
appendice \ref{schfnor}) et $G$ un groupe $p$-divisible sur $\X$ muni
d'une action de $\O$. Soient $\La$ un $\O$-module libre de rang $n$,
$X_K\ldrt \X^{rig}$ l'espace classifiant des structures de niveau $K$ 
$$
\eta : \La \iso T_p (G^{rig})\; [K]
$$
et $\X_K$ la normalisation de $\X$ dans $X_K$ (cf. appendice
\ref{norextgen}). Alors, $\X_K$ représente le foncteur $F$ défini sur 
la catégorie des $\breve{\O}$-schémas formels admissibles
normaux au dessus de
$\X$ défini par $\forall \mathfrak{Z}\xrig{\; f\;} \X$
$$
F(\mathfrak{Z}) = \{ \text{ structures de niveau } K \text{ rel. à }
\La \text{ sur } (f^* G)^{rig}\; /\mathfrak{Z}^{rig} \; \}
$$
\end{prop}
\dem 
C'est une conséquence de la propriété universelle du normalisé,
cf. l'appendice \ref{norextgen}. 
\qed

\begin{rema}\label{stDr}
Si $G$ est un groupe de Lubin-Tate et $K\subset \text{Id} + \pi^k \End
(\La )$ alors sur $\X_K$ le groupe  $G$ est muni d'une structure de niveau de
Drinfeld de niveau $k$. C'est une conséquence du fait que l'espace 
classifiant des structures de Drinfeld de niveau $k$ est fini au
dessus de $\X$, a même fibre générique que $\X_K$ (en fibre générique
toutes les définitions des structures de niveau coïncident)  
et est donc en dessous du normalisé.
\end{rema}

\begin{prop}\label{ptui}
Soit $\mathcal{P}$ un polygone de Newton i.e. la donnée pour $i$ entre
$1$ et $n-1$ de nombres rationnels $\a_i$, $0< \a_i <1$, tels que
le polygone commençant en $(0,1)$, passant par les $(q^i,\a_i)$ et
finissant en $(q^n, 0)$ soit convexe. Le foncteur défini sur la
catégorie des $\breve{\O}$-schémas formels admissibles normaux qui à
$\mathfrak{Z}$ associe l'ensemble des classes d'isomorphismes de
couples $(H,\rho)$ où $H$ est un $\O$-module formel et
$$
\rho : \mathbb{H}\times_{\Fqb} \mathfrak{Z} \text{ mod } p \ldrt 
 H  \text{ mod } p
$$ 
une quasi-isogénie de degré $0$ tels que 
$$\forall z\in\mathfrak{Z}^{rig}\;\;\;
\text{Newt} ( H[\pi]_z) \geq \mathcal{P}
$$
est représentable.
\end{prop}
\dem
Soit $\X = \spf (\breve{\O} [[x_1,\dots, x_{n-1}]])
$ l'espace de Lubin-Tate. Pour $\xs \in \X^{rig}$
$$
\text{Newt} ( [\pi]_{F_{\xs}}) \geq\mathcal{P} \ssi
\forall i \; v(x_i)\geq \a_i
$$ 
Si $\dpt{\a_i = \frac{a_i}{b_i}}$ où $a_i,b_i \in \N$ soit 
$$
\mathfrak{Y}= \spf \left (
\breve{\O} <x_1,\dots, x_{n-1}, T_1,\dots, T_{n-1}>/(x_i^{b_i} -
\pi^{a_i} T_i )_i \right )
$$
Alors $\mathfrak{Y}^{\text{normalisé}}$ convient. En effet, 
si $\mathfrak{Z}$ est normal et $(H,\rho)$ est définie sur
$\mathfrak{Z}$ $\;\exists ! f : \mathfrak{Z}\ldrt \X$ tel que 
$$
(H,\rho)= f^*(H^{univ}, \rho^{univ})
$$
Si de plus $\forall z\in\mathfrak{Z}^{rig}\;\;\text{Newt} ( H[\pi]_z) \geq \mathcal{P}$ alors 
$$
\forall i\;\forall z \in \mathfrak{Z}^{rig}\;\;
\frac{|f^* x_i^{b_i} (z)|}{\pi^{a_i}} \leq 1 
$$
$$
\limpl 
\frac{f^*x_i^{b_i}}{\pi^{a_i}} \in \O_{\mathfrak{Z}}
$$ 
puisque $\mathfrak{Z}$ est normal
(cf. appendice \ref{schfnor}). Et donc 
$
\mathfrak{Z} \drt \X \text{ se factorise en }$ 
$$
\xymatrix{
 & \mathfrak{Y}^{\text{normalisé}} \ar[d] \\
\mathfrak{Z} \ar@{-->}[ru]^{\exists !}  \ar[r] \ar[rd] & \mathfrak{Y} \ar[d] \\
 & \X
}
$$
où la flèche en pointillés résulte de ce que $\mathfrak{Z}$ est
normal.
\qed

\subsubsection{Modèles entiers des cellules $\mathcal{D}_{[\La,M],K}$}

\begin{prop}
Soit $(\La,M)$ dans l'immeuble de $\GL_n\times D^\times$. Soit
$K\subset \GL (\La)$ un sous-groupe compact ouvert. Le foncteur qui à
un $\breve{\O}$-schéma formel admissible normal  $\mathfrak{Z}$  associe
les classes d'isomorphisme de triplets $(H,\rho,\eta)$ où 
$H$ est un $\O$-module formel sur $\mathfrak{Z}$, $\rho$ une
rigidification de degré $[M:\O_D]$ et $\eta$ une structure de niveau
$K$ relativement à $\La$ sur $T_p (H^{rig})$ tels que
$$
\forall z\in \mathfrak{Z}^{rig}\;\; \text{Newt} (H[\pi]_z)\geq 
\text{ le polyogne de Gross-Hopkins}
$$
est représentable. De plus il ne dépend canoniquement que de la classe $[\La,M]$
dans $\mathcal{I}$.
\end{prop}
\dem 
C'est une conséquence des deux propositions précédentes, puisque
quitte à translater par une puissance de $\Pi$ on peut supposer que
la rigidification est de degré $0$. 

Quant à la dernière assertion, il suffit de constater que l'application
naturelle $(H,\rho,\eta)\mapsto (H/H[\pi],\rho\circ h,h_*\circ \eta)$
où $h : H\twoheadrightarrow H/H[\pi]$ induit un isomorphisme canonique 
entre le foncteur associé à $(\La,M)$ et celui associé à $(\pi^{-1}
\La, \pi M)$. 
\qed 

\begin{defi}
On note $\mathbb{D}_{[\La, M],K}$ le $\breve{\O}$-schéma formel
admissible normal défini dans la proposition précédente.
\end{defi}

Il y a des isomorphismes naturels
 $$\forall 
(g,d)\in  \GL_n (F) \times D^\times 
\;\;\; g\times d : \mathbb{D}_{[\La,M],K}\iso \mathbb{D}_{[g^{-1}\La,
    d.M],g^{-1} K g}
$$
via 
$$
(H,\rho,\eta)\longmapsto (H,\rho \circ d^{-1}, \eta\circ g)
$$

\begin{rema}
On peut calculer explicitement $\mathbb{D}_{[\La,M],GL (\La)}$. Il
s'agit de 
$$
\spf (\breve{\O}<x_1,\dots,x_{n-1},T_1,\dots,T_{n-1}>/(x_i^n - \pi^{n-i} T_i))^{\text{normalisé}}
$$
L'algèbre le définissant est  
 engendrée par l'algèbre 
$\breve{\O}<x_1,\dots,x_{n-1},T_1,\dots,T_{n-1}>/(x_i^n - \pi^{n-i} T_i)$ à laquelle on a rajoutée les
$$
\frac{{x_i}^{\left \lceil \frac{n}{n-i} \right \rceil}}{\pi},
\frac{{x_i}^{\left \lceil 2 \frac{n}{n-i} \right \rceil}}{\pi^2},
\dots,
\frac{{x_i}^{\frac{n}{n\wedge i}}}{\pi^{\frac{n-i}{n\wedge i}}}
$$
Il est de la forme $\overline{T}_\l$ (avec les notations de  \cite{Faltings6}) où $\overline{T}$ est la variété torique formelle (le complété $\pi$-adique d'une variété torique sur $\breve{\O}$) 
$$
\overline{T}=\spf (\breve{\O} <x_i,T_i,Z>/ (x_i^n- Z^{n-i} T_i))^{\text{normalisé }}
$$
avec $\l=Z$, $\overline{T}_\l= V(\l-\pi)$.
\end{rema}

\subsubsection{Bord des cellules}

\begin{prop}
Pour $i$ un entier vérifiant $1\leq i \leq n-1$ on note
$$
\partial_i\mathbb{D}_{[\La,M],K} \hookrightarrow \mathbb{D}_{[\La,M],K}
$$ 
le sous-foncteur de $ \mathbb{D}_{[\La,M],K}$ défini par l'ensemble
des $(H,\rho,\eta)$ tels que pour tout $z\in \mathfrak{Z}^{rig}$ le
polygone 
 $\text{Newt} ( H[\pi]_z)$ passe par le
point $\dpt{(q^i,1-\frac{i}{n})}$. Ce sous-foncteur est un ouvert de 
$  \mathbb{D}_{[\La,M],K}$.
\end{prop}
\dem
Avec les coordonnées explicites choisies sur
les espaces de Lubin-Tate dans la démonstration de la proposition \ref{ptui}
, celui-ci est défini par l'inégalité
$v(x_i)\geq 1-\frac{i}{n}$ sur l'espace rigide 
et est donc l'ouvert $T_i\neq 0$ avec les notations de la proposition \ref{ptui}.
\qed

\begin{rema}
Bien sûr cet ouvert est obtenu par image réciproque de son homologue
en niveau $K=\GL (\La)$. 
\end{rema}

\subsubsection{Décomposition du bord}

\begin{lemm}
Soit $\X$ un schéma formel admissible quasicompact normal tel que
$$
\X^{rig} = U_1 \coprod U_2
$$
Il existe alors des modèles entiers $\mathcal{U}_1, \mathcal{U}_2$ de $U_1$ et $U_2$ tels que 
$$
\X = \mathcal{U}_1 \coprod \mathcal{U}_2
$$
\end{lemm}
\dem
Le schéma formel $\X$ étant admissible 
$$
\GG ( \X^{rig}, \O_{\X^{rig}}) = \GG (\X,\O_{\X})[\frac{1}{\pi}]
$$
et donc la fonction rigide  valant $0$ sur $U_1$ et $1$ sur $U_2$
définit un élément $e$ de $\GG (\X,\O_{\X})[\frac{1}{\pi}]$ vérifiant 
$e^2=e$. Mais $\X$ étant normal $\GG (\X,\O_{\X})$ est intégralement
fermé dans $ \GG (\X,\O_{\X})[\frac{1}{\pi}]$ (Fait \ref{oktukok}
de l'appendice \ref{Normschf}). Donc $e\in   \GG (\X,\O_{\X})$.
\qed

\begin{coro}
Si $K\subset \text{Id}+ \pi \End (\La)$ 
 il y a une décomposition
$$
\partial_i\mathbb{D}_{[\La,M],K} = \coprod_{E\subset \pi^{-1} \La/\La}
\partial_{i,E} \mathbb{D}_{[\La,M],K}
$$
où $\partial_{i,E} \mathbb{D}_{[\La,M],K}$ représente les
$(H,\rho,\eta)$ dans $\partial_i\mathbb{D}_{[\La,M],K}$ tels qu'en tout
point de la fibre générique $\eta (E)$ soit le sous-groupe des
$q^i$-points de plus grande valuation dans $H[\pi]$.
\end{coro}

\begin{rema}
Plus généralement, soit $\underline{i}= ( 0< i_1 <\dots < i_r < n )$
et $\partial_{\underline{i}}\mathbb{D}_{[\La,M],K} =
\dpt{\bigcap_{a=1}^r \partial_{i_a} \mathbb{D}_{[\La,M],K}}$. Alors, 
$$
\partial_{\underline{i}}\mathbb{D}_{[\La,M],K} = \coprod_{E^\bullet}
\partial_{\underline{i},E^\bullet} \mathbb{D}_{[\La,M],K}
$$
où $E^\bullet$ parcourt les drapeaux de $\Fq$-e.v. dans $\pi^{-1}
\La/\La$ de ``type'' $\underline{i}$ (i.e. $E^1\subset \dots\subset
E^r$ avec $\dim_{\Fq} E^a= i_a$) qui correspondent via $\eta$ à des
drapeaux de sous-groupes canoniques dans $H[\pi]$. 
\end{rema}

\subsubsection{Applications de recollement}

\begin{prop}\label{quca}
Soient $[\La,M]\in \mathcal{I}$ et $E\subset \pi^{-1} \La/\La$ de
dimension $i$. 
Si $p: \pi^{-1}\La \twoheadrightarrow \pi^{-1} \La/\La$, si
$$K\subset \left ( \text{Id} + \pi \End (\La)\right ) \cap  \left (
\text{Id} + \pi \End (p^{-1} (E))\right )
$$
il y a un isomorphisme 
$$
\partial_{i,E} \mathbb{D}_{[\La,M],K} \iso \partial_{n-i,(\pi^{-1}
  \La/\La)/E}
\mathbb{D}_{[p^{-1} (E), \Pi^i M],K}
$$
induit par le quotient par $\eta (E)$.
\end{prop}
\dem 
D'après la remarque \ref{stDr}, sur $\mathbb{D}_{[\La,M],K}$ le
$\O$-module formel universel $H$ possède une structure de niveau de
Drinfeld de niveau $1$ étendant la structure de niveau $\eta$ sur la
fibre générique. Il résulte alors du lemme \ref{letHT} qu'il existe un unique
sous-groupe plat fini $\eta (E)\subset H[\pi]$ induisant
ponctuellement sur $\partial_{i,E}$ 
en chaque point de la fibre  générique le sous-groupe des $q^i$-points
de plus grande valuation dans $H[\pi]$.

D'après les calculs effectués dans la 
 première partie de la démonstration de la proposition
\ref{pofgh} on vérifie que $H/\eta (E)$ est dans $\partial_{n-i}$
et qu'en chaque point de la fibgre générique 
 l'image par l'isogénie $h: H\twoheadrightarrow H/\eta (E)$
envoie un facteur directe de $\eta (E)$ dans $H[\pi]$ sur les 
$q^{n-i}$-points de plus grande valuation.
 On
 en déduit aussitôt que $(H/\eta (E), h\circ \rho , h_*\circ\eta)$
 définit un élément de $\partial_{n-i,(\pi^{-1}
  \La/\La)/E}
\mathbb{D}_{[p^{-1} (E), \Pi^i M],K}$. Cela définit le morphisme 
 $$
\partial_{i,E} \mathbb{D}_{[\La,M],K} \ldrt \partial_{n-i,(\pi^{-1}
  \La/\La)/E}
\mathbb{D}_{[p^{-1} (E), \Pi^i M],K}
$$
Mais en remplaçant $E$ par $(\pi^{-1}
  \La/\La)/E$ on obtient un morphisme dans l'autre sens. La composée
  des deux est le quotient par $H[\pi]$ qui est donc l'identité.
\qed

\subsubsection{Réinterprétation en termes des arêtes orientées de $\mathcal{I}$}\label{apinf}

\begin{defi}
Soit $a\drt a'$ une arête de $\mathcal{I}$ où $a=[\La,M],
a'=[\La',M']$ et 
$$
\La \varsubsetneq \La'  \varsubsetneq \pi^{-1} \La\;\;\;
M'=\Pi^{[\La: \La']} M
$$
Soit $K\subset \text{Id}+ \pi\End (\La)\cap  \text{Id}+ \pi\End
(\La')$. 
On note 
$$
\mathbb{D}_{a\drt a',K} = \partial_{i, \La'/\La} \mathbb{D}_{a,K}
$$
\end{defi}

D'après la proposition \ref{quca} il y a un isomorphisme naturel
$$
\mathbb{D}_{a\drt a',K}\iso \mathbb{D}_{a'\drt a,K}
$$
d'où deux immersions ouvertes
$$
\xymatrix@R=6mm{
 & \mathbb{D}_{a,K} \\
\mathbb{D}_{a\drt a',K} \ar@{^(->}[ru]  \ar@{^(->}[rd] \\
 & \mathbb{D}_{a',K}
}
$$
Ces applications de face sont équivariantes pour l'action de $\GL_n
(F) \times D^\times$. 

\begin{rema} 
Plus généralement
soit 
$$
a_0 \drt a_2\drt \dots \drt a_r
$$
un simplexe orienté de $\mathcal{I}$. Mettons le sous la forme
$a_i=[\La_i,M_i]$ où 
$$
\La_0\subset \La_1 \subset \dots \subset \La_r  \subset \pi^{-1} \La
$$
et $M_i= \Pi^{[\La_0 : \La_i]} M_0$. Cela définit un type
$\underline{i}$ avec $i_a =  \dim_{\Fq} \La_i/\La_0$ et un drapeau
$E^\bullet$ de type $\underline{i}$ dans $\pi^{-1} \La_0 /\La_0$. 

Posons 
$$
\mathbb{D}_{a_0\drt \dots \drt a_r,K} = 
 \partial_{\underline{i},E^\bullet} \mathbb{D}_{[\La_0,M_0],K}
$$
Alors,
pour tout $i$ compris entre $1$ et $r$ il y a des isomorphismes 
$$
\mathbb{D}_{a_0\drt \dots \drt a_r,K} \iso \mathbb{D}_{a_i \drt
  \dots \drt a_r \drt a_0\drt \dots \drt a_{i-1},K}
$$
qui se composent de façon naturelle. Il y a également des applications
de face naturelles : pour tout simplexe orienté $\s$ et tout sous-simplexe
$\s'\subset \s$ : 
$$
\mathbb{D}_{\s,K} \hookrightarrow \mathbb{D}_{\s',K}
$$
Ces applications sont des immersions ouvertes.
\end{rema}

\subsection{Niveau infini}

\begin{defi}
Soit $a\in \mathcal{I}$. On pose 
$$
\mathbb{D}_{a,\infty} = \underset{K}{\limp} \mathbb{D}_{a,K}
$$
($K$ suffisamment petit) dans la catégorie des schémas formels
$p$-adiques.
\end{defi}

Étant donné que les morphismes de transition $\mathbb{D}_{a,K_2} \drt
\mathbb{D}_{a,K_1}$ pour $K_2\subset K_1$ sont affines une telle
limite
existe. Si $\mathbb{D}_{a,K}=\spf (A_{K})$ alors
$$
\mathbb{D}_{a,\infty} =  \spf \left (  (
    \underset{K}{\limi} A_{K}) \widehat{\;} \right )
$$

\begin{rema}
Si $\dpt{\mathcal{A}_K = A_K [\frac{1}{\pi}]}$ est l'algèbre de Banach affinoïde
$p$-adique munie de sa norme infini $|f|_\infty=\underset{x\in
  \text{Max } (\mathcal{A}_K)}{\sup} |f (x) |$, soit $\mathcal{A} =
\bigcup_K \mathcal{A}_K$ munie de la norme infini 
$$
|f |_\infty = \underset{ x : \mathcal{A} \drt \Qpb}{\sup} |f(x)|
$$
(les morphismes de
transition $sp (\mathcal{A}_K) \ldrt sp (\mathcal{A}_{K'})$ pour $K'
\subset K$ étant finis l'inclusion $\mathcal{A}_K \subset
\mathcal{A}_{K'}$ est isométrique). Soit $\widehat{\mathcal{A}}$ le complété de
$\mathcal{A}$, une algèbre de Banach $p$-adique. Alors, d'après le
fait \ref{sntt} de l'appendice \ref{Normschf}
$$
\mathbb{D}_{a,\infty} = \spf \left ( \text{ boule unité de } \mathcal{A}\;
\right )
$$
\end{rema}

\begin{defi}
Pour $\s$ un simplexe orienté de $\mathcal{I}$, $\s=(a\drt\dots)$ on
définit
$ \mathbb{D}_{\s,\infty}$ comme étant l'ouvert correspondant de
$\mathbb{D}_{a,\infty}$. 
\end{defi}

On a donc $ \mathbb{D}_{\s,\infty}= \underset{K}{\limp}
 \mathbb{D}_{\s,K}$. Ceux-ci sont munis d'une action
$$\forall (g,d)\in \GL_n (F)\times D^\times \;\;\;
\mathbb{D}_{\s,\infty} \iso \mathbb{D}_{(g,d).\s,\infty}
$$
et d'applications de face 
$$
\mathbb{D}_{\s,\infty} \hookrightarrow  \mathbb{D}_{\s',\infty}
$$
pour $\s'\subset \s$ des simplexes orientés.

\subsection{Donnée de descente}

La donnée de descente de Rapoport-Zink permet de définir des données
pour tout simplexe orienté  $\s_0$ pour
tout $K$ (éventuellement $K=\infty$) 
$$
\mathbb{D}_{\s_0,K} \ldrt  \mathbb{D}_{\Pi^{-1}.\s_0,K}^{(\s)}
$$
où $\s$ désigne le Frobenius arithmétique de $\widehat{F^{nr}}|F$. 

\section{Le schéma formel recollé en niveau fini}

Soit $A \subset \mathcal{I} $ un sous-ensemble comme dans la partie
\ref{decellrig}. Posons 
$$
\X_{0,A}= \coprod_{a\in A} \mathbb{D}_{a,K} 
$$
et 
$$
\X_{1,A}= \coprod_{a,a'\in A \atop a\drt a'}
\mathbb{D}_{a\drt a',K }
$$
La section \ref{apinf} permet  de définir un diagramme
$$
\xymatrix{ 
\X_{1,A,K}  \ar@<1ex>[r]\ar@<-1ex>[r]  & \X_{0,A,K} }
$$

\begin{prop}
Le diagramme ci dessus définit une donnée de recollement effective
pour la topologie de Zariski et définit donc un schéma formel
$\X_{A,K}$ localement de type fini sur $\breve{\O}$ tel que $\X_{A,K} /
\pi^\Z$ soit de type fini. Ce schéma formel est un modèle entier de
l'ouvert rigide $U_{A,K}$ de $\M_K$ défini dans la proposition
\ref{rigpi}.  Comme dans le cas rigide tout est équivariant sous 
$\GL_n( F) \times D^\times$ pour des $A$ et $K$ variants. 
\end{prop}
\dem
C'est une conséquence de ce que les applications de bord sont
naturelles au sens où le diagramme suivant commute 
$$
\xymatrix{
  \mathbb{D}_{a \drt a' \drt a'',K} \ar[r] \ar[d] &
  \mathbb{D}_{a\drt a",K} \ar[d] \\ 
\mathbb{D}_{a'\drt a",K} \ar[r] & \mathbb{D}_{a",K}
}
$$
(ce qui n'est rien d'autre que l'égalité $(H/\eta(E''))/(\eta
(E')/\eta (E'')) = H/\eta (E')$). 
\qed

\section{Le schéma formel en niveau infini}

\begin{prop}
Le diagramme $\GL_n (F)\times D^\times$-équivariant 
$$\xymatrix{
\dpt{\coprod_{a\drt a'}  \mathbb{D}_{a \drt a',\infty }}
\ar@<2ex>[r] \ar@<0ex>[r] & \dpt{\coprod_a \mathbb{D}_{a,\infty}}
}
$$
définit un schéma formel $p$-adique recollé $\X_\infty$ muni d'une
action de $\GL_n (F)\times D^\times$. De plus
$$
\X_\infty = \underset{A}{\limi} \underset{K}{\limp} \X_{A,K} 
$$
Ce schéma formel est muni d'une donnée de descente de
$\breve{\O}$ à $\O$ qui est effective sur les quotients 
$\X_{\infty}/\pi^{a\Z}$ pour $a\in \N^*$.  
\end{prop}
\dem
Elle ne pose pas de problème.
\qed

\begin{rema}
Soit $\ell\neq p$ et $\La\in \{ \Z/\ell^n \Z, \Zl \}$. Soit $R\GG_c (-,\La)$ un foncteur défini sur la catégorie des $\breve{\O}$-schémas formels $p$-adiques quasi-séparés sans $p$-torsion et à valeurs dans $\mathbb{D}^+ (\La)$. Supposons que :
\begin{itemize}
\item En restriction à la catégorie des $\breve{\O}$-schémas formels de type fini ce foncteur coïncide avec le foncteur $\X\longmapsto R\GG_c (\X^{an}\hat{\otimes} \Cp,\La)$ où $\X^{an}$ désigne la fibre générique au sens des espaces de Berkovich et la cohomologie est la cohomologie étale de torsion ou bien $\ell$-adique.
\item Si $\X=\underset{i\in \N}{\limp} \X_i$ où $\X_i\ldrt \X_{i+1}$ est fini alors 
$$
\underset{i\in \N}{\limi} R\GG_c (\X_i,\La) \iso R\GG_c (\X,\La)
$$
\end{itemize}
Alors 
$$
R\GG_c (\X_\infty,\La)\iso \underset{K}{\limi} R\GG_c (\M_K\hat{\otimes} \Cp, \La)
$$
Dans \cite{iso5} nous construisons un tel foncteur (et ses versions équivariantes). La condition de compatibilité aux limites projectives sera une conséquence du théorème d'approximation d'Elkik. Cela nous permet de comparer la cohomologie des tours de Lubin-Tate et de Drinfeld.
\end{rema}

\section{Décomposition cellulaire écrasée en niveau fini}

Dans la section \ref{decellrig} on a défini pour un niveau $K$ une décomposition cellulaire d'un 
ouvert admissible de $\M_K$ d'image quasi-compacte dans $\M_K/\pi^\Z$. On explique ici que quitte à modifier les cellules on a une décomposition cellulaire de tout $\M_K$ paramétrée par un quotient de l'immeuble et des cellules modifiées.

\begin{defi}[Cellules écrasées]
Soit $\s$ un simplexe orienté de $\mathcal{I}$ comme précédemment et $K\subset \GL_n (F)$ un sous-groupe compact ouvert. Soit $K'\lhd K$ compact ouvert tel que si $a=[\La,M]$ est un sommet de $\s$ alors $K'\subset \text{Id} + \pi \End (\La)$. Posons 
$$
\mathcal{D}_{\s,K} = \mathcal{D}_{\s, K'} / \text{Stab}_K (\s)
$$
où $\text{Stab}_K (\s)$ désigne l'intersection des stabilisateurs des sommets de $\s$.
\end{defi}

Ce quotient est bien défini car $\text{Stab}_K (\s)/K'$ est un groupe fini.
On a alors une décomposition cellulaire :
$$
\xymatrix{
\coprod_{\{a\drt a' \}/K} \mathcal{D}_{a\drt a',K}  \ar@<1ex>[r]\ar@<-1ex>[r]  & \coprod_{a\in \mathcal{I}/K} \mathcal{D}_{a,K} \ar[r] & \M_K
}
$$

\begin{rema}
Pour tout $K$ il existe $n-1$ sous-ensembles $K$-stables $(U_{i,K})_{1\leq i \leq n-1}$ dans $\mathcal{I}/\pi^\Z$ tels que 
\begin{itemize}
\item $\dpt{\mathcal{I}/\pi^\Z \setminus \bigcup_{1\leq i\leq n-1} U_{i,K}}$ est relativement compact dans $\mathcal{I}/\pi^\Z$
\item Si pour $1\leq i \leq n-1$ $\;P_i$ désigne le sous-groupe parabolique maximal dans $\GL_n$ stabilisateur de $F^i \oplus (0)^{n-i}$ alors
$$
\pi_0 (U_{i,K}/K) \simeq P_i (F) \bc \GL_n (F) /K \times \Z/n\Z
$$ 
\end{itemize}
On en déduit qu'il existe  des ouverts admissible $(V_{i,K})_{1\leq i
  \leq n-1}$  dans $\M/\pi^\Z$ tels que $\left ( \bigcup_i V_{i,K}
\right )^c$ est relativement compact (i.e. contenu dans une boule de
rayon $<1$) et dont l'image réciproque en niveau $K$ se scinde en une
union disjointe indexée par $P_i (F)\bc \GL_n (F) /K$ i.e. est
``induite parabolique''. On peut donner une interprétation de ce
phénomène en termes de sous-groupes canoniques.
\end{rema}

\appendix 
\section{Normalisé d'un schéma formel dans une extension de sa fibre générique}
\label{Normschf}

Dans cet appendice $K$ désigne un corps valué complet non-archimédien. 
Tous les espaces rigides considérés sont quasi-séparés et définis sur $K$. 
Par point de $X$ on entend ici les points classiques c'est à dire ceux associés aux
spectres maximaux des algèbres affinoïdes (mais tous les résultats énoncés restent valable en considérant plus
généralement les points au sens de Berkovich).

\subsection{Généralités sur les espaces rigides}
 
Nous commençons par collecter quelques résultats difficilement disponibles dans la littérature
sur les espaces rigides en donnant quelques indications sur les démonstrations dans certains cas. Nous n'auront pas besoin de tout ces énoncés dans la suite, néanmoins ils répondent à des questions qui sont venues naturellement à l'esprit de l'auteur lors de la rédaction des énoncés concernant la normalisation des schémas formels.

\begin{Fait}
Soit $X$ un espace rigide. Sont équivalents 
\begin{itemize}
\item $\forall x\in X$ l'anneau $\O_{X,x}$ est réduit
\item Il existe un recouvrement admissible affinoïde $(U_i)_i$ de $X$ tel que $\forall i$ l'anneau $\O_X (U_i)$ soit réduit
\item Pour tout ouvert admissible $U$ de $X$ l'anneau $O_X (U)$ est réduit
\end{itemize}
Si l'une des trois conditions précédentes est vérifiée $X$ est dit réduit. 
En général le faisceau qui à $U$ ouvert admissible quasicompact associe l'idéal des éléments nilpotents dans $\O_X (U)$ est cohérent et définit un sous-espace rigide Zariski fermé  réduit $X_{red}$ dans $X$. 
\end{Fait}

\begin{Fait}
Soit $X$ un espace rigide réduit. Sont équivalents
\begin{itemize}
\item $\forall x\in X$ l'anneau $\O_{X,x}$ est intègre
\item Il existe un recouvrement admissible affinoïde 
$(U_i)_i$ de $X$ tel que $\forall i\; \forall x\in U_i\;$ l'anneau
$\O_X (U_i)_{\mathfrak{m}_x}$ est intègre unibranche. 
\item Pour tout ouvert admissible connexe $U$ de $X$ l'anneau $\O_X (U)$ 
est intègre et si de plus $U$ est affinoïde
$\forall \mathfrak{P} \in \spec (\O_X (U))$ l'anneau $\O_X (U)_{\mathfrak{P}}$ est unibranche. 
\end{itemize}
Si ces conditions sont vérifiées
$X$ sera dit localement intègre. 
\end{Fait}

\dem
On utilise librement les propriétés de base concernant les anneaux locaux des espaces rigides telles que dans le chapitre 2.1 de \cite{Berk1}.
Soient $x\in X$, $U$ un ouvert admissible affinoïde contenant $x$, $\mathcal{A}=\O_X (U)$ et $\mathfrak{m}$ l'idéal de $\mathcal{A}$ associé à $x$.
Il y a un isomorphisme entre anneaux locaux noethériens 
$$
\widehat{\mathcal{A}_\mathfrak{m}} \simeq \widehat{\O_{X,x}}
$$
L'anneau $\mathcal{A}$ est réduit excellent (\cite{Kiehl1} ou \cite{Vala2}).
 On en déduit (EGA IV 7.8.3 (vii)) que $\mathcal{A}_{\mathfrak{m}}$ est intègre unibranche ssi
$\widehat{\mathcal{A}_\mathfrak{m}}$ est intègre. 
Puisque $\O_{X,x} \subset \widehat{\O_{X,x}}$ 
la seconde assertion de l'énoncé 
 entraîne la première. La troisième entraînant clairement la seconde il suffit de voir que la première implique la troisième. 
 Il suffit alors de montrer que $\O_{X,x}$ intègre implique $\mathcal{A}_{\mathfrak{m}}$ intègre unibranche
 (dans la troisième assertion $\mathfrak{m}$ n'est pas supposé maximal, mais on peut s'y ramener puisque $\mathcal{A}$ étant excellent 
 l'ensemble des points de $\spec (\mathcal{A})$ où il est unibranche est constructibe
  (EGA IV 9.7.10) et $\mathcal{A}$ est un anneau de Jacobson). Si l'on savait que
 $\O_{X,x}$ est excellent cela serait facile car $\O_{X,x}$ étant hensélien il est unibranche et donc son complété serait intègre.
 Malheureusement cela n'est pas connu. Le lemme qui suit appliqué à $\mathcal{A}_\mathfrak{m} \subset \O_{X,x}$ 
 permet néanmoins de conclure. 
 \qed

\begin{lemm}
Soit $(A,\mathfrak{m})\drt (B,\mathfrak{n})$ un morphisme local injectif d'anneaux locaux. Supposons $(B,\mathfrak{n})$ hensélien intègre. Alors $(A,\mathfrak{m})$ est unibranche.
\end{lemm}
\dem La démonstration est identique à celle de la proposition 18.6.12 de EGA IV. \qed

\begin{exem}
Si $X=\text{Spm} (\mathcal{A})$ avec $\mathcal{A}= K<x,y>/ (x^2- y^2 (1-x))$, l'algèbre $\mathcal{A}$ est intègre mais l'ouvert $|x|\leq |p|$ bien que connexe défini par l'algèbre $\mathcal{A}<\frac{x}{p}>$  est tel que 
  $\mathcal{A}<\frac{x}{p}>$ n'est pas intègre puisque si $|x|\leq |p|$ alors $\sqrt{1-x} = \sum_{k\geq 0} \left ( 1/2 \atop k \right ) (-1)^k x^k \in 
\mathcal{A}<\frac{x}{p}>$ et donc $x^2- y^2 (1-x)=( x-y \sqrt{1-x})( x+y \sqrt{1-x})$. 
Cela peut également se voir en disant que la cubique nodale n'est pas unibranche en sa singularité. 
\end{exem}

Dans l'énoncé suivant, un anneau nothérien est dit normal s'il est un produit
(nécessairement fini) 
d'anneaux intègres intégralement fermés dans leur corps des fractions, ce qui est
encore équivalent à dire que tous ses localisés en ses idéaux premiers sont
intègres intégralement clos dans leur corps des fractions.

\begin{Fait}
Soit $X$ un espace rigide réduit. Sont équivalents 
\begin{itemize}
\item Il existe un recouvrement affinoïde admissible $(U_i)_i$ de $X$ tel que $\forall i\; \O_X (U_i)$ est normal
\item Pour tout ouvert admissible connexe $U$ de $X$ $\;\; \O_X (U)$ est intègre 
intégralement clos 
\item $\forall x\in X\; \O_{X,x}$ est intègre intégralement clos dans son corps des fractions
\end{itemize} 
Si ces conditions sont vérifiées l'espace $X$ est dit normal.
\end{Fait}

\begin{Fait}
Soit $X$ un espace rigide réduit. Le normalisé $\widetilde{X}$ de $X$ est bien défini et le morphisme $\widetilde{X}\drt X$ est fini. 
\end{Fait}

Par normalisé on entend ici un objet représentant le foncteur $\Hom (-,X)$ 
restreint à la catégorie des espace rigides normaux.

Si $f: \widetilde{X} \drt X$,
le faisceaux $f_*\O_{\widetilde{X}}$ est le normalisé de $\O_X$ dans le faisceaux des
fonctions méromorphes 
$\mathcal{M}_X$. Ce faisceau des fonctions méromorphes est celui qui à $U$ ouvert affinoïde associe le corps total des fractions de $\O_X (U)$. Plus précisément,
 utilisant que $\forall U$ affinoïde $\O_X (U)$ est excellent donc Japonais et que donc $\O_{\widetilde{X}} (U)$ est fini sur $\O_X (U)$ on en déduit que $\O_{\widetilde{X}} (U)$ est une algèbre de Tate. Le normalisé de $\O_X$ dans $\mathcal{M}_X$ est un faisceau
cohérent qui définit donc un espace rigide affinoïde au dessus de $X$ 
qui est $\widetilde{X}$.

\begin{rema}
Ainsi,  $X$ réduit est localement intègre ssi
pour tout ouvert admissible connexe de $X$ son image réciproque dans $\widetilde{X}$ est connexe. 
\end{rema}

\begin{rema}
Pour un espace rigide $X$ on peut définir les composantes irréductibles de $X$
comme étant les images des composantes connexes du normalisé de $X_{red}$. Elles 
sont donc Zariski fermées puisque le morphisme de normalisation est fini.
Si pour un ouvert admissible quasicompact $U$ on note $\text{Irr} (U)$  l'ensemble de
ses composantes irréductilbles, $\text{Irr}$ est un faisceau muni d'un épimorphisme 
$\text{Irr} \twoheadrightarrow \pi_0$ associé au morphisme $\underline{\Z} \drt f_* \underline{\Z}$ où $f : \widetilde{X} \drt X$. 
\end{rema}

\subsection{Schémas formels normaux} \label{schfnor}

Désormais $K$ sera supposé de valuation discrète. Rappelons qu'un $\O_K$-schéma formel admissible
est un $\O_K$-schéma formel  de type fini sans $\pi$-torsion (où $\pi$ est une uniformisante de $K$). On les supposera toujours quasi-séparés.

\begin{Fait}\label{oktukok}
Soit $\X$ un $\O_K$-schéma formel admissible. Sont équivalents 
\begin{itemize}
\item Il existe un recouvrement ouvert affine $(\mathcal{U}_i)_i$ de $\X$ tel que $\forall i \; \O_{\X} (\mathcal{U}_i)$ est  normal
\item $\forall \mathcal{U}$ ouvert affine connexe de $X$, $\O_{\mathcal{X}} (\mathcal{U})$ est intègre intégralement clos dans son corps des fractions
\item $\forall x\in \mathcal{X}\;\; \O_{\X,x}$ est intègre intégralement clos dans son corps des fractions
\end{itemize}
\end{Fait}

\dem
On utilise le fait suivant :
 $\forall \mathcal{U}$ ouvert  affine dans $\X$ l'anneau
 $ \O_{\X} (\mathcal{U})$ est excellent (\cite{Vala1} pour le cas d'égale caractéristique et 
\cite{Vala2} pour le cas d'inégale) 
 et donc $\forall f \in  \O_{\X} (\mathcal{U})$ l'anneau 
$ \O_{\X} (\mathcal{U})[\frac{1}{f}]$ est excellent ce qui implique 
que s'il est normal alors son complété $\O_{\mathcal{U}} (D(f))$ est normal. A partir de là tout
le reste est facile.
\qed

\begin{rema}
Il n'existe pas en général de bonne notion de schéma formel réduit ou bien normal. On a en effet besoin d'utiliser des propriétés stables par complétion afin que si l'anneau $A$ possède cette
propriété, pour tout $f$ dans $A$ l'anneau $A<\frac{1}{f}>$ la possède aussi.
\end{rema}

\begin{Fait} \label{sntt}
Soit $\X$ un $\O_K$-schéma formel admissible réduit tel que
$\X^{rig}$ soit normal. Soit $\text{sp}:(\X^{rig},\O_{\X^{rig}})\drt (\X,\O_{\X})$ comme morphisme 
d'espaces annelés. Via la norme infini le faisceau $\O_{\X^{rig}}$est un faisceau en $K$-ev. normés 
: $\forall f\in \O_{\X^{rig}} (U)\; |f|_{\infty}= \sup_{x\in U} |f(x)|$. 
On note $\O_{\X^{rig}}^{\,0}$ le sous-faisceau en $\O_K$-algèbres 
des $f\in \O_{\X^{rig}}$ telles que $|f|_{\infty} \leq 1$. 
 Alors le normalisé de $\X$ est bien défini égal à 
$$
\spf (\text{sp}_*\O_{\X^{rig}}^{\,0})
$$
qui est un $\O_K$-schéma formel admissible normal fini au dessus de $\X$. 
\end{Fait}
\dem
Tout repose sur le théorème de Grauert Remmert qui assure que pour $A$ une $\O_K$-algèbre admissible,  $\mathcal{A}=A[\frac{1}{\pi}]$ l'algèbre de Tate associée, la boule unité de $\mathcal{A}$ pour la 
semi-norme infini est une algèbre admissible finie sur $A$. Nous renvoyons pour cela à la
discussion au début de la section 1 de \cite{BLR4}.
\qed

\begin{Fait}
Soit $\mathcal{U}$ un ouvert affine de $\X$ et $f_1,\dots,f_n\in \O_{X}(\mathcal{U})[\frac{1}{p}]^0$ telles que $ \O_{X}(\mathcal{U})[\frac{1}{p}]^0 = \O_{X} (\mathcal{U})[f_1,\dots,f_n]$. Soit $r\in \N$ tel que $\forall i\; \pi^r f_i \in \O_{\mathcal{U}}$. Alors, le normalisé de $\mathcal{U}$ s'identifie à l'éclatement formel admissible de l'idéal ouvert $(\pi^r,\pi^r f_1,\dots,\pi^r f_n)$.
\end{Fait}

\subsection{Normalisé dans une extension de la fibre générique}\label{norextgen}

\begin{Fait}
Soit $\X$ un $\O_K$-schéma formel admissible réduit.
 Soit $\ph: Y\ldrt \X^{rig}$ un morphisme fini d'espaces rigides tel que $Y$ soit normal. 
Alors, le normalisé $\widetilde{\X}$ de $\X$ dans $Y$ existe, est fini au dessus de $\X$ 
 et est égal à 
$$
\spf ( \text{sp}_*\ph_* \O_Y^{\, 0})
$$ Il vérifie la propriété universelle suivante : étant donné un schéma formel normal $\mathfrak{Z}$ muni d'un morphisme $\mathfrak{Z}\drt \X$ et d'un relèvement 
$$
\xymatrix{
 & Y \ar[d] \\
\mathfrak{Z}^{rig}  \ar[r]\ar[ru] & \X^{rig} 
}
$$
il existe un unique relèvement 
$$
\xymatrix{
 & \widetilde{\X} \ar[d] \\
\mathfrak{Z}  \ar[r]\ar[ru] & \X 
}
$$
tel qu'après passage à la fibre générique les deux diagrammes précédents soient compatibles.
\end{Fait}

\section{Modules de Dieudonné et cristaux des $\O$-modules $\pi$-divisibles}\label{MdcO}

Soit $F|\Qp$ une extension de degré fini. On notera $\O=\O_F$.

\subsection{Un lemme sur les $F$-cristaux $\O$-équivariants}

Soit $S$ un schéma sur lequel $p$ est localement nilpotent et $\Sigma= \spec (\Zp)$. On considère le 
gros site cristallin $\text{CRIS} (S/\Sigma)$ de \cite{BBM} ou bien le gros site cristallin nilpotent $\text{NCRIS} (S/\Sigma)$ (\cite{BBM} chapitre 1).

\begin{lemm}
Soit $\E$ un $F$-cristal (non-dégénéré) en $\O_{S/\Sigma}$-modules localement libres de rang fini sur 
$\text{CRIS} (S/\Sigma)$ ou $\text{NCRIS} (S/\Sigma)$.  
 Supposons $\E$ muni d'une action de $\O$. Alors, $\E$ est un
$\O_{S/\Sigma}\otimes_{\Zp} \O$-module 
 localement libre sur $\text{CRIS} (S/\Sigma)$.
\end{lemm}
\dem
Soit $(U\hookrightarrow T)\in \text{CRIS} (S/\Sigma)$ (resp. $\text{NCRIS} (S/\Sigma)$) 
  et $x\in U$. Les assertions suivantes sont équivalentes 
\begin{itemize}
\item[(i)] $\E_{(U\hookrightarrow T)}\otimes_{\O_T} \O_{T,x}$ est libre en tant que $\O_{T,x}\otimes_{\Zp} \O_F$-module
\item[(ii)] $\E_{(U\hookrightarrow T)}\otimes_{\O_T} k(x)$ est libre en tant que $k(x)\otimes_{\Zp} \O_F$-module
\item[(iii)] $\E_{(U\hookrightarrow T)}\otimes_{\O_T}\overline{k(x)}$ est libre en tant que $\overline{k(x)}\otimes_{\Zp} \O_F$-module
\end{itemize}
L'équivalence entre $(ii)$ et $(iii)$ résulte de ce que $\overline{k(x)}\otimes \O_F$ est une $k(x)\otimes \O_F$-algèbre fidèlement
plate. 
\\
L'implication $(i)\impl (ii)$ est claire, quant à $(ii)\impl (i)$ elle résulte du lemme de Nakayama. En effet, si 
$\E_{(U\hookrightarrow T)}\otimes_{\O_T} k(x)$ est un $k(x)\otimes \O_F$-module libre soit 
$$
e_1,\dots,e_r\in \E_{(U\hookrightarrow T)}\otimes_{\O_T} \O_{T,x}
$$
un relèvement d'une $k(x)\otimes \O_F$-base . Soit $(\epsilon_j)_j$ une base de $\O_F$ comme $\Zp$-module. Alors, d'après le lemme de Nakayama, 
$(\epsilon_i e_j)_{i,j}$ engendre $\E_{(U\hookrightarrow T)}\otimes_{\O_T} \O_{T,x}$ comme $\O_{T,x}$-module. Mais 
$$
\text{rg}_{\O_{T,x}} \left (\E_{(U\hookrightarrow T)}\otimes \O_{T,x} \right ) =\text{rg}_{k(x)} \left ( \E_{(U\hookrightarrow T)}\otimes_{\O_T} k(x)
\right ) = [F:\Qp] \;\text{rg}_{k(x)\otimes \O_F} \left ( \E_{(U\hookrightarrow T)}\otimes_{\O_T} k(x) \right )
$$
La famille $(\epsilon_i e_j)_{i,j}$ est donc une $\O_{T,x}$-base.
\\

Montrons donc l'assertion $(iii)$. Il y a des morphismes compatibles dans $\text{CRIS} (S/\Sigma)$ (resp. $\text{NCRIS} (S/\Sigma)$) 
pour $n\in \N$
$$
\xymatrix{
\spec (\overline{k(x)}) \ar@{^(->}[r]\ar[d] & \spec (W_n (\overline{k(x)})) \ar[d]  \\
U \ar@{^(->}[r] & \spec ( W_n (\O_T)) 
}
$$
La propriété de cristal implique donc que $\E_{(U\hookrightarrow T)}\otimes_{\O_T}\overline{k(x)}$ est la réduction modulo $p$ de 
$$
E= \underset{n}{\limp} \E_{(\spec ( \overline{k(x)} ) \hookrightarrow \spec ( W_n ( \overline{ k(x)}) ) }
$$
C'est un $W(\overline{k(x)})$-module libre qui est en fait un $W(\overline{k(x)})\otimes_{\Zp} \O_F$-module muni d'une isogénie $F$-linéaire (ici $F=$ Frobenius)
$$
\ph : E_\Q \iso E_\Q
$$
qui commute à l'action de $\O_F$. Or
$$
E=\bigoplus_{\tau : F^0 \hookrightarrow W ( \overline{k(x)})_\Q} E_\tau
$$
où $E_\tau$ est  un $\O_F\otimes_{\O_{F^0},\tau} W ( \overline{ k(x)})$-module libre. Le module $E$ est donc libre sur 
$W(\overline{k(x)})\otimes_{\Zp} \O_F$ ssi 
$$
\forall \tau,\tau'\;\;\; \text{rg}_{\O_F\otimes_{\O_{F^0},\tau} W ( \overline{ k(x)})} \; E_\tau = \text{rg}_{\O_F\otimes_{\O_{F^0},\tau'} W ( \overline{ k(x)})} \; E_{\tau'} 
$$
Mais 
$$
\forall \tau \;\;\; \ph : E_{\tau} \otimes \Q \iso E_{F \tau} \otimes \Q 
$$
\qed

\begin{coro}\label{liextul}
Soit $H$ un groupe $p$-divisible sur un schéma $S$ sur lequel $p$ est localement nilpotent. Supposons le muni d'une action de $\O$.  Alors l'algèbre de Lie de l'extension vectorielle universelle de $H$ est un $\O_S\otimes_{\Zp} \O$-module localement libre.
\end{coro}
\dem
C'est une conséquence de ce que cette algèbre de Lie est l'évaluation d'un des cristaux définis dans \cite{Messing1}.
\qed

\subsection{Structure du cristal de Messing d'un $\O$-module $\pi$-divisible}

\begin{defi}
Un $\O$-module $\pi$-divisible sur un $\O$-schéma est un groupe $p$-divisible muni d'une action de $\O$ qui induit l'action naturelle sur son algèbre de Lie. 
\end{defi} 

\begin{rema}
Soit $H$ un groupe $p$-divisible muni d'une action de $\O$. Alors,
$\forall n\;\; H[\pi^n]$ est un groupe plat fini. En effet, $\pi :
H\ldrt H$ est une isogénie puisque $\pi^{e-1} .\pi =\text{unité}\times
\pi$. 
On peut ainsi donner une définition analogue à celle des groupes p-divisibles pour les groupes p-divisibles munis d'une action de $\O$ en remplaçant $p$ par $\pi$ dans les définitions.  
\end{rema}

Soit $S$ un $\O$-schéma sur lequel $p$ est localement nilpotent. Soit $H$ un
$\O$-module $\pi$-divisible sur $S$.
Considérons la suite exacte 
$$  
0 \ldrt \omega_{H^D} \ldrt \text{Lie}\; E (H) \ldrt \omega_H^* \ldrt 0
$$  
où $E(H)$ désigne l'extension vectorielle universelle de $H$ de \cite{Messing1}. 
Soit $I$ l'idéal d'augmentation de $\O_S\otimes_{\Zp} \O_F \twoheadrightarrow \O_S$.

\begin{prop} \label{aolgh}
Localement sur $S$ il existe une $\O_S\otimes \O_F$-base de $\text{Lie} E(H)$
$\;\; (e_1,\dots,e_{r_1+ r_2})$ telle que 
$$
\omega_{H^D} = \O_S\otimes O_F.e_1\oplus \dots \oplus  \O_S\otimes O_F.e_{r_1} \oplus
I.e_{r_1+1}\oplus \dots \oplus I.e_{r_1+ r_2}
$$
\end{prop}
\dem
Il y a une suite exacte 
$$
0\ldrt \omega_{H^D} /I \text{Lie} E(H) \ldrt \text{Lie} E(H) /I \text{Lie} E(H) \ldrt \omega_H^* \ldrt 0 
$$
où d'après le corollaire \ref{liextul} $\text{Lie}  E(H)/I\text{Lie} E(H)$ est un $\O_S$-module localement libre.
Localement sur $S$ soit $(e_1,\dots,e_{r_1+r_2})$ un relèvement d'une base de 
$\text{Lie} E(H)/I \text{Lie} E(H)$ dans $\text{Lie} E(H)$ tel que 
$(e_{r_1+1},\dots, e_{r_1+r_2})$ s'envoie sur une base de $\omega_{H}^*$ via 
$\text{Lie} E(H) \twoheadrightarrow \omega_H^*$ (une telle base existe puisque
$\text{Lie} E(H)/I\text{Lie} E(H) \twoheadrightarrow \omega_H^*$ est localement scindé). 
L'idéal $I$ étant contenu dans le radical de Jacobson de $\O_S\otimes \O_F$ et 
$\text{Lie} E(H)$ étant libre sur $\O_S\otimes \O_F$, $(e_1,\dots,e_{r_1+r_2})$ est une 
$\O_S\otimes \O_F$-base de $\text{Lie} E(H)$ dont on vérifie aussitôt qu'elle convient.
\qed

\subsection{$\O$-extension vectorielle universelle d'un $\O$-module $\pi$-divisible}

\begin{defi}
Soit $S$ un schéma et
 $\F$ un faisceau en groupe abéliens sur $S_{fppf}$. Soit $f: \spec (\O_S [\epsilon]) \drt S$. On appelle algèbre de Lie
 de $\F$ le faisceau en $\O_S$-modules sur $S_{fppf}$
 $$
 \Lie \,\F = \ker f_* f^* \F \xrig{\;\epsilon =0 \;} \F 
 $$
\end{defi}

\begin{defi}
Soit $H$ un $\O$-module $\pi$-divisible sur un schéma $S$.
On appelle $\O$-extension vectorielle de $H$ une extension
$$
0 \ldrt \underline{V} \ldrt  E \ldrt H \ldrt 0
$$
de faisceaux en $\O_F$-modules sur $S_{fppf}$ telle que $V$ soit un $\O_S$-module cohérent,
$\underline{V}$ le faisceau fppf associé et l'action induite de $\O$ sur $\Lie \, E$ soit l'action naturelle.
\end{defi}

On remarquera que de telles extensions sont ``rigides'' puisque $\forall W$ cohérent 
$\Hom (H,\underline{W})=0$ ($p$ est localement nilpotent sur $S$)(par fainéantise on oubliera désormais de souligner les faisceaux cohérents et on notera $V$ pour $\underline{V}$).  
La proposition suivante a donc un sens.
  
\begin{prop}\label{exOver}
Tout $\O$-module $\pi$-divisible possède une $\O$-extension vectorielle universelle 
$$
0\ldrt V_\O (H) \ldrt E_\O (H) \ldrt H \ldrt 0
$$
De plus $V_\O (H)$ et $\Lie\, E_\O (H)$ sont des $\O_S$-modules 
localement libres et il y a une suite exacte
$$
0 \ldrt V_\O (H) \ldrt \Lie\; E_\O (H) \ldrt \Lie\; H \ldrt 0
$$
\end{prop}
\dem
Soit 
$$
0 \ldrt \omega_{H^D} \ldrt E(H) \ldrt H \ldrt 0
$$
L'extension vectorielle universelle de $H$ (correspondant au cas $\O_F=\Zp$). Soit 
$I=\ker \left ( \O_S\otimes \O_F\twoheadrightarrow \O_S \right )$.
Considérons le poussé en avant de l'extension précédente par le morphisme
$\omega_{H^D} \twoheadrightarrow \omega_{H^D}/I\Lie E(H)$
$$
\xymatrix{
0 \ar[r] & \omega_{H^D} \ar[r]\ar@{->>}[d] & E(H) \ar[r]\ar[d] & H \ar[r] \ar@{=}[d] & 0 \\
0\ar[r] &  \omega_{H^D}/I \Lie E(H) \ar[r] & \widetilde{E} \ar[r] & H \ar[r] & 0 
}
$$  
où
$$
 \widetilde{E} =   \omega_{H^D}/I \Lie E(H) \coprod_{\omega_{H^D}} E(H) 
=  E(H)/I \Lie E(H)
$$
Alors cette extension est une $\O$-extension vectorielle universelle de $H$. Avant de le montrer commençons par montrer quelques propriétés de cette extension. \\
D'après la proposition \ref{aolgh} 
 $ \omega_{H^D}/I \Lie E(H)$ est localement libre. De plus, si $\widetilde{E}_n$ désigne l'image réciproque de $H[\pi^n]$ dans $\widetilde{E}$, $\widetilde{E} =\underset{n}{\limi} \widetilde{E}_n$ et il y a des suites exactes
$$
0 \ldrt  \omega_{H^D}/I \Lie E(H) \ldrt  \widetilde{E}_n \ldrt H[\pi^n] \ldrt 0
$$
Donc, $ \widetilde{E}_n$ est un torseur ffpf sous l'espace affine $ \omega_{H^D}/I \Lie E(H)$ au dessus de $H[\pi^n]$ et est donc représentable par un $H[\pi^n]$-schéma lisse. Il y a donc des 
suites exactes pour tout $n$ 
$$
0\ldrt  \omega_{H^D}/I \Lie E(H) \ldrt \Lie \widetilde{E}_n \ldrt \Lie H[\pi^n] \ldrt 0
$$
et étant donné que $\Lie \widetilde{E} =\underset{n}{\limi} \widetilde{E}_n$, que pour $n>>0\;\;
\Lie H[\pi^n]=\Lie H$, pour $n>>0\;\; \Lie \widetilde{E}_n= \Lie \widetilde{E}$ et il y a une suite
exacte
$$
0 \ldrt  \omega_{H^D}/I \Lie E(H) \ldrt \Lie \widetilde{E} \ldrt \Lie H \ldrt 0
$$
et même un diagramme commutatif
$$
\xymatrix{
0 \ar[r] & \omega_{H^D} \ar[r]\ar@{->>}[d] & \Lie E(H) \ar[r]\ar[d] & \Lie H \ar[r] \ar@{=}[d] & 0 \\
0\ar[r] &  \omega_{H^D}/I \Lie E(H) \ar[r] & \Lie \widetilde{E} \ar[r] & \Lie H \ar[r] & 0 
}
$$
duquel on déduit que 
$$
\Lie \widetilde{E}= \Lie E(H)/I.\Lie E(H)
$$
Montrons maintenant l'universalité de l'extension. Soit 
$$
0\ldrt \underline{V} \ldrt E \ldrt H \ldrt 0
$$
une $\O$-extension vectorielle. Oubliant l'action de $\O$, celle-ci est induite par un unique morphisme $\O_S$-linéaire $\omega_{H^D}\ldrt V$
$$
\xymatrix{
0 \ar[r] & \omega_{H^D} \ar[r]\ar@{->>}[d] & E(H) \ar[r]\ar[d] & H \ar[r] \ar@{=}[d] & 0 \\
0 \ar[r] & \underline{V} \ar[r] & E \ar[r] & H \ar[r] & 0
}
$$
l'unicité d'un tel morphisme (qui résulte de $\Hom (H, \underline{V})=0$) implique que 
tous les morphismes dans le diagramme précédent sont $\O$-linéaires. Il y a alors un diagramme
de $\O_S\otimes \O$-modules sur $S_{fppf}$

$$
\xymatrix{
0 \ar[r] & \omega_{H^D} \ar[r]\ar@{->>}[d] & \Lie E(H) \ar[r]\ar[d] & \Lie H \ar[r] \ar@{=}[d] & 0 \\
0 \ar[r] & \underline{V} \ar[r] & \Lie E \ar[r] & \Lie H \ar[r] & 0
}
$$
duquel on déduit que $\omega_{H^D} \ldrt \underline{V}$ se factorise par 
$\omega_{H^D} \twoheadrightarrow \omega_{H^D}/I.\Lie E(H)$. 
\qed   

\begin{rema}
Il résulte de la démonstration précédente que si l'on pose $\text{ht}_\O H = \text{ht}\; H/ [F:\Qp]$ alors
$$
\text{rg}_{\O_S} V_\O (H) = \text{ht}_\O\, H - \dim H
$$
\end{rema}

\begin{rema}
On vérifie comme dans le lemme 1.18 de \cite{Messing1} que $E_\O (H)$ est formellement lisse.
\end{rema}

\begin{rema}
Comme dans la proposition 1.19 du chapitre IV de \cite{Messing1} on vérifie que le complété formel 
(au sens de la section II.1.0 de \cite{Messing1}) $\widehat{E}_\O (H)$ de $E_\O (H)$ est un $\O$-module formel (i.e. un groupe de Lie formel muni d'une action de $\O$ telle que l'action
induite sur l'algèbre de Lie soir l'action naturelle) extension du $\O$-module formel
$\widehat{H}$  par le groupe formel vectoriel 
$\widehat{V}_\O (H)$ (localement isomorphe à une somme finie de $\widehat{\mathbb{G}}_a$).
\end{rema}

\subsection{Cristal de Messing généralisé et théorie de la déformation}

\begin{theo}
Soit $\Sigma= \spec (\O_F)$ muni de l'idéal à puissances divisées $(p)$ et $\E$ le cristal
de Messing de $H$ sur $\text{NCRIS} (S/\Sigma)$ en tant que $\O_{S/\Sigma}\otimes \O_F$-module.
Soit $\mathcal{I}=\ker \left ( \O_{S/\Sigma} \otimes \O_F \twoheadrightarrow \O_{S/\Sigma}\right )$. 

Alors, 
$\E/\mathcal{I} \E$ est un cristal en $\O_{S/\Sigma}$-modules localement libres de rang fini tel que
$\forall (U\hookrightarrow T)\in \text{NCRIS} (S/\Sigma)\;\; \forall \widetilde{H} $ un relèvement de $H\times_S U$ sur $T$ comme $\O$-module $\pi$-divisible 
$$
\left ( \E/\mathcal{I}\E\right )_{(U\hookrightarrow T)}
$$
s'identifie à l'algèbre de lie de $E_\O (\widetilde{H})$, la $\O$-extension vectorielle universelle de $\widetilde{H}$.

De plus, pour un tel relèvement $\widetilde{H}$ la partie vectorielle de $E_\O (\widetilde{H})$,
$V_\O (\widetilde{H}) = \omega_{\widetilde{H}^{D}}/ I. \Lie E_\O (\widetilde{H})$ définit une
filtration localement facteur direct dans $\left ( \E/\mathcal{I} \E \right )_{(U\hookrightarrow T)}$
se réduisant modulo l'idéal de $U$ dans $T$ sur la partie vectorielle de $E_\O (H)$.
\\
Cette correspondance définit une équivalence de catégories entre la catégorie des $\O$-modules
 $\pi$-divisibles sur $T$ et celle des couples $(H_0,\text{Fil})$ où $H_0$ est un $\O$-module $\pi$-divisible sur $U$ et $\text{Fil}$ une filtration localement facteur direct dans l'évaluation du cristal précédent sur $U\hookrightarrow T$ se réduisant sur la partie vectorielle de 
$\Lie\; E_\O (H_0)$.
\end{theo}
\dem
Seule la vérification du fait que la correspondance évoquée est une équivalence de catégories
reste à vérifier.

Elle se déduit du théorème
de Messing (th. 1.6 chap. V de \cite{Messing1}) 
 en remarquant que les filtrations localement facteur direct
 $\O$-stables
$\text{Fil}$ dans $\E_{(U\hookrightarrow T)}$   et telles que l'action de $\O$ sur $\E_{(U\hookrightarrow T)}/\text{Fil}$ soit l'action naturelle sont en bijection via 
l'application
$$
\text{Fil} \longmapsto \text{Fil}/ \mathcal{I}_{(U\hookrightarrow T)} \E_{(U\hookrightarrow T)}
$$
avec les filtrations localement facteur directe dans $(\E/\mathcal{I} \E)_{(U\hookrightarrow T)}$. 
\qed

\subsection{Exponentielle $\pi$-adique}

\subsubsection{$\O$-puissances divisées (\cite{HopkinsGross} section 10, \cite{Faltings7} section 7) }
 
Soit $T$ un $\O$-schéma sur lequel $\pi$ est localement nilpotent et $J\subset \O_T$ un idéal cohérent. 

\begin{defi}
Une structure de $\O$-puissance divisées sur $J$ consiste en la donnée d'une 
application $\gamma: J \ldrt J$ vérifiant
\begin{itemize}
\item[(i)] $\forall a\in \O_S\; \forall x \in J \;\; \gamma (ax) = a^q \gamma (x)$ \item[(ii)] $\forall x\in J \; \pi \gamma (x) = x^q $
\item[(iii)] $ \forall x,y\in J\;
\gamma (x+y) = \gamma (x) + \gamma (y) + \sum_{0< i < q } \a_i \g_i (x) \g_{q_i} (y)$
où  $\a_i = \left ( \,^q_i \right )/\pi \in \O$
\end{itemize}
\end{defi}

\begin{rema}
Lorsque $\O=\Zp$, d'après la proposition A.1 de l'appendice de \cite{BerthelotOgus}, 
on retrouve la définition des puissances divisées usuelles avec $\gamma = \gamma_p$. 
 Par exemple, lorsque la base ne possède pas de $p$-torsion, $x\mapsto \frac{x^p}{p}$ permet de retrouver les $\frac{x^n}{n !}$ par composition et multiplication par des éléments de $\Zp$.
Plus précisément, $\gamma_{p^n} = \a_n \gamma^{\circ n}$ où $\a_n \in \Zp^\times$ 
 et si $n=\sum_i a_i p^i$ où $0\leq a_i \leq p-1$, $\gamma_n =\beta_n \prod_i \gamma_{p^i}\circ \gamma_{a_i}$ où $\beta_n\in \Zp^\times$.
\end{rema}

\begin{rema}
L'élément $\pi$ étant localement nilpotent l'existence d'une 
$\O$-P.D. structure sur un idéal implique que celui-ci est un nil-idéal.
\end{rema}

\begin{exem}
Si $J^2=\pi J=0$  une structure de $\O$-P.D. sur $J$  est  la même chose qu'un  morphisme additif $Frob_q$-linéaire $\gamma :J\drt J$. 
\end{exem}

\begin{defi}
Si $\g$ est une  $\O$-P.D. structure sur $J$ on pose comme dans \cite{Faltings7}
$$
\forall n\geq 1\;\; \delta_n (x) = ( \underbrace{\g\circ \dots \circ \g }_{n-\text{fois}}
) (x) . \pi^{1+q+\dots + q^{n-1}-n}
$$
et $\delta_0 (x) =x$.
\end{defi}

Ainsi, si la base est sans $\pi$-torsion $\delta_n (x) =\frac{x^{q^n}}{\pi^n}$.  On note $J^{[n]}$ l'idéal engendré par les $\gamma^{\circ a_1} (x_1)\dots \gamma^{\circ a_k} (x_k)$ où $\sum_i q^a_i \geq n$ et on dit que $\g$ est nilpotente si $J^{[n]}=0$ pour $n>>0$. Lorsque $\O=\Zp$ on vérifie facilement que l'on retrouve la définition usuelle pour les idéaux $J^{[n]}$. 

\begin{lemm}\label{extPD}
Soit $T'\drt T$ un morphisme plat et $J\subset \O_T$ muni de $\O$-P.D.. Alors celles-ci s'étendent de façon compatible à $J\O_{T'}$. 
\end{lemm}
\dem
Procéder par exemple comme dans le lemme 1.8 p.80 de \cite{Messing1}.
\qed

\begin{defi}
Pour $\Sigma= \spec (\O)$ et $S$ un $\O$-schéma sur lequel $\pi$ est localement nilpotent 
la définition des $\O$-P.D. et le lemme précédent permet de définir des sites cristallins $\text{CRIS}_\O (S/\Sigma)$ et
$\text{NCRIS}_\O (S/\Sigma)$.
\end{defi}

\subsubsection{Logarithme}

\begin{lemm}\label{wittdiv}
Supposons $J$ muni de $\O$-P.D.. Il y a alors un isomorphisme de $\O$-modules 
\begin{eqnarray*}
W_\O (J) &\iso & J^\N \\
 \text{[} x_i \text{]}_i &\longmapsto & (\underbrace{\g_i (x_0) + \g_{i-1} (x_1) + \dots + \g_0 (x_i)}_{ {''} \frac{\mathcal{W}_{\O,i} (\underline{x})}{\pi^i}''})_i
\end{eqnarray*}
Si de plus $\delta_n (J)=0$ pour $n>>0$ il y a une application 
$$
\widehat{W}_\O (J) \xrig{\; \text{log} \;} J^{(\N)} = \Lie \;\widehat{W}_\O \otimes J
$$
Si de plus les $\O$-P.D. sont nilpotentes alors l'application précédente est un isomorphisme. 
\end{lemm}
\dem
Pour la première assertion il suffit de vérifier que l'enveloppe à $\O$-P.D. de $\O[X_i]_{i\geq 0}$ par rapport à $(X_i)_{i\geq 0}$ est sans $\pi$-torsion et on se ramène alors à montrer que $\widehat{W}_\O (J)\drt J^\N$ est un morphisme lorsque $\pi$ est inversible, ce qui est clair. Le reste est facile.
\qed

\begin{rema}\label{calFV}
Du coté de $J^\N$ les opérateurs $F,V$ et $[a]$ sont donnés par 
\begin{eqnarray*}
[a]. (y_i)_i &=& (a y_i)_i \\
\,^V (y_i)_i &= &(0,y_0,y_1,\dots, y_i,\dots ) \\
\,^F ( y_i)_i &=& (\pi y_1, \pi y_2,\dots, \pi y_i ,\dots) 
\end{eqnarray*}
\end{rema}

\begin{prop}\label{plog}
Soit $G$ un $\O$-module formel sur un $\O$-schéma et $J$ un idéal muni de $\O$-puissances divisées $\gamma$ telles que localement $\delta_n (J)=0$ pour $n>>0$ (par exemple, $\pi$ est localement nilpotent). 
 Il existe alors une unique
application logarithme fonctorielle en $(G,J,\gamma)$
$$
\log_G : G(J) \ldrt  \Lie\; G\otimes J
$$
compatible avec celle définie pour $\widehat{W}_\O$. Si les $\O$-P.D. sont nilpotentes cette application est un isomorphisme.
\end{prop}
\dem 
Soit 
$$
M= \Hom (\widehat{W}_\O,G)
$$
le module de Cartier de $G$, un $\mathbb{E}_\O$-module. Alors, 
$$
G(J) \simeq \widehat{W}_\O (J) \otimes_{\mathbb{E}_\O} M 
$$
où $\mathbb{E}_\O$ agit sur $\widehat{W}_\O$ de telle manière que 
$$
x\otimes Fm = \,^V x \otimes m, \;\;\; x\otimes Vm= \,^F x \otimes m,\;\;\; x\otimes [a] m = [a] x \otimes m
$$
Posons 
$$
log_G : \widehat{W}_\O (J)\otimes_{\mathbb{E}_\O} M \ldrt J^{(\N)} \otimes_{\mathbb{E}_\O} M
$$
où $ J^{(\N)}$ est un $\mathbb{E}_\O$-module via les formules de la remarque \ref{calFV}. Alors, 
$$
J^{(\N)} = J\oplus I_{J} 
$$
où $$J = J \oplus 0 \oplus \dots \oplus 0 \oplus \dots $$
$$
I_J= 0\oplus J^{[\N)} = \,^V J^{(\N)}
$$
De plus remarquons que 
$$
\boxed{\,^F J=0 }
$$
On en déduit aussitôt que
$$
J^{(\N)} \otimes_{\mathbb{E}_\O} M = J \otimes_R M/VM
$$
\qed

\begin{rema}
L'annulation $\,^F J=0$ utilisée dans la démonstration précédente 
est à la base des liens reliant puissances divisées et théorie
de Cartier, cf. par exemple le lemme 38 de \cite{Zink2}.
\end{rema}

\begin{rema}
La notation $\log_G$ n'est pas correcte, on devrait plutôt noter $\log_{G,\gamma}$. 
\end{rema}

\begin{coro}
Soit $G$ un $\O$-module formel sur $\Sigma$ un $\O$-schéma et $S=\Sigma \text{ mod }\pi$. 
 Soit $\mathcal{J}_{S/\Sigma}$ l'idéal à puissances divisées sur $\text{NCRIS}_\O (S/\Sigma)$. Il y a alors un isomorphisme
$$ 
G(\mathcal{J}_{S/\Sigma}) \iso \Lie\, G\otimes_{\O_\Sigma} \mathcal{J}_{S/\Sigma}
$$ 
\end{coro}

\begin{prop}\label{logcar}
Soit $(G,J,\gamma)$ comme précédemment. Supposons que $J^2=\pi J= 0$. Soit 
$\a: \Lie \, G/\pi Lie\, G \ldrt \Lie \,G/\pi Lie \, G$ l'application $\text{Frob}_q$-linéaire qui à une dérivation invariante $d$ associe $d^q$. L'application $\gamma : J\drt J$ est  $\text{Frob}_q$-linéaire et il y a un diagramme commutatif naturel en $(G,J,\gamma)$ 
$$
\xymatrix{
 & \Lie \, G \otimes J \ar@{=}[r] &\Lie\, G/ \pi \Lie\, G \otimes J  \ar[dd]^{Id+ \a\otimes \gamma}_\simeq \\
G(J) \ar[ru]_{\sim}^{can} \ar[rd]_{\sim}^{log_G}  \\
 & \Lie \, G \otimes J \ar@{=}[r] & \Lie G\, / \pi \Lie\, G \otimes J
}
$$
où $can$ est l'isomorphisme définissant $\Lie \, G$.
\end{prop}
\dem Si $M$ est le module de Cartier de $G$ et $\pi=0$ sur la base,  l'application $\a$ est 
$F: M/VM\drt M/VM$ qui est bien définie puisque $FV=VF$. Pour vérifier l'assertion il 
suffit de la faire pour $\widehat{W}_\O$. En effet, après choix d'une $V$-base du module de 
Cartier de $G$ il existe une surjection $\oplus \widehat{W}_\O (J) \twoheadrightarrow G(J)$ et par naturalité du log en le groupe formel (éventuellement de dimension infinie) on conclu.

Dans le cas de $\widehat{W}_\O$, 
\begin{eqnarray*}
can : \widehat{W}_\O (J) & \iso & \Lie \widehat{W}_\O \otimes J = J^{(\N)} \\
 \sum_{i\geq 0} V^i [x_i] & \longmapsto & (x_i)_{i\geq 0}
\end{eqnarray*}
et $\a \left (  (x_i)_{i\geq 0} \right )= (0,x_1,x_2,\dots)$. Étant donné que $\pi J=0$, $\forall i\geq 2 \;\delta_i (J)=0$ donc
$$
log \left (  \sum_{i\geq 0} V^i [x_i]\right ) = ( x_0, x_1+\gamma (x_0),\dots, x_i+ \gamma (x_{i-1}),\dots )
$$
d'où le résultat.
\qed

\begin{rema}
Lorsque $J^2=0$
l'isomorphisme $can$ correspond au cas du logarithme lorsque l'on muni l'idéal des $\O$-P.D. triviales $\gamma =0$. 
\end{rema}

\subsubsection{Exponentielle}

Nous n'utiliserons pas la proposition qui suit plus tard car nous aurons besoin d'exponentier des morphismes pas seulement sur la partie formelle de nos groupes mais sur tout le groupe.
 Néanmoins cette proposition est complémentaire de celle que nous utiliserons (noter par exemple que les $\O$-P.D. ne sont pas nilpotentes
dans la proposition qui suit ce qui est compensé par le fait que nos groupes sont formels).

\begin{prop}
Soient $G$ et $H$ deux $\O$-modules formels sur $T$ et $J\subset \O_T$ un idéal cohérent
      muni de $\O$-P.D.  $ \gamma$. Il existe alors une application
$$
exp : \Hom (G, \widehat{J.\Lie H}) \ldrt \ker \left ( \Hom (G,H) \drt \Hom (G \text{ mod } J, H \text{ mod } J) \right )
$$
où si $\widehat{\Lie H}$ désigne le groupe de Lie formel associé à $\Lie H$ alors 
$ \widehat{J.\Lie H}$ est le sous-foncteur de $\widehat{\Lie H}$ défini par 
$(\widehat{J.\Lie H}) (Z) = \GG (Z,  \Lie H \otimes J (\O_Z)_{nilp} ) $.
 Cette application est naturelle en $(G,H,J,\gamma)$.
  
Si $J^2= (0)$ et $\gamma=0$ alors l'application précédente coïncide avec l'application identité donnée par 
$\Lie H \otimes J \simeq \ker ( H \drt H \text{ mod } J)$
\begin{eqnarray}\label{candiv}
\ker \left ( \Hom (G,H) \drt \Hom (G \text{ mod } J, H \text{ mod } J) \right ) = \Hom (G, J\Lie \, H)
\end{eqnarray}
\\
Si $J^2=\pi J=(0)$, soit $\a: \Lie\, H/\pi \, \Lie H \drt \Lie\, H/\pi \, \Lie H$ l'application qui à une dérivation invariante $d$ associe
$d^q$ et $\Pi=\a \otimes \gamma $ l'endomorphisme  $Frob_q$-linéaire de $J \Lie H$. Alors, naturellement en $(G,H,J,\gamma)$, via l'isomorphisme canonique (\ref{candiv}) 
$$
exp\, f =  \left ( \sum_{n\geq 0} (-1)^n \Pi^n \right ) \circ f
$$ 
où la somme infinie est finie lorsqu'elle est évaluée sur une algèbre nilpotente puisque $\text{Im} f \subset  \widehat{J.\Lie H}$.
\end{prop}

\begin{rema}
 Si $G = \widehat{\mathcal{M}}$ où $\mathcal{M} $ est un module localement libre de rang fini, il y a un morphisme 
$\Hom_{\O_T} (\mathcal{M}, J.\Lie H) \ldrt \Hom ( \widehat{\mathcal{M}}, \widehat{J.\Lie H})$, d'où une application exponentielle de source $\Hom_{\O_T} (\mathcal{M}, J.\Lie H)$. 
\end{rema}

\begin{rema}
La formule donnée pour $\exp \, f$ lorsque $J^2 = \pi J= 0$ coïncide bien avec celle de \cite{Messing1} ((2.6.6.3) page 142) lorsque $\O=\Zp$. En effet, si $(\gamma_n)_{n\geq 1}$ sont des puissances divisées alors $\gamma_{p^n} \equiv (-1)^n \gamma_p^n \text{ mod } p$. 
\end{rema}

\dem
Si $Q$ est un $\O_T$-module notons $C(Q)$ le $\mathbb{E}_\O$-module à gauche défini par
$$
C(Q)= Q^\N = \{\sum_{i\geq 0} \; V^i x_i\;|\; x_i \in Q \; \}
$$
sur lequel 
\begin{eqnarray*}
V \sum_{i\geq 0} V^i x_i &= &\sum_i V^{i+1} x_i \\
F \sum_{i\geq 0} V^i x_i &=& \pi \sum_{i\geq 1} V^{i-1} x_i \\
\text{[} a\text{]} \sum_{i\geq 0} V^i x_i &=& \sum_{i\geq 0} V^i \left ( a ^{q^i} x_i \right )
\end{eqnarray*}
Lorsque $Q$ est localement libre de rang fini $C(Q)$ est le module de Cartier de $\widehat{Q}$.
De l'isomorphisme $W_\O (J)\simeq J^\N$ du lemme \ref{wittdiv}
on déduit une décomposition $W_\O (J)=J\oplus \,^V W_\O (J)$ et donc $J\hookrightarrow \mathbb{E}_\O$. Soit $M$ le module de Cartier de $H$. L'application de réduction modulo $VM$ induit un isomorphisme 
$$
\psi : JM \iso J. (M/VM) = J\otimes M/VM
$$
Il suffit en effet de vérifier que $JM\cap VM= (0)$, mais si $(e_\a)_\a$ est une V-base de $M$ et 
$(\l_\a)_\a $ une famille d'éléments de $J$ alors $\sum_\a \l_a e\a \in VM \ssi \forall a\; \mathcal{W}_{\O,0} (\l_a)=0 \ssi \forall \a\; \l_\a=0$. 

De plus, l'inverse de cet isomorphisme induit un morphisme de $\mathbb{E}_\O$-modules 
\begin{eqnarray*}
\Xi : C ( J\otimes M/VM) &\ldrt & M \\
\sum_i V^i x_i & \longmapsto & \sum_i V^i \psi^{-1} (x_i) 
\end{eqnarray*}
En effet, cela résulte aussitôt de ce que $\,^F J= (0)$. 
L'application exponentielle s'écrit alors comme la composée
$$
\Hom (G, \widehat{J.\Lie H}) \ldrt \Hom_{\mathbb{E}_\O} \left (M_G, C ( J. M/VM)\right ) \xrig{\;\;\Xi_*\;\;} \Hom_{\mathbb{E}_\O} ( M_G,M)= \Hom (G,H)
$$
où $M_G$ est le module de Cartier de $G$ et où il faut vérifier que si 
$f:G\ldrt \widehat{J.\Lie H} \hookrightarrow \widehat{\Lie H}$ alors le morphisme 
$M_G\ldrt C (\Lie (H))$ se factorise par $C(J.\Lie H)$, mais cela résulte de la description
des modules de Cartier en termes de courbes $p$-typiques.
\\

Quant aux deux dernières assertions concernant les cas $J^2=0$ et $\gamma =0$, resp. $J^2=\pi J=0$, soient $x\in J$ et $\bar{m}\in M/VM$. Voyons $x$ comme élément de $J\subset W(J)$ grâce à $\gamma$. 
 Alors $x=\sum_i \,^{V^i} [a_i]$ où $a_0=x$, $a_i=0$ si $i>0$ lorsque $\gamma=0$ et  $a_i = (-1)^i \gamma^i (x)$ si $\pi J=0$. Dès lors 
$$
x m = \sum_{i\geq 0} (-1)^i V^i [a_i]F^i .m \in J.M 
$$ 
Or, lorsque $\pi=0$, $F^i : M/VM \ldrt M/VM$ s'identifie à $\a^i$. On en déduit facilement le résultat.
\qed

\begin{rema}
Le lien entre le logarithme de la section précédente et l'exponentielle de la proposition précédente
est que l'exponentiation des morphismes est obtenues par exponentiation des courbes $p$-typiques via
$\log_H: H (YJ[[Y]])\iso \Lie H \otimes YJ[[Y]]= C ( \Lie H \otimes J ) $.  
\end{rema}

Voici la proposition que nous utiliserons.

\begin{prop}\label{exputil}
Soit $T$ un $\O$-schéma, $J\subset \O_T$ un idéal cohérent muni de $\O$-P.D. nilpotentes. 
Soit $G$ un schéma en $\O$-modules plat sur $T$.
Soit $H$ un faisceau en $\O$-modules sur
$T_{fppf}$ tel que $\widehat{H}$ soit un $\O$-module formel. Il y a alors une application d'exponentiation 
$$
\exp : \Hom (G, J.\Lie \, H) \ldrt \ker \left ( \Hom (G,H) \drt \Hom (G \text{ mod } J, H \text{ mod } J) \right )
$$
telle que si $J^2=0$ et $\gamma =0$, via l'identification canonique de $\ker \left ( \Hom (G,H) \drt \Hom (G \text{ mod } J, H \text{ mod } J)\right )$ avec
$\Hom (G, J.\Lie\, H)$, $\exp =Id$. Si $J^2=\pi J =0$ et $\Pi = \a\otimes \gamma$ est l'endomorphisme nilpotent $\text{Frob}_q$-linéaire de $J.\Lie \, H = \Lie H/\pi \Lie H \otimes J$ alors
$$
\exp f = \left ( \sum_{n\geq 0} (-1)^n \Pi^n \right ) \circ f
$$
\end{prop}
\dem

Le schéma $G$ étant plat sur $T$ il suffit de définir $(\exp \, f)_U : G(U) \ldrt H(U)$ pour tout schéma $U$ affine et plat sur $T$ (i.e. définir $\exp \, f$ comme morphisme de  faisceaux sur le petit
site plat de $T$). Soit donc $U\ldrt T$ plat. Alors, d'après le lemme \ref{extPD}, les $\O$-P.D. $\gamma$ s'étendent à $J\O_U$, d'où un isomorphisme (proposition \ref{plog})
$$
\log_G : \widehat{H} (J\O_U) \iso \Lie\, H\otimes_{\O_T} J\O_U = J\Lie H\otimes_{\O_T} \O_U
$$
Posons
$$
(\exp \, f)_U: G(U) \ldrt \GG (U, J\Lie\, H) \xrig{\log^{-1}} \widehat{H} ( \GG ( U, J\O_U)) \subset H(U) 
$$
On vérifie aussitôt que cette définition est fonctorielle en $U$. 
Les assertions concernant les cas $J^2=0$ ou $J^2=\pi J=0$ se déduisent de la proposition \ref{logcar}.
\qed

\subsection{Extension du cristal de Messing généralisé aux $\O$-puissances divisées}
\label{sPiu}

Il y a un morphisme de sites 
$$
\Pi_\O :  \text{NCRIS}_\O (S/\Sigma) \ldrt \text{NCRIS} (S/\Sigma)
$$
puisque les puissances divisées classiques induisent des $\O$-puissances divisées.

\begin{theo}
Soit $H$ un $\O$-module $\pi$-divisible sur $S$. Soit $\mathcal{E}$ le cristal algèbre de Lie de la $\O$-extension vectorielle universelle sur $\text{NCRIS} (S/\Sigma)$ défini dans la proposition  \ref{exOver}. Ce cristal s'étend naturellement à $\text{NCRIS}_\O (S/\Sigma)$ au sens où il existe un
cristal en $\O_{S/\Sigma}$-modules localement libres $\F$ sur $\text{NCRIS}_\O (S/\Sigma)$ tel que 
$\mathcal{E}= \Pi_{\O *} \F$ (et donc nécessairement $\F= \Pi_\O^*\mathcal{E}$). De plus, si 
$(U\hookrightarrow T)\in \text{NCRIS}_\O (S/\Sigma)$ et $\widetilde{H}$ est un
relèvement de $H\times_S U$ à $T$ alors $\mathcal{F}_{(U\hookrightarrow T)}$ s'identifie à 
$\Lie \, E_\O (\widetilde{H})$. De plus $G\mapsto \F$ est fonctoriel en $G$.  
\end{theo}

\begin{rema}
En fait, on étend directement le cristal $\O$-extension vectorielle universelle à $\text{NCRIS}_\O (S/\Sigma)$. 
\end{rema}
\dem
La démonstration utilise l'application exponentielle construite dans la proposition \ref{exputil}
 et suit celle de \cite{Messing1}.
Faisons tout de même remarquer au lecteur qu'elle n'est pas indépendante de \cite{Messing1} puisque la proposition clef \ref{aolgh}
  qui permet d'affirmer que $E_\O (H)$ existe et est extension de $H$ par $V_\O (H)$ localement libre utilise \cite{Messing1} (c'est surtout le fait que $V_\O (H)$ soit localement libre qui est crucial). 

Comme dans \cite{Messing1} tout repose sur le théorème suivant analogue du théorème 2.2 page 129 de \cite{Messing1}

\begin{theo}
Soit $(S\hookrightarrow T) \in \text{NCRIS}_\O (S/\Sigma)$ et $H_1,H_2$ deux $\O$-modules $\pi$-divisibles sur $T$ de réduction $H_{1/S}, H_{2/S}$. Soit $f: H_{1/S}\ldrt H_{2/S}$. Il existe alors un unique morphisme 
$$
g: E_\O (H_1) \ldrt E_\O (H_2)
$$
tel que $\forall u : V_\O (H_1)\ldrt V_\O (H_2)$ relevant $V_\O (f) : V_\O (H_{1/S})\ldrt
V_\O (H_{2/S})$, dans le diagramme 
$$
\xymatrix{
V_\O (H_1) \ar@{^(->}[r]^i \ar[d]^u & E_\O (H_1) \ar[d]^g \\
V_\O (H_2) \ar@{^(->}[r]^i & E_\O (H_2)
}
$$ 
$$
g\circ i - i\circ u \in \exp \left ( \Hom ( V_\O (H_1), J.\Lie\, E_\O (H_2))\right )
$$
où $J$ est l'idéal de $S$ dans $T$ et $\exp$ est l'application définie dans la proposition \ref{exputil}.

\end{theo}

 La démonstration de ce théorème repose sur l'analogue du lemme 2.6.3 page 135 de \cite{Messing1} :

\begin{lemm}\label{lepr}
Soit $H$ un $\O$-module $\pi$-divisible sur $T$ et $N>>0$ tel que $\pi^N\O_T=0$ et $\omega_{H [\pi^N]^D}$ soit localement libre. Soit $I=\ker ( \O_T\otimes_{\Zp} \O \twoheadrightarrow \O_T)$ et 
$$
\a : H[\pi^N] \ldrt \omega_{H[\pi^N]^D}/I.\Lie\, E (H) = V_\O (H)
$$
l'application universelle telle que 
$$
\a_* \left ( 0\drt H [\pi^N] \drt H \xrig{\pi^N}  H \drt 0 \right ) = \left ( 0\drt V_\O ( H) \drt E_\O (H) \drt H \drt 0 \right )
$$
 Soit $\F$ un faisceau en $\O$-modules sur $T_{fppf}$ tel que $\widehat{\F}$ soit un $\O$-module formel. Soit un diagramme 
$\O$-équivariant
$$
\xymatrix{
H[\pi^N] \ar[d]^{\a} \ar[rd]^v \\
V_\O (H) \ar[r]^w & \mathcal{F}
}
$$
qui commute après réduction sur $S$ et tel que si
$$
v_* \left ( 0\drt H [\pi^N] \drt H \xrig{\pi^N}  H \drt 0 \right ) =  \left (  \F \drt \mathcal{G} \drt H \drt 0 \right )
$$
alors l'action de $\O$ sur $\Lie\, \mathcal{G}$ est l'action naturelle donnée par $\O\drt \O_S$. 
Il existe alors un unique morphisme $\O$-équivariant 
 $ w': V_\O (H) \ldrt \F$ tel qu'en remplaçant $w$ par $w'$ le diagramme précédent commute, $w'\equiv w $ sur $S$ et $w'-w$ soit une exponentielle.
\end{lemm}

\dem la démonstration est identique à celle de \cite{Messing1} :
 on dévisse au cas $J^2=\pi J=0$ puis on utilise la formule explicite donnée dans la 
proposition \ref{exputil} pour $\exp$ dans ce cas là. Ce qu'il faut vérifier c'est que les dévisages n'affectent pas l'hypothèse de l'action
de $\O$ sur l'objet noté $\Lie\, \mathcal{G}$ dans l'énoncé. Mais cela ne pose pas de problème puisque la seule modification faite au morphisme 
$v$ est $v\mapsto v- w\circ \a$, or $w\circ \a$ est tel que si $\xi$ est l'extension $ 0\drt H [\pi^N] \drt H \xrig{\pi^N}  H \drt 0$ alors
$(w\circ \a)_* \xi= w_* (\a_* \xi)$ mais $\O$ agit naturellement via $\O\drt \O_S$ sur $\Lie\, (\a_* \xi)$ et on en déduit le résultat d'après le 
lemme \ref{lecla} qui suit.

 La fin de l'argument page 143 de \cite{Messing1} se modifie en utilisant que sur $S$, $\a: H[\pi^N]_{/S} \ldrt V_\O (H_{/S})$ est universel pour les $\O$-morphismes $\beta$ de $H[\pi^N]_{/S}$ vers un $\O_S$-module cohérent $\mathcal{M}$ tel que si $\xi$ est l'extension précédente alors 
 $\O$ agisse naturellement sur $\Lie\, \beta_*\xi$. Il faut également utiliser de nouveau le lemme \ref{lecla} appliqué à $(\sum (-1)^n \Pi^n)^{-1}$. 

\begin{lemm}\label{lecla}
Les classes des extensions  $0\drt E_1 \drt E_2 \drt E_3\drt 0$ de $\O$-modules sur $S_{fppf}$ telles que 
$\O$ agisses naturellement sur $\Lie \, E_2$ via $\O\drt \O_S$ forment un sous $\O$-module de $\text{Ext}_\O (E_1, E_2)$. Ces classes sont stables par image directe
via un morphisme $\O$-linéaire $u: E_1 \drt {E'}_1$  tel que $\O$ agisse naturellement sur $\Lie\, E'_1$ et 
induisant $u_* : \text{Ext}_\O (E_1, E_2) \drt \text{Ext}_\O (E'_1, E_2)$.
\end{lemm}
\dem
Elle est facile et laissée au lecteur.
\qed
\\

Expliquons maintenant quels sont les arguments à adapter dans la démonstration du théorème. La principale chose à
vérifier est que l'on peut appliquer le lemme \ref{lepr} dans la démonstration de l'équivalent du lemme 2.7.4 de \cite{Messing1}. Plus précisément,
avec les notations de ce lemme dans \cite{Messing1}, il faut vérifier que si $v': G\ldrt E_\O (H)$ est l'unique relèvement
de $v'_0$ alors, si $$\xi= \left ( 0\drt G[\pi^N] \drt G \xrig{\pi^N} G \drt 0 \right )$$  $\O$ agit naturellement sur 
$\Lie \, (v'_{|G[\pi^N]})_*\xi$.
Mais $(v'_{|G[\pi^N]})_*\xi$ est une extension scindée. Son algèbre de Lie est donc une extension scindée de 
$\Lie\, E_\O (H)$ par $\Lie \, G$, d'où l'assertion.

Le reste de la démonstration en bas de la page 145 jusqu'à la page 146 s'adapte immédiatement en 
remplaçant les $\Hom$ et $Ext$ par des $\Hom_\O$ et $\Ext_\O$.
\qed

\subsection{Théorie de la déformation des $\O$-modules $\pi$-divisibles}

\begin{theo}
Soit $(S\hookrightarrow T)\in \text{NCRIS}_\O (S/\Sigma)$. Soit $\mathcal{C}$ la catégorie
des $\O$-modules $\pi$-divisibles sur $T$. Soit $\mathcal{D}$ la catégorie des $\O$-modules $\pi$-divisibles sur $S$ munis d'une filtration localement facteur directe de leur cristal généralisé
évalué sur $S\hookrightarrow T$ (cf. théorème précédent) déformant la partie vectorielle de leur $\O$-extension vectorielle universelle. Soit $F:\mathcal{C} \drt \mathcal{D}$ le foncteur qui à $H$ sur $T$ associe 
$\left ( H\times_T S,V_\O (H) \hookrightarrow \mathcal{F}_{(S\hookrightarrow T)}\right )$. Alors,
$F$ induit une équivalence de catégories.
\end{theo}
\dem
La démonstration est absolument identique à celle du théorème 1.6 page 151 du chapitre V de \cite{Messing1}. 
\qed

\subsection{Théorie de Dieudonné ``classique'' des $\O$-modules $\pi$-divisibles}\label{okklu}

Soit $k$ un corps parfait extension du corps résiduel de $\O$ et $H$ un $\O$-module $\pi$-divisible sur $k$. Soit $\s$ le Frobenius de $W(k)$. Soit $\mathbb{D} (H)$ le module de Dieudonné de $H$, un $W(k)$-module libre de rang $\text{ht} \,H$ muni de $\ph: \mathbb{D} (H) \ldrt \mathbb{D} (H)$ une application $\s$-linéaire induisant un isomorphisme après inversion de $p$. 

L'action de $\O$ sur $\mathbb{D} (H)$ muni celui-ci d'une structure de $W(k)\otimes_{\Zp} \O$-module. Donc, 
$$
\mathbb{D} (H)= \bigoplus_{\tau : F^0 \hookrightarrow W(K)_\Q} \mathbb{D} (H)_\tau
$$
où $ \mathbb{D} (H)_\tau$ est un $\O_F\otimes_{\O_{F^0},\tau} W(k)$-module libre. De plus,
$$
\ph : \mathbb{D} (H)_\tau \ldrt \mathbb{D} (H)_{\s\tau}
$$
Rappelons que $k$ est une extension du corps résiduel de $\O$. Cela définit canoniquement
un $\tau_0: F^0 \hookrightarrow W(k)_\Q$ qui définit lui-même un isomorphisme
$$
\O_F\otimes_{\O_{F^0},\tau} W(k) \iso W_\O (k)
$$
Soit $\s_\O$ le Frobenius de $W_\O(k)$.

\begin{defi}
On note $\mathbb{D}_\O (H)= \mathbb{D} (H)_{\tau_0}$, un $W_\O (k)$-module libre de rang 
$\text{ht}_\O \,H = \frac{\text{ht}\, H}{[F:\Qp]}$. Il est muni d'un opérateur $\s_\O$-linéaire 
$\ph = \frac{\pi}{p^{[F^0 :\Qp]}} \ph^{[F^0 : \Qp]}$.
\end{defi}

\begin{prop}
Le foncteur $H\longmapsto  ( \mathbb{D}_\O (H),\ph_\O)$ induit une équivalence entre la catégorie 
des $\O$-modules $\pi$-divisibles sur $k$ et celle des couples $(M,\ph)$ où $M$ est un $W_\O (k)$-module
libre de rang fini et $\ph:M\drt M$ est une application $\s_\O$-linéaire telle que $\pi M \subset \ph (M)$.
\end{prop}
\dem
Il suffit de voir que $(D,\ph)\mapsto (D_{\tau_0}, \frac{\pi}{p^{[F^0 :\Qp]}} \ph^{[F^0 : \Qp]})$ induit une équivalence entre les couples sur $W(k)$ $\;(D,\ph)$ tels que $p D \subset \ph (D) \subset D$ et si
$V=p\ph^{-1}$ alors $\O$ agit naturellement sur $D/VD$ et les $(M,\ph)$ comme dans l'énoncé.
Mais si $D=\bigoplus_{\tau } D_\tau$, 
$$
V: D_\tau \ldrt D_{\s^{-1} \tau} \;\text{ et } \; D/VD = \bigoplus_{\tau} D_\tau/V D_{\s \tau}
$$
et donc, 
$$
\forall \tau\neq \tau_0\;\; V: D_{\s \tau } \iso D_\tau
$$
est un isomorphisme tandis que 
$$
D/VD= D_{\tau_0} / V D_{\s \tau_0}
$$
On en déduit que $V D_{\s \tau_0}= V^{[F^0 : \Qp]} D_{\tau_0}$ et en posant $V_\O= V^{[F^0 : \Qp]}$ on en déduit facilement que $(D,V)\mapsto (D_{\tau_0}, V_\O)$ est une équivalence avec les
couples tels que $\pi D \subset V_\O D$. On a alors $\ph_\O V_\O = \pi$.
\qed

\begin{prop}
Supposons $p\neq 2$. 
Soit $\E$ le cristal de $H$ algèbre de Lie de la $\O$-extension vectorielle universelle sur
$\text{NCRIS}_\O (\spec (k)/\spec (\O))$. Le morphisme de Frobenius $F:H\ldrt H^{(q)}$ commute à
l'action de $\O$ et définit donc $F:\text{Fr}_q^* \E \ldrt \E$. Alors, $\mathbb{D}_\O (H)$ s'identifie à 
$$
\GG \left ( \spec (k)\hookrightarrow \spec (\O),\E \right )
$$
et $\ph$ au morphisme induit par $F$.
\end{prop}
\dem
Soit $\E'$ le cristal de $H$ algèbre de Lie de l'extension vectorielle universelle sur $\text{NCRIS} ( \spec (k)/ \spec (W(k)))$.
 Alors,
$$
\mathbb{D} (H) \simeq \GG \left ( \spec (k) \hookrightarrow \spec (W(k)),\E' \right )
$$
et par la propriété de cristal
$$
\mathbb{D} (H)\otimes_{W(k)} \O \simeq \GG \left ( \spec (\O/p\O)\hookrightarrow \spec (\O),\E' \right )
$$
et donc
\begin{eqnarray*}
\mathbb{D} (H)_{\tau_0} &= &\left ( \mathbb{D} (H)\otimes_{W(k)} \O \right )_{\O\otimes_{\Zp} \O} \O \\
& \simeq & \GG \left ( \spec (\O/p\O) \hookrightarrow \spec (O), \Pi_{\O*} \E \right ) \\
&=& \GG \left ( \spec   (\O/p\O) \hookrightarrow \spec (O),\E \right ) \\
&=& \GG \left ( \spec   (k) \hookrightarrow \spec (O),\E \right )
\end{eqnarray*}
\qed

\bibliographystyle{plain}
\bibliography{biblio}
\end{document}

%% file: DrLT1_immeuble2.pstex_t
\begin{picture}(0,0)%
\includegraphics{DrLT1_immeuble2.pstex}%
\end{picture}%
\setlength{\unitlength}{3355sp}%
\begingroup\makeatletter\ifx\SetFigFont\undefined%
\gdef\SetFigFont#1#2#3#4#5{%
  \reset@font\fontsize{#1}{#2pt}%
  \fontfamily{#3}\fontseries{#4}\fontshape{#5}%
  \selectfont}%
\fi\endgroup%
\begin{picture}(3660,3324)(4876,-6235)
\put(7813,-3748){\makebox(0,0)[lb]{\smash{{\SetFigFont{8}{9.6}{\rmdefault}{\mddefault}{\updefault}{\color[rgb]{0,0,0}boules de rayon $\frac{1}{2}$}%
}}}}
\put(7813,-3973){\makebox(0,0)[lb]{\smash{{\SetFigFont{8}{9.6}{\rmdefault}{\mddefault}{\updefault}{\color[rgb]{0,0,0}dans l'arbre}%
}}}}
\put(6763,-4048){\makebox(0,0)[lb]{\smash{{\SetFigFont{8}{9.6}{\rmdefault}{\mddefault}{\updefault}{\color[rgb]{0,0,0}$\frac{1}{2}$}%
}}}}
\end{picture}%

%% file: Newton2.pstex_t
\begin{picture}(0,0)%
\includegraphics{Newton2.pstex}%
\end{picture}%
\setlength{\unitlength}{3947sp}%
\begingroup\makeatletter\ifx\SetFigFont\undefined%
\gdef\SetFigFont#1#2#3#4#5{%
  \reset@font\fontsize{#1}{#2pt}%
  \fontfamily{#3}\fontseries{#4}\fontshape{#5}%
  \selectfont}%
\fi\endgroup%
\begin{picture}(2637,2117)(4576,-4866)
\put(4576,-3961){\makebox(0,0)[lb]{\smash{{\SetFigFont{10}{12.0}{\rmdefault}{\mddefault}{\updefault}{\color[rgb]{0,0,0}$\lambda_1$}%
}}}}
\put(5626,-4821){\makebox(0,0)[lb]{\smash{{\SetFigFont{10}{12.0}{\rmdefault}{\mddefault}{\updefault}{\color[rgb]{0,0,0}$q^i$}%
}}}}
\put(6811,-4821){\makebox(0,0)[lb]{\smash{{\SetFigFont{10}{12.0}{\rmdefault}{\mddefault}{\updefault}{\color[rgb]{0,0,0}$q^n$}%
}}}}
\put(5696,-3901){\makebox(0,0)[lb]{\smash{{\SetFigFont{10}{12.0}{\rmdefault}{\mddefault}{\updefault}{\color[rgb]{0,0,0}$1-\frac{i}{n}$}%
}}}}
\end{picture}%

%% file: Cellulaire.bbl
\begin{thebibliography}{10}

\bibitem{Berk1}
V.G. Berkovich.
\newblock Étale cohomology for non-archimedean analytic spaces.
\newblock {\em Inst. Hautes Études Sci. Publ. Math.}, 78:5--161, 1993.

\bibitem{BBM}
P.~Berthelot, L.~Breen, and W.~Messing.
\newblock {\em Théorie de Dieudonné cristalline. II}, volume 930.
\newblock Springer-Verlag, Berlin, Lecture Notes in Mathematics, 1982.

\bibitem{BerthelotOgus}
P.~Berthelot and A.~Ogus.
\newblock {F}-isocrystals and de {R}ham cohomology. {I}.
\newblock {\em Invent. Math.}, 72(2):159--199, 1983.

\bibitem{BLR4}
S.~Bosch, W.~Lütkebohmert, and M.~Raynaud.
\newblock Formal and rigid geometry. {IV}. the reduced fibre theorem.
\newblock {\em Invent. Math.}, 119(2):361--398.

\bibitem{DeJong1}
A.~J. de~Jong.
\newblock Crystalline dieudonné module theory via formal and rigid geometry.
\newblock {\em Inst. Hautes Études Sci. Publ. Math.}, 82:5--96, 1995.

\bibitem{Drinfeld1}
V.~G. Drinfel'd.
\newblock Elliptic modules.
\newblock {\em Mat. Sb. (N.S.)}, 136(94):594--627, 1974.

\bibitem{DrinfeldOmega}
V.G. Drinfeld.
\newblock Coverings of $p$-adic symmetric domains.
\newblock {\em Functional Analysis and its Applications}, 10(2):29--40, 1976.

\bibitem{Faltings6}
G.~Faltings.
\newblock Almost étale extensions.
\newblock In {\em Cohomologies $p$-adiques et applications arithmétiques, II.},
  volume 279 of {\em Astérisque}, pages 185--270, 2002.

\bibitem{Faltings7}
G.~Faltings.
\newblock Group schemes with strict $\mathcal{O}$-action.
\newblock {\em Mosc. Math. J.}, 2(2):249--279, 2002.

\bibitem{Faltings8}
G.~Faltings.
\newblock A relation between two moduli spaces studied by {V}. {G}. {D}rinfeld.
\newblock In {\em Algebraic number theory and algebraic geometry}, volume 300
  of {\em Contemp. Math.}, pages 115--129, 2002.

\bibitem{Rami}
L.~Fargues.
\newblock Application de {H}odge-{T}ate duale d'un groupe de {L}ubin-{T}ate,
  immeuble de {B}ruhat-{T}its du groupe linéaire et filtrations de
  ramification.
\newblock {\em preprint}.

\bibitem{iso5}
L.~Fargues.
\newblock L'isomorphisme entre les tours de lubin-tate et de drinfeld :
  comparaison de la cohomologie des deux tours.
\newblock {\em preprint}.

\bibitem{iso4}
L.~Fargues.
\newblock L'isomorphisme entre les tours de {L}ubin-{T}ate et de {D}rinfeld :
  démonstration du résultat principal.
\newblock {\em preprint}.

\bibitem{Points}
L.~Fargues.
\newblock L'isomorphisme entre les tours de {L}ubin-{T}ate et de {D}rinfeld au
  niveau des points.
\newblock {\em preprint}.

\bibitem{Laurent1}
L.~Fargues.
\newblock Cohomologie des espaces de modules de groupes p-divisibles et
  correspondances de langlands locales.
\newblock In {\em Variétés de Shimura, espaces de Rapoport-Zink et
  Correspondances de Langlands locales, Asterisque 291}, 2004.

\bibitem{Har5}
M.~Harris.
\newblock The local {L}anglands correspondence: Notes of (half) a course at the
  {IHP}, {S}pring 2000.
\newblock {\em A paraître dans Asterisque}.

\bibitem{Har4}
M.~Harris, R.~Taylor.
\newblock {\em The geometry and cohomology of some simple {S}himura varieties},
  volume 151 of {\em Annals of Mathematics Studies}.
\newblock Princeton University Press, Princeton, NJ, 2001.

\bibitem{HopkinsGross}
B.~H. Hopkins, M.~Gross.
\newblock Equivariant vector bundles on the {L}ubin-{T}ate moduli space.
\newblock In {\em Topology and representation theory (Evanston, IL, 1992)},
  volume 158 of {\em Contemp. Math.}, pages 23--88. Amer. Math.
  Soc.,Providence, RI, 1994.

\bibitem{Hu1}
R.~Huber.
\newblock {\em Étale cohomology of rigid analytic varieties and adic spaces.}
\newblock Aspects of Mathematics. Friedr. Vieweg and Sohn, Braunschweig, 1996.

\bibitem{Kiehl1}
R.~Kiehl.
\newblock Ausgezeichnete ringe in der nichtarchimedischen analytischen
  geometrie.
\newblock {\em J. Reine Angew. Math.}, 234:89--98, 1969.

\bibitem{Lubin1}
J.~Lubin.
\newblock Canonical subgroups of formal groups.
\newblock {\em Trans. Amer. Math. Soc.}, 251:103--127, (1979).

\bibitem{Messing1}
W.~Messing.
\newblock The crystals associated to {B}arsotti-{T}ate groups: with
  applications to abelian schemes.
\newblock {\em Lecture Notes in Mathematics}, 264, 1972.

\bibitem{Oort1}
Frans Oort.
\newblock Newton polygons and formal groups: conjectures by {M}anin and
  {G}rothendieck.
\newblock {\em Ann. of Math.}, 152(1):183--206, 2000.

\bibitem{RZ}
M.~Rapoport, Th.~Zink.
\newblock {\em Period spaces for $p$-divisible groups.}
\newblock Number 141 in Annals of Mathematics Studies. Princeton University
  Press, Princeton, NJ, 1996.

\bibitem{Ray2}
Michel Raynaud and L.~Gruson.
\newblock Critères de platitude et de projectivité. {T}echniques de
  platification d'un module.
\newblock {\em Invent. Math.}, 13:1--89, 1971.

\bibitem{Vala1}
P.~Valabrega.
\newblock On the excellent property for power series rings over polynomial
  rings.
\newblock {\em J. Math. Kyoto Univ.}, 15(2):387--395, 1975.

\bibitem{Vala2}
P.~Valabrega.
\newblock A few theorems on completion of excellent rings.
\newblock {\em Nagoya Math. J.}, 61:127--133, 1976.

\bibitem{Yu}
Jiu-Kang Yu.
\newblock On the moduli of quasi-canonical liftings.
\newblock {\em Compositio Math.}, 96(3):293--321, 1995.

\bibitem{Zink2}
Thomas Zink.
\newblock The display of a formal $p$-divisible group.
\newblock In {\em Cohomologies $p$-adiques et applications arithmétiques (I),
  Astérisque 278}, pages 127--248, 2002.

\end{thebibliography}
